\documentclass[12pt,oneside,english]{amsart}
\usepackage[T1]{fontenc}
\usepackage[latin9]{inputenc}
\usepackage{textcomp}
\usepackage{amstext}
\usepackage{amsthm}
\usepackage{amssymb}
\usepackage{graphicx}

\makeatletter
\numberwithin{equation}{section}
\numberwithin{figure}{section}
\theoremstyle{plain}
\newtheorem{thm}{\protect\theoremname}
\theoremstyle{remark}
\newtheorem{rem}[thm]{\protect\remarkname}
\theoremstyle{definition}
\newtheorem{defn}[thm]{\protect\definitionname}
\theoremstyle{definition}
\newtheorem{example}[thm]{\protect\examplename}

\usepackage{pdflscape}

\makeatother

\usepackage{babel}
\providecommand{\definitionname}{Definition}
\providecommand{\examplename}{Example}
\providecommand{\remarkname}{Remark}
\providecommand{\theoremname}{Theorem}

\begin{document}

\title{Multilevel Contours on Bundles of Complex Planes}

\maketitle
\vspace{0.3cm}

\begin{center}

\author{\textbf{Arni S.R. Srinivasa Rao}}
\begin{center}
Laboratory for Theory and Mathematical Modeling, 
\par\end{center}

\begin{center}
Medical College of Georgia, 
\par\end{center}

\begin{center}
and
\par\end{center}

\begin{center}
Department of Mathematics, 
\par\end{center}

\begin{center}
Augusta University, Georgia, USA 
\par\end{center}

\begin{center}
Email: arrao@augusta.edu (OR) arni.rao2020@gmail.com
\par\end{center}

\end{center}

\vspace{0.2cm}
\begin{abstract}
A new concept called multilevel contours is introduced through this
article by the author. Theorems on contours constructed on a bundle
of complex planes are stated and proved. Multilevel contours can transport
information from one complex plane to another. Within a random environment,
the behavior of contours and multilevel contours passing through the
bundles of complex planes are studied. Further properties of contours
by a removal process of the data are studied. The concept of 'islands'
and 'holes' within a bundle is introduced through this article. These
all constructions help to understand the dynamics of the set of points
of the bundle. Further research on the topics introduced here will
be followed up by the author. These include closed approximations
of the multilevel contour formations and their removal processes.
The ideas and results presented in this article are novel.
\end{abstract}



\vspace{1.0cm}

\dedicatory{\textbf{\large{}Dedication: This article is dedicated to my friend
and collaborator Professor Steven G. Krantz, Washington University,
St. Louis, U.S.A on completion of his 70th Birthday.}}

\vspace{0.8cm}

\textbf{\tableofcontents{}}

\pagebreak

\textbf{\listoffigures
}\pagebreak

\section{\textbf{Introduction}}

\begin{landscape}

\begin{figure}
\includegraphics[scale=0.8]{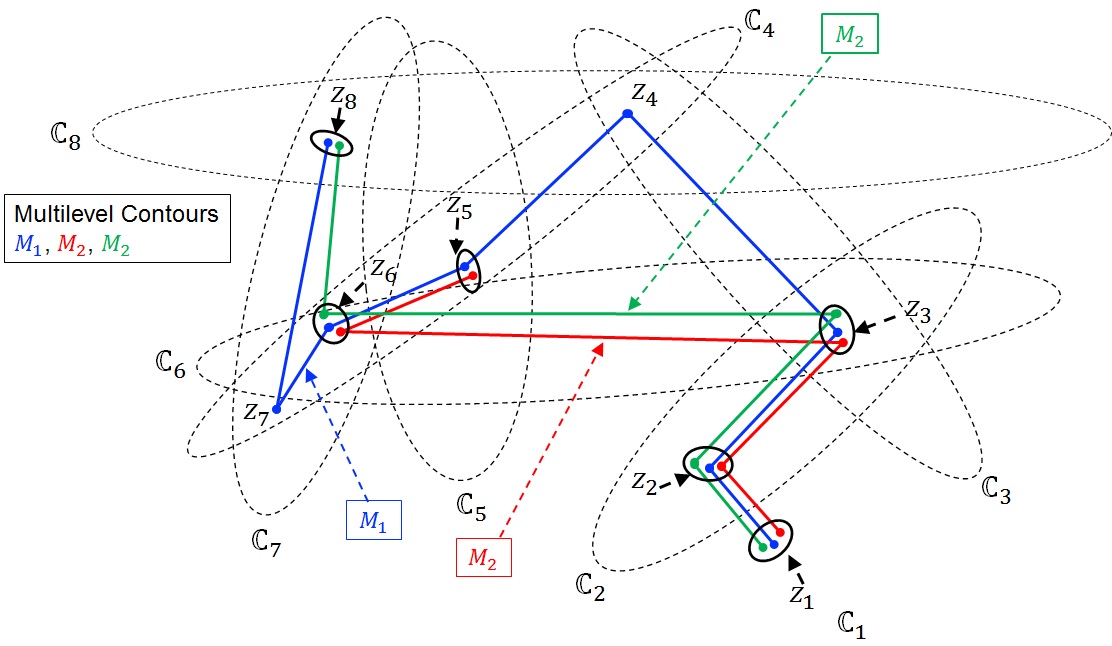}\caption{\label{fig:A-bundle-ofplanes}A bundle of complex planes and multilevel
contours $M_{1}$, $M_{2}$, $M_{3}$.}
\end{figure}

\end{landscape}

Let us consider a bundle of eight complex planes $\mathbb{C}_{1},$
$\mathbb{C}_{2}$, $\mathbb{C}_{3}$, $\mathbb{C}_{4}$, $\mathbb{C}_{5}$,
$\mathbb{C}_{6}$, $\mathbb{C}_{7}$, and $\mathbb{C}_{8}$ as shown
in Figure \ref{fig:A-bundle-ofplanes}. These planes are considered
such that one plane is parallel to any other plane in the bundle or
they could intersect with each other at some angle. Let $\gamma_{1}$
be an arc constructed from the points generated by $z(t_{1})\in\mathbb{\mathbb{C}}_{1}$
for $a_{11}\leq t_{1}\leq b_{11}$ such that $z(a_{11})=z_{1}$ and
$z(b_{11})=z_{2}.$ Here $a_{11}$, $b_{11}\mathbb{\in R}.$ Point
$z_{2}$ is located at the intersection of the planes $\mathbb{C}_{1}$
and $\mathbb{C}_{2}$. We allow constructing an arc in the plane $\mathbb{C}_{2}$
from $z_{2}$ to $z_{3}$ for $z_{3}\in\mathbb{C}_{2}$, $\mathbb{C}_{3}$,
and $\mathbb{C}_{6}$. Let $\gamma_{2}$ be an arc constructed from
$z(t_{2})\in\mathbb{\mathbb{C}}_{2}$ for $a_{12}\leq t_{2}\leq b_{12}$
($a_{12}$, $b_{12}\mathbb{\in R}.$) such that $z(a_{12})=z_{2}$
and $z(b_{12})=z_{3}.$ The arc $\gamma_{i}$ is constructed by joining
$\ensuremath{z(a_{1i})=z_{i}}$ and $\ensuremath{z(b_{1i})=z_{i+1}}$
generated by the set of points $z(t_{i})\text{ for }a_{1i}\leq t_{i}\leq b_{1i}$
for $i=3,4,...,7$ and $a_{1i}$, $b_{1i}\mathbb{\in R}.$ We allow
the possibility to construct an arc from an ending point of an arc
in a plane to a point located in a different plane if that ending
point of an arc is located at the intersection of two or more complex
planes. We saw above a few points lying at the intersection of two
or more planes. The other points, for example, $z_{4}\in\mathbb{C}_{3}$,
$\mathbb{C}_{4}$, and $\mathbb{C}_{8}$, $z_{5}\in$ $\mathbb{C}_{4}$,
and $\mathbb{C}_{8}$, $z_{6}\in\mathbb{C}_{4}$, $\mathbb{C}_{6}$,
and $\mathbb{C}_{8}$, $z_{7}\in\mathbb{C}_{4}$ and $\mathbb{C}_{8}$,
$z_{8}\in\mathbb{C}_{7}$ and $\mathbb{C}_{8}$. 

Let us form a contour by piecewise joining of arcs $\gamma_{i}$ for
$i=1,2,...,7$ and call this $M_{1}.$ Let us rename the arcs corresponding
to the contour $M_{1}$ be $\gamma_{i}^{M_{1}}$ for $i=1,2,...,7.$
Two more sample contours $M_{2}$ and $M_{3}$ are constructed using
the points $\{z_{1},z_{2},z_{3},z_{6},z_{5}\}$, and $\{z_{1},z_{2},z_{3},z_{6},z_{8}\}$,
respectively. See Figure \ref{fig:A-bundle-ofplanes}. Let the piecewise
arcs corresponding to the contour $M_{2}$ be $\gamma_{i}^{M_{2}}$
for $i=1,2,3,4,$ and the piecewise arcs corresponding to the contour
$M_{3}$ be $\gamma_{i}^{M_{3}}$ for $i=1,2,3,4.$ For the sake of
visualization, we have separated a single point at the intersecting
planes as two or more points in different colors a smaller oval-shaped
object in Figure \ref{fig:A-bundle-ofplanes}. Suppose a contour $M_{2}$
is constructed using a set of values $z(s_{i})$ for $a_{2i}\leq s_{i}\leq a_{2i}$
with corresponding arcs $\gamma_{i}^{M_{2}}$ for $i=1,2,3,4$, and
another contour $M_{3}$ is constructed $z(u_{i})$ for $a_{3i}\leq u_{i}\leq a_{3i}$
with corresponding arcs $\gamma_{i}^{M_{3}}$ for $i=1,2,3,4$. Here
$a_{2i}$, $a_{3i}$, $b_{2i}$, $b_{3i}\in\mathbb{R}.$ Let 
\[
t_{i}=\xi_{i}(\tau)\text{ \ensuremath{\left(\alpha_{1i}\leq\tau\leq\beta_{1i}\right)}}
\]
for $i=1,2,...,7$ be the parametric representation for a real-valued
function $\xi_{i}$ mapping $[\alpha_{1i},\beta_{1i}]$ onto $[a_{1i},b_{1i}].$
Then the length of the contour $M_{1}$, say, $L(M_{1})$ is computed
through the integral

\begin{equation}
L(M_{1})=\sum_{i=1}^{7}\int_{\alpha_{1i}}^{\beta_{1i}}\left|z^{'}[\xi_{i}(\tau)]\right|\xi_{i}^{'}(\tau)d\tau.\label{eq:L(M1)}
\end{equation}

Let $s_{i}=\phi_{i}(\tau)$ $(\alpha_{2i}\leq\tau\leq\beta_{2i})$
for $i=1,2,3,4$ be the parametric representation for a real-valued
function $\phi_{i}$ mapping $[\alpha_{2i},\beta_{2i}]$ onto $[a_{2i},b_{2i}]$,
and $u_{i}=\psi_{i}(\tau)$ $(\alpha_{3i}\leq\tau\leq\beta_{3i})$
for $i=1,2,3,4$ be the parametric representation for a real-valued
function $\psi_{i}$ mapping $[\alpha_{3i},\beta_{3i}]$ onto $[a_{3i},b_{3i}].$ 

Then the lengths of the contours $M_{2}$ and $M_{3}$ can be computed
as 
\begin{equation}
L(M_{2})=\sum_{i=1}^{4}\int_{\alpha_{2i}}^{\beta_{2i}}\left|z^{'}[\phi_{i}(\tau)]\right|\phi_{i}^{'}(\tau)d\tau,\label{eq:L(M2)}
\end{equation}

\begin{equation}
L(M_{3})=\sum_{i=1}^{4}\int_{\alpha_{3i}}^{\beta_{3i}}\left|z^{'}[\psi_{i}(\tau)]\right|\psi_{i}^{'}(\tau)d\tau.\label{eq:L(M3)}
\end{equation}
The contours $M_{1},$ $M_{2}$, $M_{3}$ are located on bundle of
complex planes, which we term here as multilevel contours. One can
draw several such multilevel contours on a bundle of complex planes
as shown in Figure \ref{fig:A-bundle-ofplanes}. We have considered
eight complex planes and three multilevel contours as an example,
but one can extend these examples to demonstrate the intersection
of many more complex planes and contours passing through them. Although
multilevel contours are newly introduced here in this article, the
principles associated with contours on a single complex plane can
be found in any standard textbooks, see for example \cite{Ahlfors-book1978,Churchil-Brwon-book,Krantz-geometric,Rudin-Real-Complex}.

\section{\textbf{Infinitely Many Bundles of Complex Planes}}

\begin{figure}
\includegraphics[scale=0.19]{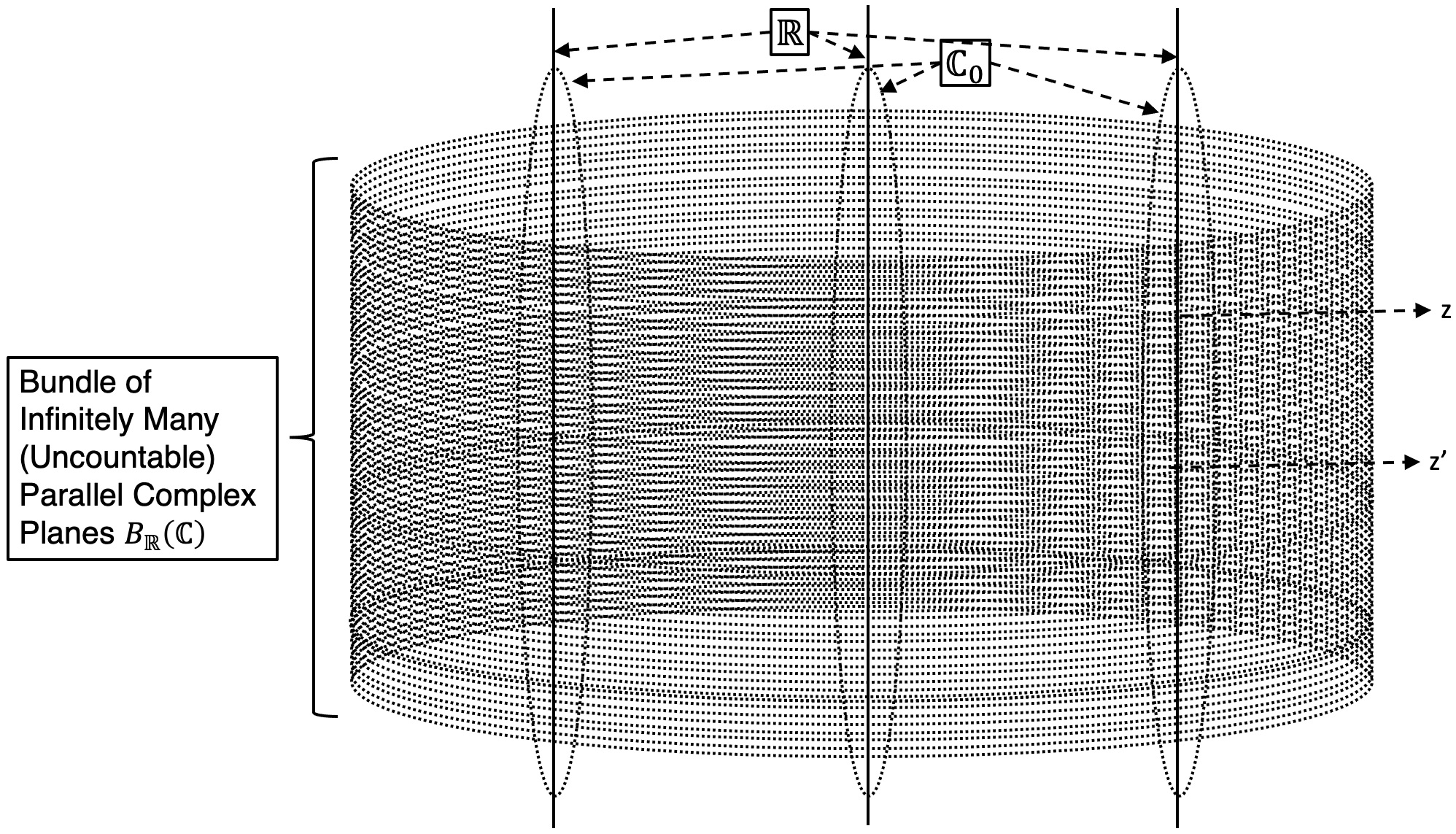}

\caption{\label{fig:Bundle-of-infinitelymanyparallelplanes}Bundle of infinitely
many (uncountable) complex planes and shortest contour passing through
them. }
\end{figure}

Let us consider infinitely many (uncountable) complex planes parallel
to each other as shown in Figure \ref{fig:Bundle-of-infinitelymanyparallelplanes}
and call this $B_{\mathbb{R}}(\mathbb{C}).$ We try to form contours
passing through these bundles and understand the behavior of the contours
at the intersection of other planes. A contour passing through the
points (complex numbers) lying in the intersection planes are given
a feature to switch a plane. The points at intersections are assumed
to possess special features also behavior of those points under a
random environment that we will see in this article. Before we understand
other properties, let us prove a Theorem.

\begin{thm}
\label{thm:The-shortest-possible}The shortest possible multilevel
contour passing through $B_{\mathbb{R}}(\mathbb{C})$ is the real
line. 
\end{thm}

\begin{proof}
Consider a complex plane that is located perpendicular to bundle $B_{\mathbb{R}}(\mathbb{C})$
such that it is at $90$ degrees with the $x-$axis. Call this $\mathbb{C}_{0}$.
Using $\mathbb{C}_{0}$ we slice $B_{\mathbb{R}}(\mathbb{C})$ vertically
at an arbitrary location as shown in Figure \ref{fig:Bundle-of-infinitelymanyparallelplanes}.
Each slice will have uncountable lines distinct from each other, and
these lines are parallel to each other. 
\end{proof}

\begin{proof}
Let $l$ and $q$ be two arbitrary lines among these uncountable lines
formed out of the above slicing method. Suppose we chose a point $z_{1}$
on $l$. There exists a point $z_{2}$ on $q$ such that the $x-$coordinate
and $y-$coordinate of both $z_{1}$ and $z_{2}$ are the same. Note
that $z_{1},z_{2}\in\mathbb{C}_{0}$. Depending upon how we visualize
the $xy-$axes of $\mathbb{C}_{0}$, the following possibilities for
the values of $z_{1}$ and $z_{2}$ will arise: 

(i) $z_{1}=(0,0)$ and $z_{2}=(0,0)$ or $z_{1}\neq(0,0)$ and $z_{2}\neq(0,0)$. 

(ii) If $z_{1},z_{2}\neq(0,0)$, then both $z_{1}$ and $z_{2}$ will
have either a non-zero $x-$coordinate (and $y-$coordinate as zero)
or a non-zero $y-$coordinate (and $x-$coordinate zero). We first
assume that both $z_{1}$ and $z_{2}$ are on the $y-$coordinate.
Let $C(z_{1},z_{2})$ be a contour described by the equation $z(t_{1})$
for $a_{1}\leq t_{1}\leq a'_{1},$ where $z(a_{1})=z_{1}$ for $z_{1}\in\mathbb{C}_{0},l$
and $z(a'_{1})=z_{2}$ for $z_{2}\in\mathbb{C}_{0},q.$ Here $C(z_{1},z_{2})$
is a multilevel contour because $z_{1}$ and $z_{2}$ are points on
parallel lines $l$ and $q$ on different planes but both these points
also belong to $\mathbb{C}_{0}.$ Suppose 

\[
t_{1}=\varepsilon_{1}(\tau)\text{ \ensuremath{\left(\delta_{1}\leq\tau\leq\delta'_{1}\right)}}
\]
be the parametric representation for $C(z_{1},z_{2})$, where $\varepsilon_{1}$
is a real-valued function mapping $[\delta_{1},\delta'_{1}]$ onto
the interval $[a_{0},a'_{0}].$ The length of the contour $C(z_{1},z_{2})$
is obtained by 

\begin{equation}
L\left[C(z_{1},z_{2})\right]=\int_{\delta_{1}}^{\delta'_{1}}\left|z^{'}[\varepsilon_{1}(\tau)]\right|\varepsilon_{1}^{'}(\tau)d\tau\label{eq:LCz1z2}
\end{equation}
Suppose we consider a point $z_{3}$ on the same plane in which the
point $z_{2}$ lies but not on the line $q$ such that $z_{3}\notin\mathbb{C}_{0}.$
That means, $z_{3}\neq z_{2}$. Let $C(z_{2},z_{3})$ be a contour
described by the equation $z(t_{2})$ for $a_{2}\leq t_{2}\leq a'_{2},$
where $z(a_{2})=z_{2}$ for $z_{2}\in\mathbb{C}_{0},q$ and $z(a'_{2})=z_{3}$
for $z_{3}\notin\mathbb{C}_{0}\text{ and }z_{3}\notin q.$ Here $C(z_{2},z_{3})$
is not a multilevel contour. Suppose 

\[
t_{2}=\varepsilon_{2}(\tau)\text{ \ensuremath{\left(\delta_{2}\leq\tau\leq\delta'_{2}\right)}}
\]
be the parametric representation for $C(z_{2},z_{3})$, where $\varepsilon_{2}$
is a real-valued function mapping $[\delta_{2},\delta'_{2}]$ onto
the interval $[a_{2},a'_{2}].$ The length of the contour $C(z_{2},z_{3})$
can be obtained by 

\begin{equation}
L\left[C(z_{2},z_{3})\right]=\int_{\delta_{2}}^{\delta'_{2}}\left|z^{'}[\varepsilon_{2}(\tau)]\right|\varepsilon_{2}^{'}(\tau)d\tau.\label{eq:LCz2z3}
\end{equation}
Since $z_{3}$ is not in $\mathbb{C}_{0}$, we cannot draw a contour
directly from $z_{1}$ to $z_{3}$. To draw a contour to $z_{3}$
from $z_{1}$, we can have piecewise arcs passing through $z_{2}$
or through any other points of the line $q$to $z_{3}.$ If the contour
from $z_{1}$ to $z_{3}$ passing through $z_{2}$, then
\begin{equation}
\begin{array}{c}
\int_{\delta_{1}}^{\delta'_{1}}\left|z^{'}[\varepsilon_{1}(\tau)]\right|\varepsilon_{1}^{'}(\tau)d\tau<\\
\int_{\delta_{1}}^{\delta'_{1}}\left|z^{'}[\varepsilon_{1}(\tau)]\right|\varepsilon_{1}^{'}(\tau)d\tau+\int_{\delta_{2}}^{\delta'_{2}}\left|z^{'}[\varepsilon_{2}(\tau)]\right|\varepsilon_{2}^{'}(\tau)d\tau,
\end{array}\label{eq:Lcz1z2<l+L}
\end{equation}

If the contour from $z_{1}$ to $z_{3}$ passing through an arbitrary
point, say, $z_{4}$ $(z_{4}\neq z_{2})$ for $z_{4}\in q$, then

\begin{equation}
\int_{\delta_{1}}^{\delta'_{1}}\left|z^{'}[\varepsilon_{1}(\tau)]\right|\varepsilon_{1}^{'}(\tau)d\tau<\int_{\delta_{3}}^{\delta'_{3}}\left|z^{'}[\varepsilon_{3}(\tau)]\right|\varepsilon_{3}^{'}(\tau)d\tau,\label{eq:Lz1z2<Lz1z4}
\end{equation}
where the R.H.S. of the inequality (\ref{eq:Lz1z2<Lz1z4}) is the
length of the contour $C(z_{1},z_{4})$ described by the equation
$z(t_{3})$ for $a_{3}\leq t_{3}\leq a'_{3},$ where $z(a_{3})=z_{1}$
for $z_{1}\in\mathbb{C}_{0},l$ and $z(a'_{3})=z_{4}$ for $z_{4}\in\mathbb{C}_{0}.$
Here $t_{3}=\varepsilon_{3}(\tau)\text{ \ensuremath{\left(\delta_{3}\leq\tau\leq\delta'_{3}\right)}}$
is the parametric representation for $C(z_{1},z_{4})$ and $\varepsilon_{3}$
is a real-valued function mapping $[\delta_{3},\delta'_{3}]$ onto
the interval $[a_{3},a'_{3}].$ From (\ref{eq:Lz1z2<Lz1z4}) and (\ref{eq}),
we have

\[
\int_{\delta_{1}}^{\delta'_{1}}\left|z^{'}[\varepsilon_{1}(\tau)]\right|\varepsilon_{1}^{'}(\tau)d\tau+\int_{\delta_{2}}^{\delta'_{2}}\left|z^{'}[\varepsilon_{2}(\tau)]\right|\varepsilon_{2}^{'}(\tau)d\tau
\]
\begin{equation}
<\int_{\delta_{3}}^{\delta'_{3}}\left|z^{'}[\varepsilon_{3}(\tau)]\right|\varepsilon_{3}^{'}(\tau)d\tau+\int_{\delta_{4}}^{\delta'_{4}}\left|z^{'}[\varepsilon_{4}(\tau)]\right|\varepsilon_{4}^{'}(\tau)d\tau,\label{eq:LL<LL}
\end{equation}
where the second term of R.H.S. of the inequality (\ref{eq:LL<LL})
is the length of the contour $C(z_{4},z_{3})$ described by the equation
$z(t_{4})$ for $a_{4}\leq t_{4}\leq a'_{4},$ where $z(a_{4})=z_{4}$
for $z_{4}\in\mathbb{C}_{0},q$ and $z(a'_{4})=z_{3}$ for $z_{3}\neq\mathbb{C}_{0}.$
Here $t_{4}=\varepsilon_{4}(\tau)\text{ \ensuremath{\left(\delta_{4}\leq\tau\leq\delta'_{4}\right)}}$
is the parametric representation for contour $C(z_{4},z_{3})$ and
$\varepsilon_{4}$ is a real-valued function mapping $[\delta_{4},\delta'_{4}]$
onto the interval $[a_{4},a'_{4}].$ From (\ref{eq}) to (\ref{eq:LL<LL})
we conclude that $L\left[C(z_{1},z_{2})\right]$ is the shortest contour
and it is a line. We can chose a point $z$ on a line $p$ for $p\in\mathbb{C}_{0}$
and find a corresponding point $z'$ on a line, say, $p'$ for $p'\in\mathbb{C}_{0}$
and construct an argument as above to see that $L\left[C(z,z')\right]$
is the shortest. The $x$ and $y$ coordinates of $z$ and $z'$ are
the same. We can construct infinitely many contours between various
points $z$, $z'$ lying on the slice such that $L\left[C(z,z')\right]$
is the shortest. If $p$ and $p'$ are the lines from adjacent planes
then 

\begin{equation}
\bigcup_{z,z'}\int_{\delta_{p}}^{\delta'_{p}}\left|z^{'}[\varepsilon_{p}(\tau)]\right|\varepsilon_{p}^{'}(\tau)d\tau=\infty\label{eq:=00003Dinfty}
\end{equation}
and 

\[
\bigcup_{z,z'}C\left(z,z'\right)\sim\mathbb{R}.
\]
In (\ref{eq:=00003Dinfty}), the integral on the LH.S. is the length
of the contour $C(z,z')$ described by the equation $z(t_{p})$ for
$a_{p}\leq t_{p}\leq a'_{p},$ where $z(a_{p})=z$ for $z\in\mathbb{C}_{0},p$
and $z(a'_{p})=z'$ for $z'\in p',\mathbb{C}_{0}.$ Here $t_{p}=\varepsilon_{p}(\tau)\text{ \ensuremath{\left(\delta_{p}\leq\tau\leq\delta'_{p}\right)}}$
is the parametric representation for contour $C(z,z')$ and $\varepsilon_{p}$
is a real-valued function mapping $[\delta_{p},\delta'_{p}]$ onto
the interval $[a_{p},a'_{p}].$
\end{proof}
\begin{rem}
\label{rem:limzn=00003Dz}Let $z_{n}$ be a sequence on the slice
$\mathbb{C}_{0}$. then
\end{rem}

\[
\lim_{n\rightarrow\infty}z_{n}=z
\]
since $z$ is equal to $z_{n}$ $\forall n.$
\begin{rem}
Remark (\ref{rem:limzn=00003Dz}) is true for each slice on the bundle
that is parallel to $\mathbb{C}.$ But the limits of convergence on
each slice are different.
\end{rem}

$ $
\begin{rem}
The distances between each pair of adjacent points on every contour
created by slicing parallel to $\mathbb{C}_{0}$ are equal. 
\end{rem}

Suppose we remove the space created by $\mathbb{C}_{0}$ from the
bundle $B_{\mathbb{R}}(\mathbb{C})$, then the bundle formed on the
left of $\mathbb{C}_{0}$ be denoted by $B_{\mathbb{R}}^{L}(\mathbb{C}_{0})$
and the bundle formed on the right of $\mathbb{C}_{0}$ be denoted
by $B_{\mathbb{R}}^{R}(\mathbb{C}_{0}).$ Then

\begin{equation}
B_{\mathbb{R}}^{L}(\mathbb{C}_{0})\cup\mathbb{C}_{0}\cup B_{\mathbb{R}}^{R}(\mathbb{C}_{0})=B_{\mathbb{R}}(\mathbb{C})\label{eq:firstsliceunion}
\end{equation}
The set $B_{\mathbb{R}}(\mathbb{C})-\mathbb{C}_{0}$ is defined as

\[
B_{\mathbb{R}}(\mathbb{C})-\mathbb{C}_{0}=\left\{ z:z\in B_{\mathbb{R}}(\mathbb{C})\text{ and }z\notin\mathbb{C}_{0}\right\} 
\]
forms a disconnected set because

\begin{equation}
B_{\mathbb{R}}(\mathbb{C})-\mathbb{C}_{0}=B_{\mathbb{R}}^{L}(\mathbb{C}_{0})\cup B_{\mathbb{R}}^{R}(\mathbb{C}_{0})\label{eq:firstdisconnected}
\end{equation}
and $B_{\mathbb{R}}^{L}(\mathbb{C}_{0})$ and $B_{\mathbb{R}}^{R}(\mathbb{C}_{0})$
are disjoint and non-empty sets within $B_{\mathbb{R}}(\mathbb{C}).$
No multilevel contours can be drawn passing through various planes
of the bundle $B_{\mathbb{R}}^{L}(\mathbb{C}_{0})$ unless $B_{\mathbb{R}}^{L}(\mathbb{C}_{0})$
is sliced similarly as we did for the bundle $B_{\mathbb{R}}(\mathbb{C})$.
Under similar circumstances, no multilevel contours can be drawn passing
through the planes of $B_{\mathbb{R}}^{R}(\mathbb{C}_{0}).$ Suppose
we slice the bundle $B_{\mathbb{R}}^{L}(\mathbb{C}_{0})$ with a complex
plane that was kept parallel to $\mathbb{C}_{0}$ and call this new
plane $\mathbb{C}_{r}$ ($r>0$) that was used for slicing. Now $\mathbb{C}_{1}$
intersects with each and every plane of $B_{\mathbb{R}}^{L}(\mathbb{C}_{0}).$
Let us remove the space created by $\mathbb{C}_{r}$ from $B_{\mathbb{R}}^{L}(\mathbb{C}_{0})$
to form a new disjoint bundles $B_{\mathbb{R}}^{L}(\mathbb{C}_{r})$
and $B_{\mathbb{R}}^{R}(\mathbb{C}_{r})$ such that

\begin{equation}
B_{\mathbb{R}}^{L}(\mathbb{C}_{r})\cup\mathbb{C}_{r}\cup B_{\mathbb{R}}^{R}(\mathbb{C}_{r})=B_{\mathbb{R}}^{L}(\mathbb{C}_{0}),\label{eq:connected-BLRC0}
\end{equation}
where $B_{\mathbb{R}}^{L}(\mathbb{C}_{r})$ is the bundle formed on
the left of $\mathbb{C}_{r}$ and $B_{\mathbb{R}}^{R}(\mathbb{C}_{r})$
is the bundle formed on the right of $\mathbb{C}_{r}$ due to removal
of the space $\mathbb{C}_{r}$ from $B_{\mathbb{R}}^{L}(\mathbb{C}_{0}).$
Using the similar argument of (\ref{eq:firstdisconnected}), we write
below $B_{\mathbb{R}}^{L}(\mathbb{C}_{0})-\mathbb{C}_{r}$ as an union
of two disconnected sets

\begin{equation}
B_{\mathbb{R}}^{L}(\mathbb{C}_{0})-\mathbb{C}_{r}=B_{\mathbb{R}}^{L}(\mathbb{C}_{r})\cup B_{\mathbb{R}}^{R}(\mathbb{C}_{r})\:\text{}(r>0)\label{eq:DISCONNECTED-BLRc0}
\end{equation}
Although we could not draw a multilevel contour passing through all
the planes of the bundle $B_{\mathbb{R}}^{L}(\mathbb{C}_{r})$, one
can draw such a contour through the plane $\mathbb{C}_{r}$ while
it intersects the bundle $B_{\mathbb{R}}^{L}(\mathbb{C}_{r}).$ One
can write another disjoint set 

\begin{equation}
B_{\mathbb{R}}^{L}(\mathbb{C}_{r})-\mathbb{C}_{r'}=B_{\mathbb{R}}^{L}(\mathbb{C}_{r'})\cup B_{\mathbb{R}}^{R}(\mathbb{C}_{r'})\:\text{}(r,r'>0),\label{eq:disconnected BLRCr}
\end{equation}
where $\mathbb{C}_{r'}$ is complex plane used to slice $B_{\mathbb{R}}^{L}(\mathbb{C}_{r'})$
and $\mathbb{C}_{r}$ is complex plane used to slice $B_{\mathbb{R}}^{L}(\mathbb{C}_{0}).$ 

Let us now consider right side of the plane $\mathbb{C}_{0}$ within
the bundle $B_{\mathbb{R}}(\mathbb{C})$ i.e. $B_{\mathbb{R}}^{R}(\mathbb{C}_{0}).$
Suppose we slice the bundle $B_{\mathbb{R}}^{R}(\mathbb{C}_{0})$
using a plane parallel to $\mathbb{C}_{0}$ and call this $\mathbb{C}_{s}$
($s>0)$. Let us remove the space created by $\mathbb{C}_{s}$ from
$B_{\mathbb{R}}^{R}(\mathbb{C}_{0})$ such that

\begin{equation}
B_{\mathbb{R}}^{R}(\mathbb{C}_{0})-\mathbb{C}_{s}=B_{\mathbb{R}}^{L}(\mathbb{C}_{s})\cup B_{\mathbb{R}}^{R}(\mathbb{C}_{r})\:\text{}(s>0).\label{eq:firstdisconnected-RIGHT}
\end{equation}
From the above constructions (\ref{eq:firstsliceunion}) through (\ref{eq:firstdisconnected-RIGHT}),
the set of points of the bundle $B_{\mathbb{R}}(\mathbb{C})$ with
intersecting planes $\mathbb{C}_{r}$, $\mathbb{C}_{0}$, $\mathbb{C}_{s}$
are written as

\begin{align}
B_{\mathbb{R}}(\mathbb{C}) & =B_{\mathbb{R}}^{L}(\mathbb{C}_{0})\cup\mathbb{C}_{0}\cup B_{\mathbb{R}}^{R}(\mathbb{C}_{0})\nonumber \\
 & =\left[B_{\mathbb{R}}^{L}(\mathbb{C}_{r})\cup\mathbb{C}_{r}\cup B_{\mathbb{R}}^{R}(\mathbb{C}_{r})\right]\cup\mathbb{C}_{0}\cup\left[B_{\mathbb{R}}^{L}(\mathbb{C}_{s})\cup\mathbb{C}_{s}\cup B_{\mathbb{R}}^{R}(\mathbb{C}_{s})\right]\label{eq:totalconnected}
\end{align}
$\implies$

\begin{align}
B_{\mathbb{R}}(\mathbb{C})-\left(\mathbb{C}_{0}\cup\mathbb{C}_{r}\cup\mathbb{C}_{r}\right) & =\nonumber \\
 & \left[B_{\mathbb{R}}^{L}(\mathbb{C}_{r})\cup B_{\mathbb{R}}^{R}(\mathbb{C}_{r})\right]\cup\left[B_{\mathbb{R}}^{L}(\mathbb{C}_{s})\cup B_{\mathbb{R}}^{R}(\mathbb{C}_{s})\right]\label{eq:eq:totaldisconnected}
\end{align}

\begin{figure}
\includegraphics[scale=0.25]{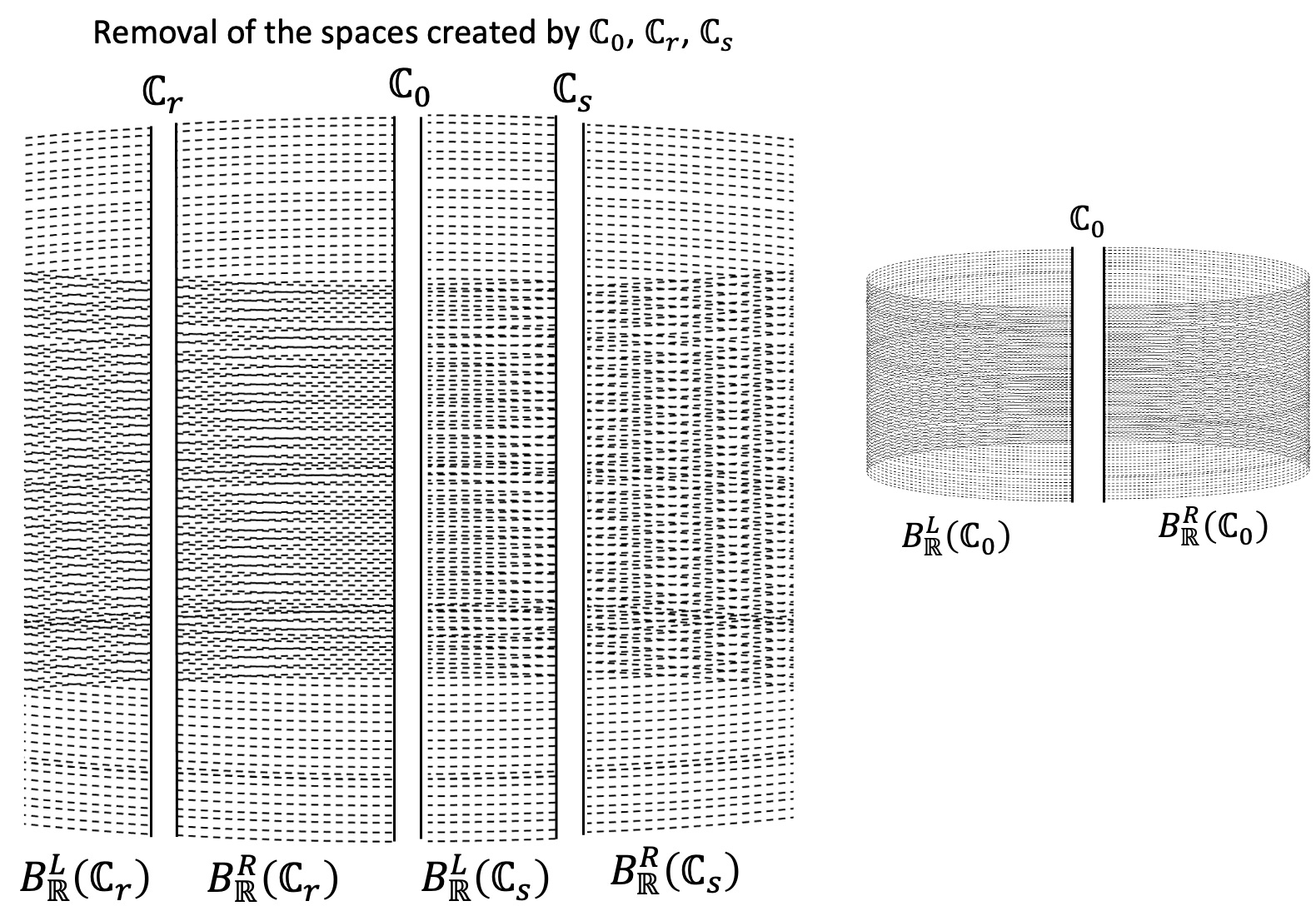}

\caption{\label{fig:Four-disconnectedsets}Creation of disconnected sets out
of the bundle $B_{\mathbb{R}}(\mathbb{C})$ due to the slicing and
removal of the spaces created by the complex planes $\mathbb{C}_{0}$,
$\mathbb{C}_{r}$, $\mathbb{C}_{s}$.}
\end{figure}

\begin{landscape}

\begin{figure}
\includegraphics[scale=0.3]{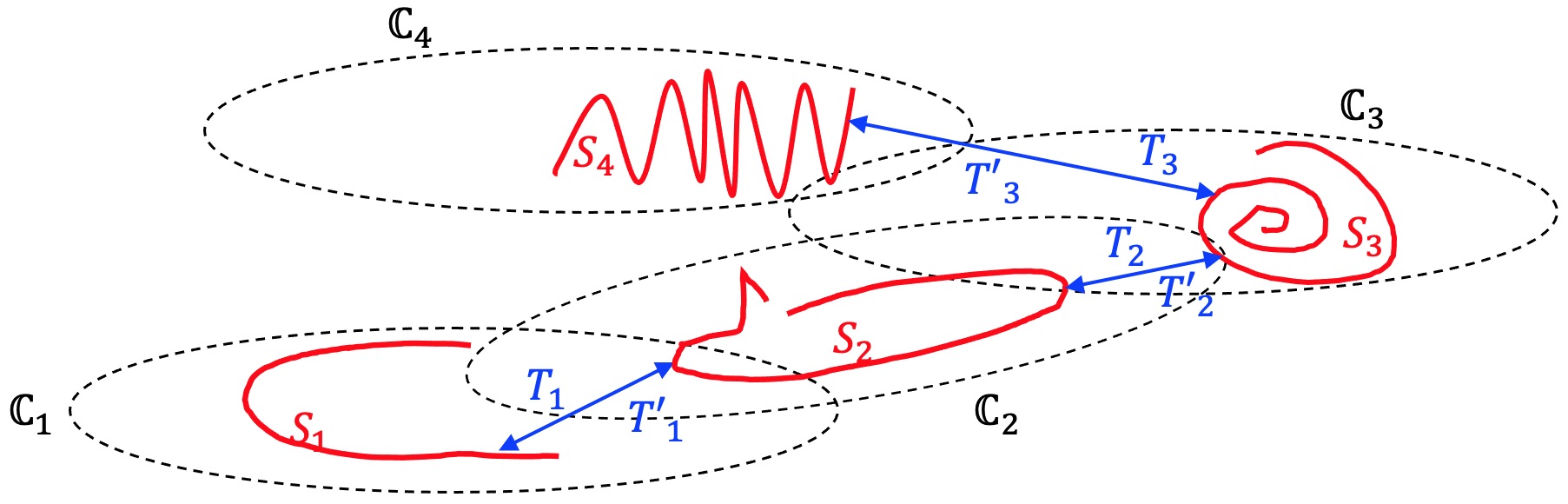}

\caption{\label{fig:Transportation-of-information}Transportation of information
through multilevel contours $S_{1}$,$...,$$S_{7}$ passing through
four intersecting complex planes $\mathbb{C}_{1}$, $\mathbb{C}_{2}$,
$\text{\ensuremath{\mathbb{C}}}_{3}$, and $\mathbb{C}_{4}$.}
\end{figure}

\end{landscape}

The four disconnected sets in (\ref{eq:eq:totaldisconnected}) can
be used for forming infinitely many disconnected sets. Due to the
removal of the spaces as shown in Figure \ref{fig:Four-disconnectedsets},
it is impossible to draw multilevel contours within and between these
four disconnected sets in (\ref{eq:eq:totaldisconnected}). One of
the advantages of multilevel contours is to develop trees of arcs,
paths that can be used for transportation of information between two
or more interacting (intersecting) complex planes. The bundle B(C)
have parallel planes however the contours drawn within each of these
planes need not be similar. One can construct a functional mapping
such that a contour or an arc drawn in one plane could be mapped to
a contour or an arc in another plane. But these two sets of contours
could be used for transporting information continuously only if these
contours in two or more planes are path-connected (Figure \ref{fig:Transportation-of-information}). 

Suppose $S_{1}$, $S_{2}$, $S_{3}$, $S_{4}$ are four contours drawn
for specific purposes in four different complex planes $\mathbb{C}_{1},$
$\mathbb{C}_{2}$, $\mathbb{C}_{3}$, $\mathbb{C}_{4}$, respectively.
Let $S_{i}$ be described by $z(t)$ for $a_{s_{i}}\leq t_{s_{1}}\leq a'_{s_{i}}$
where $a_{s_{i}},$ $a'_{s_{i}}\in\mathbb{R}$ for $i=1,2,3,4.$ See
Figure \ref{fig:Transportation-of-information}. Had these four planes
have no intersecting set of points, then one couldn't construct multilevel
contours passing through these planes. Let $z(a_{s_{i}})=z_{1}^{S_{i}}$
be the starting point and $z(a'_{s_{i}})=z_{2}^{S_{i}}$ be the ending
point of the contour $S_{i}$ for $i=1,2,3,4.$ Suppose each independent
contour $S_{i}$ has a specific information stored in it. Information
stored in $S_{1}$ is transferred to the contour $S_{2}$ using a
contour $T_{1}$. Suppose $\mathbb{C}_{1}\cap\mathbb{C}_{2}=\left\{ z:z\in\mathbb{C}_{1}\text{ and }z\mathbb{\in C}_{2}\right\} $.
Then for the structure of the planes and contours $S_{1}$ through
$S_{4}$ in Figure \ref{fig:Transportation-of-information}, we have
$\mathbb{C}_{1}\cap\mathbb{C}_{2}\neq\phi$ (empty set) and $S_{1}\cap\left(\mathbb{C}_{1}\cap\mathbb{C}_{2}\right)=\phi$.
Let $T_{1}$ be a contour drawn from a point in $S_{1}$ to a point
in $S_{2}$ through a set of points for which $\mathbb{C}_{1}\cap\mathbb{C}_{2}\neq\phi.$
$T_{1}$ is described by $z(u_{t_{1}})$ for $a_{t_{1}}\leq u_{t_{1}}\leq a'_{t_{1}}$
where $a_{t_{1}},$ $a'_{t_{1}}\in\mathbb{R}$. The starting point
of $T_{1}$ lies on $S_{1}$ and the ending point of $T_{1}$ lies
on $S_{2}$. We call $T_{1}$ a transporting contour. Using $T_{1}$
the information stored in $S_{1}$ and $S_{2}$ can be communicated.
We will discuss later more features on information transfer. The length
of the multilevel contour due to $S_{1}$, $T_{1}$, $S_{2}$ , say
$L[S_{1},S_{2}]$ is computed as 
\begin{align}
L[S_{1},S_{2}] & =\int_{\delta_{s_{1}}}^{\delta'_{s_{1}}}\left|z[\varepsilon_{s_{1}}(\tau)]\right|\varepsilon_{s_{1}}^{'}(\tau)d\tau+\int_{\omega_{t_{1}}}^{\omega_{t_{1}}'}\left|z[\eta_{t_{1}}(\tau)]\right|\eta_{t_{1}}^{'}(\tau)d\tau\nonumber \\
 & +\int_{\delta_{s_{2}}}^{\delta'_{s_{2}}}\left|z[\varepsilon_{s_{2}}(\tau)]\right|\varepsilon_{s_{2}}^{'}(\tau)d\tau\label{eq:LS1S2}
\end{align}
Here $t_{s_{i}}=\varepsilon_{s_{i}}(\tau)\text{ \ensuremath{\left(\delta_{s_{i}}\leq\tau\leq\delta'_{s_{i}}\right)}}$
for $i=1,2,3,4$ is the parametric representation for contour $S_{i}$
with a real-valued function $\varepsilon_{s_{i}}$ mapping $[\delta_{s_{i}},\delta'_{s_{i}}]$
onto the interval $[a_{s_{i}},a'_{s_{i}}],$ and $u_{t_{1}}=\eta_{t_{1}}(\tau)\text{ \ensuremath{\left(\omega_{t_{1}}\leq\tau\leq\omega'_{t_{1}}\right)}}$
is the parametric representation for contour $T_{1}$ with a real-valued
function $\eta_{t_{1}}$ mapping $[\omega_{t_{1}},\omega'_{t_{1}}]$
onto the interval $[a_{t_{1}},a'_{t_{1}}]$. Since the total information
stored in $S_{1}$ and $S_{2}$ are exchanged we have considered total
lengths of $S_{1}$ and $S_{2}$ even though $T_{1}$ could be connected
with any point of $S_{1}$ and $S_{2}$. Since $z(a_{t_{1}})\in S_{1}$
and $z(a_{t_{2}})\in S_{2}$ and the $T_{1}$ describes a contour
from $z(a_{t_{1}})$ to $z(a_{t_{2}})$, the length of the middle
integral in R.H.S. of (\ref{eq:LS1S2}) is not constant. The transportation
contour can be used to transport information from $S_{2}$ to $S_{1}$.
Let us denote this by the contour $T'_{1}$. In that case, the starting
point of $T'_{1}$ lies on $S_{2}$, and the ending point of $T'_{1}$
lies on $S_{1}$. $T'_{1}$ is described by $z'(u'_{t_{1}})$ for
$a_{t_{1}}\leq u_{t_{1}}\leq a'_{t_{1}}$ where $a_{t_{1}},$ $a'_{t_{1}}\in\mathbb{R}$,
and $u'_{t_{1}}=\eta*_{t_{1}}(\tau)\text{ \ensuremath{\left(\omega_{t_{1}}\leq\tau\leq\omega'_{t_{1}}\right)}}$
is the parametric representation for contour $T'_{1}$ with a real-valued
function $\eta*_{t_{1}}$ mapping $[\omega_{t_{1}},\omega'_{t_{1}}]$
onto the interval $[a_{t_{1}},a'_{t_{1}}]$. When we measure the length
$S_{2}$ to $S_{1}$, sat $L(S_{2},S_{1}]$, the orientation of the
transportation contour changes, and it is computed as
\begin{align}
L[S_{2},S_{1}] & =\int_{\delta_{s_{1}}}^{\delta'_{s_{1}}}\left|z[\varepsilon_{s_{1}}(\tau)]\right|\varepsilon_{s_{1}}^{'}(\tau)d\tau+\int_{\omega'_{t_{1}}}^{\omega_{t_{1}}}\left|z[\eta*_{t_{1}}(\tau)]\right|\eta*_{t_{1}}^{'}(\tau)d\tau\nonumber \\
 & +\int_{\delta_{s_{2}}}^{\delta'_{s_{2}}}\left|z[\varepsilon_{s_{2}}(\tau)]\right|\varepsilon_{s_{2}}^{'}(\tau)d\tau\label{eq:LS2S1}
\end{align}

The values of the middle integrals in the R.H.S. of (\ref{eq:LS1S2})
and (\ref{eq:LS2S1}) need not be the same unless the below condition
is satisfied:

\begin{equation}
z(a_{t_{1}}),\text{ }z'(a'_{t_{1}})\in S_{1}\text{ and }z'(a_{t_{1}})\text{ }z(a'_{t_{1}})\in S_{2}.\label{eq:twoconditionsofMIDDLEterm}
\end{equation}

The transportation contour $T_{2}$ joins $S_{2}$ to $S_{3}$ and
the transportation contour $T'_{2}$ joins $S_{3}$ to $S_{2}.$ The
contour $T_{2}$ is described by $z(u_{t_{2}})$ for $a_{t_{2}}\leq u_{t_{2}}\leq a'_{t_{2}}$
where $a_{t_{2}},$ $a'_{t_{2}}\in\mathbb{R}$. The starting point
of $T_{2}$ is in the set $S_{2}$ and the ending point of $T_{2}$
is in the set $S_{3}$. The function $u_{t_{2}}=\eta_{t_{2}}(\tau)\text{ \ensuremath{\left(\omega_{t_{2}}\leq\tau\leq\omega'_{t_{2}}\right)}}$
is the parametric representation for contour $T_{2}$ with a real-valued
function $\eta_{t_{2}}$ mapping $[\omega_{t_{2}},\omega'_{t_{2}}]$
onto the interval $[a_{t_{2}},a'_{t_{2}}]$. The contour $T'_{2}$
is described by $z(u'_{t_{2}})$ for $a_{t_{2}}\leq u'_{t_{2}}\leq a'_{t_{2}}$
where $a_{t_{2}},$ $a'_{t_{2}}\in\mathbb{R}$, and the function $u'_{t_{2}}=\eta*_{t_{2}}(\tau)\text{ \ensuremath{\left(\omega_{t_{2}}\leq\tau\leq\omega'_{t_{2}}\right)}}$
is the parametric representation for contour $T'_{2}$ with a real-valued
function $\eta*_{t_{2}}$ mapping $[\omega_{t_{2}},\omega'_{t_{2}}]$
onto the interval $[a_{t_{2}},a'_{t_{2}}]$. The lengths of these
transportation contours can be computed and 

\[
\int_{\omega_{t_{2}}}^{\omega_{t_{2}}'}\left|z[\eta_{t_{2}}(\tau)]\right|\eta_{t_{2}}^{'}(\tau)d\tau=\int_{\omega'_{t_{2}}}^{\omega_{t_{2}}}\left|z[\eta*_{t_{2}}(\tau)]\right|\eta*_{t_{2}}^{'}(\tau)d\tau
\]
if, and only, if

\[
z(a_{t_{2}}),\text{ }z'(a'_{t_{2}})\in S_{2}\text{ and }z'(a_{t_{2}})\text{ }z(a'_{t_{2}})\in S_{3}.
\]
The multilevel contour lengths of $S_{2}$ to $S_{3}$ and $S_{3}$
to $S_{2}$ are computed as \linebreak{}
\begin{align}
L[S_{2},S_{3}] & =\int_{\delta_{s_{2}}}^{\delta'_{s_{2}}}\left|z[\varepsilon_{s_{2}}(\tau)]\right|\varepsilon_{s_{2}}^{'}(\tau)d\tau+\int_{\omega_{t_{2}}}^{\omega_{t_{2}}'}\left|z[\eta_{t_{2}}(\tau)]\right|\eta_{t_{2}}^{'}(\tau)d\tau\nonumber \\
 & +\int_{\delta_{s_{3}}}^{\delta'_{s_{3}}}\left|z[\varepsilon_{s_{3}}(\tau)]\right|\varepsilon_{s_{3}}^{'}(\tau)d\tau\label{eq:Ls2S3}
\end{align}
\begin{align}
L[S_{3},S_{2}] & =\int_{\delta_{s_{2}}}^{\delta'_{s_{2}}}\left|z[\varepsilon_{s_{2}}(\tau)]\right|\varepsilon_{s_{2}}^{'}(\tau)d\tau+\nonumber \\
 & \int_{\omega'_{t_{2}}}^{\omega_{t_{2}}}\left|z[\eta*_{t_{2}}(\tau)]\right|\eta*_{t_{2}}^{'}(\tau)d\tau+\nonumber \\
 & \int_{\delta_{s_{3}}}^{\delta'_{s_{3}}}\left|z[\varepsilon_{s_{3}}(\tau)]\right|\varepsilon_{s_{3}}^{'}(\tau)d\tau\label{eq:LS3S2}
\end{align}

The transportation contour $T_{3}$ joins $S_{3}$ to $S_{4}$ and
the return transportation contour $T'_{3}$ joins $S_{4}$ to $S_{4}.$
The contour $T_{3}$ is described by $z(u_{t_{3}})$ for $a_{t_{3}}\leq u_{t_{3}}\leq a'_{t_{3}}$
where $a_{t_{3}},$ $a'_{t_{3}}\in\mathbb{R}$. The starting point
of $T_{3}$ is in the set $S_{3}$ and the ending point of $T_{3}$
is in the set $S_{4}$. The function $u_{t_{3}}=\eta_{t_{3}}(\tau)\text{ \ensuremath{\left(\omega_{t_{3}}\leq\tau\leq\omega'_{t_{3}}\right)}}$
is the parametric representation for $T_{3}$ with a real-valued function
$\eta_{t_{3}}$ mapping $[\omega_{t_{3}},\omega'_{t_{3}}]$ onto the
interval $[a_{t_{3}},a'_{t_{3}}]$. The contour $T'_{3}$ is described
by $z(u'_{t_{3}})$ for $a_{t_{3}}\leq u'_{t_{3}}\leq a'_{t_{3}}$
where $a_{t_{3}},$ $a'_{t_{3}}\in\mathbb{R}$, and the function $u'_{t_{3}}=\eta*_{t_{3}}(\tau)\text{ \ensuremath{\left(\omega_{t_{3}}\leq\tau\leq\omega'_{t_{3}}\right)}}$
is the parametric representation for $T'_{3}$ with a real-valued
function $\eta*_{t_{3}}$ mapping $[\omega_{t_{3}},\omega'_{t_{3}}]$
onto the interval $[a_{t_{3}},a'_{t_{3}}]$. The lengths of these
transportation contours can be computed and 

\[
\int_{\omega_{t_{3}}}^{\omega_{t_{3}}'}\left|z[\eta_{t_{3}}(\tau)]\right|\eta_{t_{3}}^{'}(\tau)d\tau=\int_{\omega'_{t_{3}}}^{\omega_{t_{3}}}\left|z[\eta*_{t_{3}}(\tau)]\right|\eta*_{t_{3}}^{'}(\tau)d\tau
\]
if, and only, if

\[
z(a_{t_{3}}),\text{ }z'(a'_{t_{3}})\in S_{3}\text{ and }z'(a_{t_{3}})\text{ }z(a'_{t_{3}})\in S_{4}.
\]
The multilevel contour lengths of $S_{2}$ to $S_{3}$ and $S_{3}$
to $S_{2}$ are computed as 
\begin{align}
L[S_{3},S_{4}] & =\int_{\delta_{s_{3}}}^{\delta'_{s_{3}}}\left|z[\varepsilon_{s_{3}}(\tau)]\right|\varepsilon_{s_{3}}^{'}(\tau)d\tau+\int_{\omega_{t_{3}}}^{\omega_{t_{3}}'}\left|z[\eta_{t_{3}}(\tau)]\right|\eta_{t_{3}}^{'}(\tau)d\tau\nonumber \\
 & +\int_{\delta_{s_{4}}}^{\delta'_{s_{4}}}\left|z[\varepsilon_{s_{4}}(\tau)]\right|\varepsilon_{s_{4}}^{'}(\tau)d\tau\label{eq:LS3S4}
\end{align}
\begin{align}
L[S_{4},S_{3}] & =\int_{\delta_{s_{3}}}^{\delta'_{s_{3}}}\left|z[\varepsilon_{s_{3}}(\tau)]\right|\varepsilon_{s_{3}}^{'}(\tau)d\tau+\int_{\omega'_{t_{3}}}^{\omega_{t_{3}}}\left|z[\eta*_{t_{3}}(\tau)]\right|\eta*_{t_{3}}^{'}(\tau)d\tau\nonumber \\
 & +\int_{\delta_{s_{4}}}^{\delta'_{s_{4}}}\left|z[\varepsilon_{s_{4}}(\tau)]\right|\varepsilon_{s_{4}}^{'}(\tau)d\tau.\label{eq:LS4S3}
\end{align}
The total lengths of multilevel contours from $S_{1}$ to $S_{4}$
and from $S_{4}$ to $S_{1}$ are obtained by
\begin{align}
L[S_{1},S_{4}]=\nonumber \\
 & \sum_{i=1}^{4}\int_{\delta_{s_{i}}}^{\delta'_{s_{i}}}\left|z[\varepsilon_{s_{i}}(\tau)]\right|\varepsilon_{s_{i}}^{'}(\tau)d\tau+\sum_{i=1}^{3}\int_{\omega_{t_{i}}}^{\omega_{t_{i}}'}\left|z[\eta_{t_{i}}(\tau)]\right|\eta_{t_{i}}^{'}(\tau)d\tau\label{eq:LS1S4}
\end{align}
\begin{align}
L[S_{4},S_{1}]=\nonumber \\
 & \sum_{i=1}^{4}\int_{\delta_{s_{i}}}^{\delta'_{s_{i}}}\left|z[\varepsilon_{s_{i}}(\tau)]\right|\varepsilon_{s_{i}}^{'}(\tau)d\tau+\sum_{i=1}^{3}\int_{\omega'_{t_{i}}}^{\omega_{t_{i}}}\left|z[\eta*_{t_{i}}(\tau)]\right|\eta*_{t_{i}}^{'}(\tau)d\tau.\label{eq:LS4s1}
\end{align}

Given the fixed shapes of contours on different planes, as shown arbitrarily
in Figure \ref{fig:Transportation-of-information}, the above constructions
of lengths and transportation contours are to be treated as an example
of the usefulness of multilevel contours. Such constructions can be
extended for several other practical situations arising from the data.
Note that the contours $T_{i}$ and $T'_{i}$ pass through the line
created by $\mathbb{C}_{i}\cap\mathbb{C}_{i+1}$ for $i=1,2,3$. So
the corresponding integrals of the second term of the R.H.S. of (\ref{eq})
represent combined lengths created due to traveling of the contour
$T_{i}$ from a point in $S_{i}$ to a point in $\mathbb{C}_{i}\cap\mathbb{C}_{i+1}$
and then traveling from a point in $\mathbb{C}_{i}\cap\mathbb{C}_{i+1}$
to $S_{i+1}.$ Similarly, the integrals of the second term of the
R.H.S. of (\ref{eq}) represent combined lengths created due to traveling
of the contour $T'_{i}$ from a point in $S_{i+1}$ to a point in
$\mathbb{C}_{i}\cap\mathbb{C}_{i+1}$ and then traveling from a point
in $\mathbb{C}_{i}\cap\mathbb{C}_{i+1}$ to $S_{i}.$ Next, we will
see how this combined integral can be subdivided into smaller integrals
while computing the shortest distance. 

Since $T_{i}$ was described by $z(u_{t_{i}})$ for $a_{t_{i}}\leq u_{t_{i}}\leq a'_{t_{i}}$,
we further partition the $T_{i}$ into three contours mentioned in
the previous paragraph. Suppose $z(a_{t_{i}}^{m_{1}})$ be the point
on $S_{i}$ for $a_{t_{i}}^{m_{1}}\in[a_{t_{i}},a'_{t_{i}}]$ that
is the closest to the point on the line (say, $z(a_{t_{i}}^{m_{2}})$
created due to $\mathbb{C}_{i}\cap\mathbb{C}_{i+1}$, and $z(a_{t_{i}}^{m_{3}})$
be the point on the line created due to $\mathbb{C}_{i}\cap\mathbb{C}_{i+1}$
that is closest to $S_{i+1}$. The point at which the contour from
$z(a_{t_{i}}^{m_{3}})$ joins $S_{i+1}$ say, $z(a_{t_{i}}^{m_{4}})$
for $a_{t_{i}}^{m_{4}}\in[a_{t_{i}},a'_{t_{i}}]$ that is the closest
from a point on the line created due to $\mathbb{C}_{i}\cap\mathbb{C}_{i+1}$.
Here $z(a_{t_{i}}^{m_{1}})$ and $z(a_{t_{i}}^{m_{4}})$ are complex
numbers on different complex planes. The complex numbers $z(u_{t_{i}})$
partitioned as

\[
\begin{array}{cc}
z(u_{t_{i}})= & \left\{ \begin{array}{cc}
z(u_{t_{i}}^{I}) & (a_{t_{i}}\leq u_{t_{i}}<a_{t_{i}}^{I})\\
\\
z(u_{t_{i}}^{II}) & (a_{t_{i}}^{I}\leq u_{t_{i}}<a_{t_{i}}^{II})\\
\\
z(u_{t_{i}}^{III}) & (a_{t_{i}}^{II}\leq u_{t_{i}}\leq a'_{t_{i}})
\end{array}\right),\end{array}
\]
such that 

\[
(a_{t_{i}}\leq u_{t_{i}}<a_{t_{i}}^{I})\cup(a_{t_{i}}^{I}\leq u_{t_{i}}<a_{t_{i}}^{II})\cup(a_{t_{i}}^{II}\leq u_{t_{i}}\leq a'_{t_{i}})=[a_{t_{i}},a'_{t_{i}}].
\]
Let us re-define the function $u_{t_{i}}$ below to represent the
three partitions mentioned above.
\[
\begin{array}{cc}
u_{t_{i}}= & \left\{ \begin{array}{cc}
\eta_{t_{i}}^{I} & (\omega_{t_{i}}\leq\tau<\omega_{t_{i}}^{I})\\
\\
\eta_{t_{i}}^{II} & (\omega_{t_{i}}^{I}\leq\tau<\omega_{t_{i}}^{II})\\
\\
\eta_{t_{i}}^{I} & (\omega_{t_{i}}^{II}\leq\tau\leq\omega'_{t_{i}})
\end{array}\right),\end{array}
\]
such that
\[
(\omega_{t_{i}}\leq\tau<\omega_{t_{i}}^{I})\cup(\omega_{t_{i}}^{I}\leq\tau<\omega_{t_{i}}^{II})\cup(\omega_{t_{i}}^{II}\leq\tau\leq\omega'_{t_{i}})=[\omega_{t_{i}},\omega'_{t_{i}}].
\]
The three shortest distances arise out of above partitions are, say,
$L(I:T_{i})$, $L(II:T_{i})$, and $L(III:T_{i})$ . These shortest
distances are given by

\begin{equation}
L(I:T_{i})=\int_{\omega_{t_{i}}}^{\omega_{t_{i}}^{I}}\left|z[\eta_{t_{i}}^{I}(\tau)]\right|\left(\eta_{t_{i}}^{I}\right)'(\tau)d\tau\label{eq:L-I}
\end{equation}

\begin{equation}
L(II:T_{i})=\int_{\omega_{t_{i}}^{I}}^{\omega_{t_{i}}^{II}}\left|z[\eta_{t_{i}}^{II}(\tau)]\right|\left(\eta_{t_{i}}^{II}\right)'(\tau)d\tau\label{eq:LII}
\end{equation}

\begin{equation}
L(I:T_{i})=\int_{\omega_{t_{i}}^{II}}^{\omega'_{t_{i}}}\left|z[\eta_{t_{i}}^{III}(\tau)]\right|\left(\eta_{t_{i}}^{III}\right)'(\tau)d\tau\label{eq:LIII}
\end{equation}
The shortest transportation contour from $S_{i}$ to $S_{i+1}$ is
computed by using (\ref{eq:L-I}) through (\ref{eq:LIII})

\begin{equation}
L(T_{i}^{m})=L(I:T_{i})+L(II:T_{i})+L(III:T_{i}).\label{eq:shortestTi}
\end{equation}
 The shortest transportation contour from $S_{i+1}$ to $S_{i}$ would
be the same as in (\ref{eq:shortestTi}). Next we compute the farthest
transportation contour that joins $S_{i+1}$ from $S_{i}.$ We partition
the $T_{i}$ described by $z(u_{t_{i}})$ for $a_{t_{i}}\leq u_{t_{i}}\leq a'_{t_{i}}$
into three contours that represent longest contour drawn from $S_{i}$
to $S_{i+1}.$ Suppose $z(a_{t_{i}}^{M_{1}})$ be the point on $S_{i}$
for $a_{t_{i}}^{M_{1}}\in[a_{t_{i}},a'_{t_{i}}]$ that is the farthest
to the point on the line (say, $z(a_{t_{i}}^{M_{2}})$ created due
to $\mathbb{C}_{i}\cap\mathbb{C}_{i+1}$, and $z(a_{t_{i}}^{M_{3}})$
be the point on the line created due to $\mathbb{C}_{i}\cap\mathbb{C}_{i+1}$
that is farthest to $S_{i+1}$. The point at which the contour from
$z(a_{t_{i}}^{M_{3}})$ joins $S_{i+1}$ say, $z(a_{t_{i}}^{M_{4}})$
for $a_{t_{i}}^{M_{4}}\in[a_{t_{i}},a'_{t_{i}}]$ that is the farthest
from a point on the line created due to $\mathbb{C}_{i}\cap\mathbb{C}_{i+1}$.
The complex numbers $z(u_{t_{i}})$ partitioned as

\[
\begin{array}{cc}
z(u_{t_{i}})= & \left\{ \begin{array}{cc}
z(u_{t_{i}}^{IV}) & (a_{t_{i}}\leq u_{t_{i}}<a_{t_{i}}^{IV})\\
\\
z(u_{t_{i}}^{V}) & (a_{t_{i}}^{IV}\leq u_{t_{i}}<a_{t_{i}}^{V})\\
\\
z(u_{t_{i}}^{VI}) & (a_{t_{i}}^{V}\leq u_{t_{i}}\leq a'_{t_{i}})
\end{array}\right),\end{array}
\]
such that 

\[
(a_{t_{i}}\leq u_{t_{i}}<a_{t_{i}}^{IV})\cup(a_{t_{i}}^{IV}\leq u_{t_{i}}<a_{t_{i}}^{V})\cup(a_{t_{i}}^{VI}\leq u_{t_{i}}\leq a'_{t_{i}})=[a_{t_{i}},a'_{t_{i}}].
\]
Let us re-define the function $u_{t_{i}}$ below to represent the
three partitions mentioned above.
\[
\begin{array}{cc}
u_{t_{i}}= & \left\{ \begin{array}{cc}
\eta_{t_{i}}^{IV} & (\omega_{t_{i}}\leq\tau<\omega_{t_{i}}^{IV})\\
\\
\eta_{t_{i}}^{V} & (\omega_{t_{i}}^{IV}\leq\tau<\omega_{t_{i}}^{V})\\
\\
\eta_{t_{i}}^{VI} & (\omega_{t_{i}}^{V}\leq\tau\leq\omega'_{t_{i}})
\end{array}\right),\end{array}
\]
such that
\[
(\omega_{t_{i}}\leq\tau<\omega_{t_{i}}^{IV})\cup(\omega_{t_{i}}^{IV}\leq\tau<\omega_{t_{i}}^{V})\cup(\omega_{t_{i}}^{VI}\leq\tau\leq\omega'_{t_{i}})=[\omega_{t_{i}},\omega'_{t_{i}}].
\]
The three longest distances arise out of above partitions are, say,
$L(IV:T_{i})$, $L(V:T_{i})$, and $L(VI:T_{i})$ . These longest
distances are given by

\begin{equation}
L(IV:T_{i})=\int_{\omega_{t_{i}}}^{\omega_{t_{i}}^{IV}}\left|z[\eta_{t_{i}}^{IV}(\tau)]\right|\left(\eta_{t_{i}}^{IV}\right)'(\tau)d\tau\label{eq:L-I-1}
\end{equation}

\begin{equation}
L(V:T_{i})=\int_{\omega_{t_{i}}^{I}}^{\omega_{t_{i}}^{V}}\left|z[\eta_{t_{i}}^{V}(\tau)]\right|\left(\eta_{t_{i}}^{V}\right)'(\tau)d\tau\label{eq:LII-1}
\end{equation}

\begin{equation}
L(VI:T_{i})=\int_{\omega_{t_{i}}^{II}}^{\omega'_{t_{i}}}\left|z[\eta_{t_{i}}^{VI}(\tau)]\right|\left(\eta_{t_{i}}^{VI}\right)'(\tau)d\tau\label{eq:LIII-1}
\end{equation}
The shortest transportation contour from $S_{i}$ to $S_{i+1}$ is
computed by using (\ref{eq:L-I-1}) through (\ref{eq:LIII-1})

\begin{equation}
L(T_{i}^{M})=L(IV:T_{i})+L(V:T_{i})+L(VI:T_{i}).\label{eq:shortestTi-1}
\end{equation}

The shortest distance and the longest distance of a multilevel contour
give us an idea about the transportation contour $T_{i}$, $T'_{i}$,
and the range of times that they carry information from one contour
to another contour. The shorter the distance the quicker is the information
transformed between two contours in different planes and the longer
the transportation contour, the longer is the time for transporting
the information. The time taken to reach from $z_{1}$ to $z_{2}$
in a plane is assumed to be proportional to the distance between $z_{1}$
and $z_{2}$. A transportation contour is assumed here to carry the
information on the shape of the contour, this carries the location
of data of points on the contour in that plane and joins with another
contour in another plane. In this way, all the contours lying in different
places are joined so that combined information on contours (multilevel
contours) is constructed. The contours lying in different planes are
otherwise disjoint. By combining contours in different planes an information
tree is attained that has locations of all the points lying in various
contours of a bundle. 
\begin{thm}
\label{thm:(Spinning-of-bundle)}(Spinning of bundle theorem) Suppose
the bundle $B_{\mathbb{R}}(\mathbb{C})$ is rotated anti-clockwise
such that the line passing through $(0,0)$ of all the planes within
$B_{\mathbb{R}}(\mathbb{C})$ forms an angle $\theta$ $(\theta>0)$
with $y-$axis. Suppose the rotation is continued for each $\theta\in(0,360]\subset\mathbb{R}^{+}.$
Then the space created due to such a rotation forms a $1-1$ correspondence
with $B_{\mathbb{R}}(\mathbb{C}).$
\end{thm}

\begin{proof}
Suppose we make a copy of the bundle $B_{\mathbb{R}}(\mathbb{C})$
combined with the positioning of $\mathbb{C}_{0}$ and place it on
the bundle such that these two bundles occupy exactly the same space.
Let us call the original bundle with the positioning of $\mathbb{C}_{0}$
as $B_{\mathbb{R}}^{o}(\mathbb{C})$ and its as $B_{\mathbb{R}}^{c}(\mathbb{C})$.
Suppose we tilt $B_{\mathbb{R}}^{c}(\mathbb{C})$ to the left such
that $B_{\mathbb{R}}^{c}(\mathbb{C})$ inclined at an angle $\theta$
for $\theta>0$ with $y-$axis. See Figure \ref{fig:Angle-between-two}.
The points (complex numbers) on the plane $\mathbb{C}_{0}$ do not
change with this tilting. So as the points in the space created by
$B_{\mathbb{R}}^{c}(\mathbb{C}).$ Let us consider a plane $\mathbb{C}_{p}$
before tilting for $\mathbb{C}_{p}\in B_{\mathbb{R}}^{o}(\mathbb{C}).$
The same $\mathbb{C}_{p}$ in $B_{\mathbb{R}}^{c}(\mathbb{C})$ is
now inclined away at an angle $\theta.$ Let us call the copied bundle
$B_{\mathbb{R}}^{c}(\mathbb{C})$ that is inclined at an angle $\theta$
be $B_{\mathbb{R}}^{c}(\mathbb{C},\theta)$. Each value of $\mathbb{C}_{p}$
that was there when $\mathbb{C}_{p}\in B_{\mathbb{R}}^{o}(\mathbb{C},\theta)$
is still there after inclination. However, $\mathbb{C}_{p}$ in $B_{\mathbb{R}}^{c}(\mathbb{C},\theta)$
intersects with infinite set of planes of $B_{\mathbb{R}}^{o}(\mathbb{C}).$
Next, we show that $\mathbb{C}_{p}$ in $B_{\mathbb{R}}^{c}(\mathbb{C},\theta)$
has the same points which are a subset of $B_{\mathbb{R}}^{o}(\mathbb{C}).$ 

Since $\mathbb{C}_{p}$ in $B_{\mathbb{R}}^{c}(\mathbb{C},\theta)$
intersects with infinitely many (uncountable) planes. So by Theorem
\ref{thm:The-shortest-possible}, one can draw a contour that passes
through all the infinite planes of $B_{\mathbb{R}}^{o}(\mathbb{C}).$
There is no point of $\mathbb{C}_{p}$ that does not intersect with
points of the bundle $B_{\mathbb{R}}^{o}(\mathbb{C}).$ This brings
the conclusion that 
\begin{equation}
\mathbb{C}_{p}\in B_{\mathbb{R}}^{c}(\mathbb{C},\theta)\implies\mathbb{C}_{p}\subset B_{\mathbb{R}}^{o}(\mathbb{C}).\label{eq:Cpare same}
\end{equation}
This means, 
\begin{equation}
\text{for every }z\in\mathbb{C}_{p}\text{ implies }\text{\ensuremath{z\in}}B_{\mathbb{R}}^{o}(\mathbb{C}).\label{eq:zinCpimplieszinbundle}
\end{equation}

\begin{figure}
\includegraphics[scale=0.19]{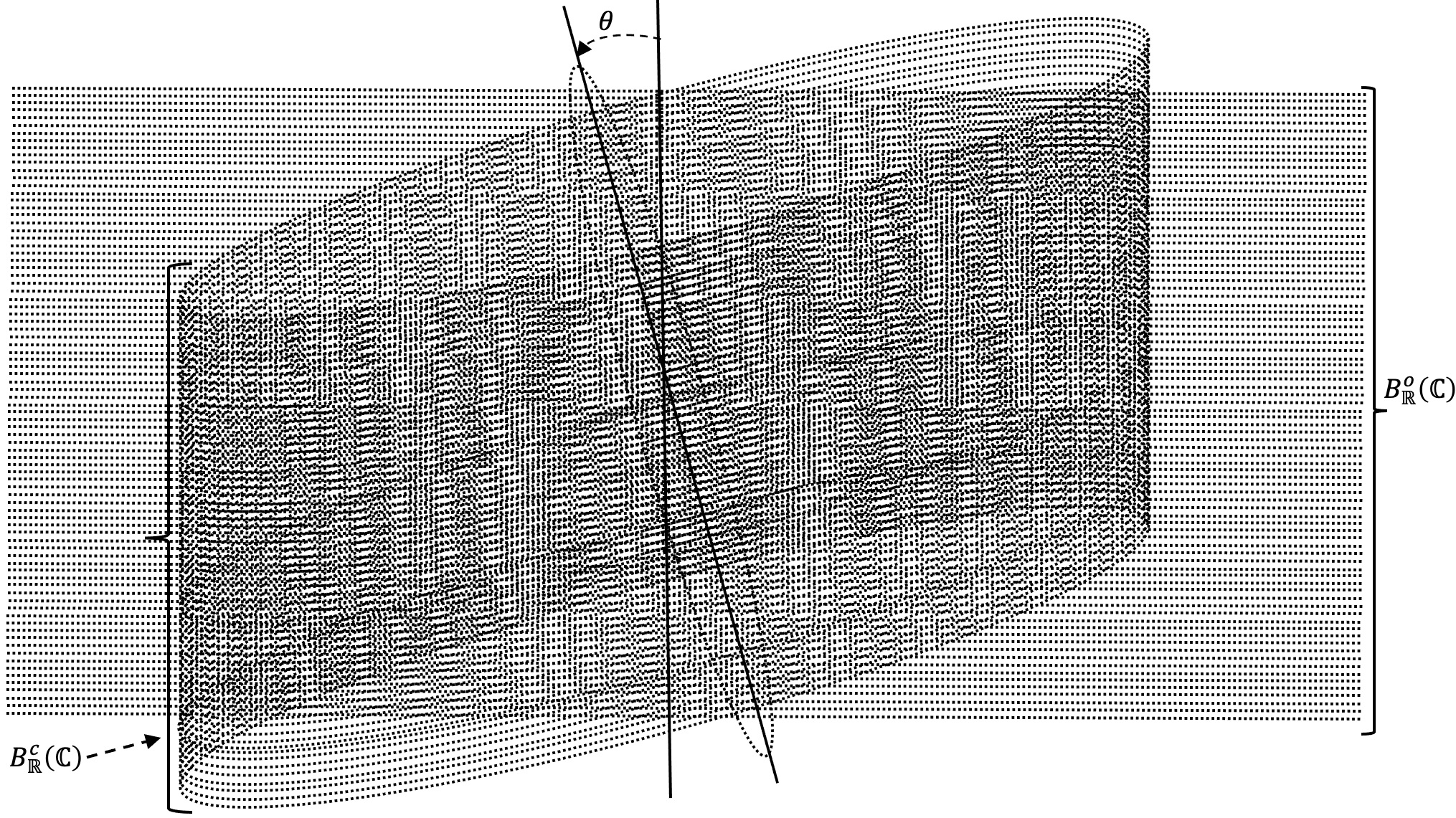}

\caption{\label{fig:Angle-between-two}Angle between two bundles $B_{\mathbb{R}}^{o}(\mathbb{C})$
and $B_{\mathbb{R}}^{c}(\mathbb{C})$ and rotation of the bundle $B_{\mathbb{R}}^{c}(\mathbb{C})$
over $y-$axis.}

\end{figure}

Statement (\ref{eq:zinCpimplieszinbundle}) is true for every $\theta>0$
and every $B_{\mathbb{R}}^{c}(\mathbb{C},\theta).$ Since $\theta$
is arbitrary, the statement (\ref{eq:zinCpimplieszinbundle}) is true
for each $\theta\in(0,360].$ Hence, the space created by $B_{\mathbb{R}}^{c}(\mathbb{C},\theta)$
for all $\theta$ values for $\theta\in(0,360^{0}]$ is the same as
the space of the bundle $B_{\mathbb{R}}(\mathbb{C}).$
\end{proof}
\begin{rem}
The $n-$dimensional complex plane $\mathbb{C}^{n}$ is also subset
of the space created by the rotation in the Theorem \ref{thm:(Spinning-of-bundle)}. 
\end{rem}

\section{\textbf{Multilevel Contours in a Random Environment}}

Let us consider the bundle $B_{\mathbb{R}}(\mathbb{C})$ combined
with the $\mathbb{C}_{0}$ intersecting with the bundle as in Section
2. Let the intersection of $\mathbb{C}_{0}$ on the bundle be arbitrary.
Choose a plane $\mathbb{C}_{l}$ for $\mathbb{C}_{l}\in B_{\mathbb{R}}(\mathbb{C}).$
Let $X\left(\gamma_{l},\mathbb{C}_{l}\right)$ be a random variable
describing the position of the complex number on a contour $\gamma_{l}$
at time $t$ in the plane $\mathbb{C}_{l}$. Here $\gamma_{l}$ is
described by $z_{l}(t)$ $(t_{0}\leq t\leq b_{0})$ for $t_{0},b_{0}\in\mathbb{R}$
and $z_{l}(t_{0})=$ $z_{l_{0}}$, say. Let $X\left(z_{l_{0}},\mathbb{C}_{l}\right)$
be the position of the random variable at $t_{0}$ in $\mathbb{C}_{l}$.
Here $X\left(\gamma_{l},\mathbb{C}_{l}\right)$ can be treated like
a standard measurable function. Suppose $X\left(z_{l_{0}},\mathbb{C}_{l}\right)$
picked arbitrarily in the plane $\mathbb{C}_{l}$. Here $X\left(z_{l}(t),\mathbb{C}_{l}\right)$
for $t\in[t_{0},\infty)$ and $z_{l}(t)\in\mathbb{C}_{l}$ is a stochastic
process on the complex plane $\mathbb{C}_{l}$. Once $X\left(z_{l_{0}},\mathbb{C}_{l}\right)$
is chosen, the location of $X\left(z_{l_{1}},\mathbb{C}_{l}\right)$
for $t_{1}>t_{0}$ could be anywhere within an open disc $D\left(z_{l_{0}}.r_{0}\right)$,
where $X\left(z_{l_{1}},\mathbb{C}_{l}\right)$ is the complex number
picked by $X\left(z_{l}(t),\mathbb{C}_{l}\right)$ at $t_{1}$ and
\[
D\left(z_{l_{0}}.r_{0}\right)=\left\{ z:\left|z-z_{l_{0}}\right|<r_{0}\text{ for \ensuremath{r_{0}>0} and }z\in\mathbb{C}_{l}\right\} .
\]
\begin{figure}
\includegraphics[scale=0.5]{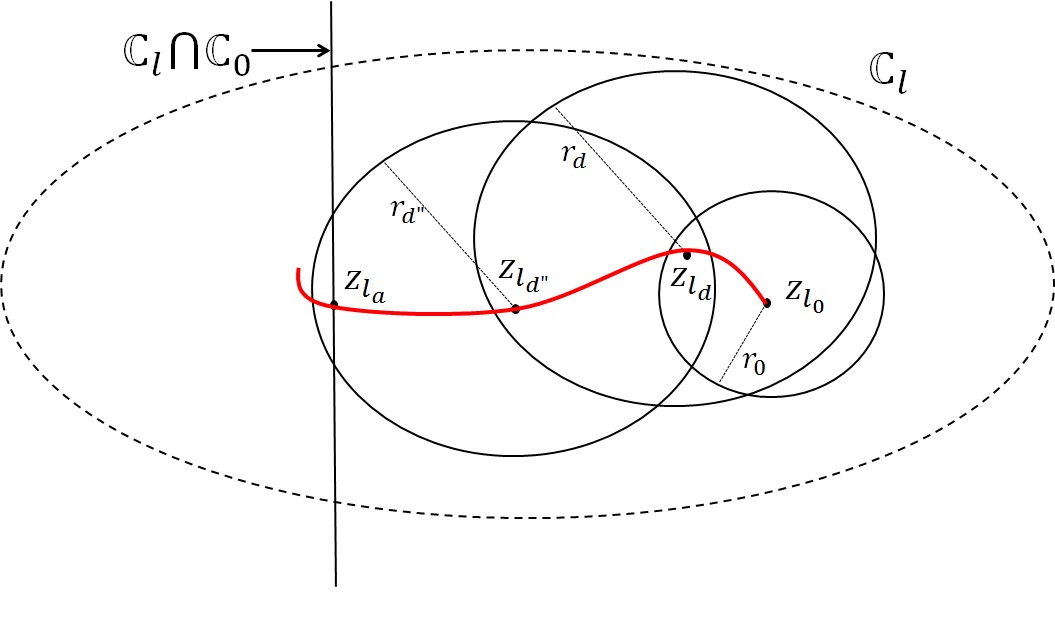}

\caption{\label{fig:FormationofcontoursDISCS}Formation of contours through
two-step randomness.}

\end{figure}
Suppose $X\left(z_{l_{0}},\mathbb{C}_{l}\right)=z_{l_{1}}$. The selection
of radius is done randomly at each step of selecting a complex number.
Once $X\left(z_{l},\mathbb{C}_{l}\right)$ chooses a number at $t_{0}$
then during $[t_{0},t_{1}]$ a value for $r_{0}$ for $r_{0}>0$ is
chosen randomly to build $D\left(z_{l_{0}}.r_{0}\right).$ Once the
set of points of $D\left(z_{l_{0}}.r_{0}\right)$ are available then
$X\left(z_{l_{0}},\mathbb{C}_{l}\right)$ will choose the second number
during $[t_{0},t_{1}].$ See Figure \ref{fig:FormationofcontoursDISCS}.
This procedure of two-step randomness continues forever in the intervals
\[
\left\{ [t_{0},t_{a}],(t_{a},t_{a'}],(t_{a'},t_{b}],(t_{b},t_{b'],...}\right\} 
\]
where $a,a',b,b',...\in\mathbb{R}$ and $t_{0}<t_{a}<t_{a'}<t_{b}<t_{b'<...}.$
We can draw a contour from $z_{l_{0}}$to $z_{l_{1}}$ using the set
of points generated by $z_{l}(t)$ within the time interval $[t_{0},t_{1})$.
The set of all possible values the set $\left\{ X\left(z_{l},\mathbb{C}_{l}\right)\right\} $
can take over the interval $t_{0},\infty)$ is called the state-space
of $X\left(z_{l}(t),\mathbb{C}_{l}\right).$ This continuous-time
stochastic process does not stop once the initial value on the plane
$\mathbb{C}_{l}$ is chosen and its state-space is $\mathbb{C}_{l}$.
The described process is a continuous time and continuous space stochastic
process. Let $p\left(z_{l}(t_{1}),\mathbb{C}_{l}\right)$ be the probability
function that is associated with $X\left(z_{l}(t_{1}),\mathbb{C}_{l}\right)$
such that $p\left(z_{l}(t_{1}),\mathbb{C}_{l}\right)$ describes the
probability that $X\left(z_{l}(t_{1}),\mathbb{C}_{l}\right)$ picks
$z_{l_{1}}$for $z_{l_{1}}=z_{l}(t_{1})$ $(t_{1}>t_{0}).$ The probability
function defined as

\begin{equation}
p\left(z_{l}(t_{1}),t_{1}\right)=\text{Prob}\left[X\left(z_{l}(t_{1}),\mathbb{C}_{l}\right)=z_{l_{1}}\right]\text{ for }t_{1}>t_{0}.\label{eq:DEFprobability function}
\end{equation}
Once $X\left(z_{l}(t_{1}),\mathbb{C}_{l}\right)$ is generated, then
one can draw contour from $z_{l_{0}}$to $z_{l_{1}}$and compute $L(z_{l_{0}},z_{l_{1}})$
the length from $z_{l_{0}}$to $z_{l_{1}}$. The contour $\gamma_{l}$
is described by 
\[
\gamma_{l}=z_{l}(t)\text{ for }t\in[t_{0},\infty)
\]
 and 
\[
t=v_{l}(\tau)\text{ }\left(a_{0}^{l_{0}}\leq t<a_{0}^{l_{1}}\right)
\]
is the parametric representation for $\gamma_{l}$ with a real valued
function $v_{l}$ mapping $[a_{0}^{l_{0}},a_{0}^{l_{1}})$ onto the
interval $[t_{0},t_{1}).$ This gives us,
\begin{equation}
L(z_{l_{0}},z_{l_{1}})=\int_{a_{0}^{l_{0}}}^{a_{0}^{l_{1}}}\left|z\left[v_{l}(\tau)\right]\right|v'_{l}(\tau)d\tau.\label{eq:lengthLzl0zl1}
\end{equation}
Our main focus here in this article is on multilevel contours. So
initially we assume that once the variable $X\left(z_{l}(t_{1}),\mathbb{C}_{l}\right)$
is picked a value at $t_{0}$, it can not reach the same value for
$t>t_{0},$ i.e.,

\begin{equation}
X\left(z_{l}(t_{0}),\mathbb{C}_{l}\right)\neq X\left(z_{l}(t),\mathbb{C}_{l}\right)\text{ for }t>t_{0}.\label{eq:zloNEQzl1}
\end{equation}
For simplicity in understanding, we might denote (at certain specific
places in the article) a complex number generated at each iteration
with an integer index. However, in reality the numbers chosen by $X\left(z_{l}(t),\mathbb{C}_{l}\right)$
are uncountably infinite. 
\begin{defn}
\textbf{\label{def:Distinct-complex-numbers}Distinct complex numbers
by $X\left(z_{l}(t),\mathbb{C}_{l}\right)$}: Suppose $X\left(z_{l}(t_{0}),\mathbb{C}_{l}\right)=z_{l_{0}}$
for $z_{l_{0}}\in\mathbb{C}_{l}.$ Suppose $X\left(z_{l}(t),\mathbb{C}_{l}\right)=z_{l_{\alpha}}$
for $t\in[t_{0},t_{\alpha}]$ and $X\left(z_{l}(t),\mathbb{C}_{l}\right)=z_{l_{\beta}}$
for $t\in(t_{\alpha},t_{\beta}]$ for all $t_{0}<t_{\alpha}<t_{\beta}$,
and $\alpha,\beta,$$t_{0},t_{\alpha},t_{\beta}\in\mathbb{R}$ then
$z_{l_{0}}\neq z_{l_{\alpha}}\neq z_{l_{\beta}}.$ This property assures
that $X\left(z_{l}(t),\mathbb{C}_{l}\right)$ cannot choose the same
number that it already chose in any of the previous time intervals
after the initial number is chosen (including the initial number). 
\end{defn}

When a random variable $X\left(z_{l}(t),\mathbb{C}_{l}\right)$ chooses
a number within the disc created and that value (number) has been
chosen already and was part of the contour, then, $X\left(z_{l}(t),\mathbb{C}_{l}\right)$
will choose another number in the disc. This procedure continues until
a distinct number is chosen by $X\left(z_{l}(t),\mathbb{C}_{l}\right)$.
Such an assumption in (\ref{def:Distinct-complex-numbers}) or in
(\ref{eq:zloNEQzl1}) will allow quicker formation of multilevel contour.
We will later see the consequences of relaxing the assumption in (\ref{def:Distinct-complex-numbers}).
We have

\begin{equation}
z_{l_{1}}\in D\left(z_{l_{0}}.r_{0}\right)\text{ and }z_{l_{0}}\neq z_{l_{1}}.\label{eq:zloNEQzl1inthedisc}
\end{equation}
Let $z_{l}(t_{2})$ be the value of $X\left(z_{l}(t_{1}),\mathbb{C}_{l}\right)$
at $t_{2}$ for $t_{2}>t_{1}$ and $z_{l}(t_{2})\in D\left(z_{l_{1}}.r_{1}\right)$
for $r_{1}>0$ such that 
\[
p\left(z_{l}(t_{2}),t_{2}\right)=\text{Prob}\left[X\left(z_{l}(t_{1}),\mathbb{C}_{l}\right)=z_{l_{2}}\right]\text{ for }t_{2}>t_{1}.
\]
The contour $\gamma_{l}$ with a new parametric representation 
\[
t=v_{l}(\tau)\text{ }\left(a_{0}^{l_{1}}\leq t<a_{0}^{l_{2}}\right)
\]
with a real valued function $v_{l}$ mapping $[a_{0}^{l_{1}},a_{0}^{l_{2}})$
onto the interval $[t_{1},t_{2})$ helps us to compute the length
$L(z_{l_{1}},z_{l_{2}})$ from $z_{l_{1}}$to $z_{l_{2}}=z_{l}(t_{2})$
from $z_{l_{1}}$ for $t_{2}>t_{1}.$ That is,

\begin{align*}
p\left(z_{l}(t_{2}),t_{2}\right)= & \text{Prob}\left[X\left(z_{l}(t_{2}),\mathbb{C}_{l}\right)=z_{l}(t_{2})/X\left(z_{l}(t_{0}),\mathbb{C}_{l}\right)=z_{l_{0}},\right.\\
 & \left.X\left(z_{l}(t_{1}),\mathbb{C}_{l}\right)=z_{l}(t_{1})\right]
\end{align*}
\begin{equation}
L(z_{l_{1}},z_{l_{2}})=\int_{a_{0}^{l_{1}}}^{a_{0}^{l_{2}}}\left|z\left[v_{l}(\tau)\right]\right|v'_{l}(\tau)d\tau.\label{eq:lengthLzl1zl2}
\end{equation}
Note that, we are not drawing contours from $z_{l_{0}}$to $z_{l_{2}}$
because $X\left(z_{l}(t_{1}),\mathbb{C}_{l}\right)$ will change to
$z_{l_{2}}$for $t_{2}>t_{1}.$ In fact, under the construction explained,
the contour will start $z_{l_{0}}$and reach $z_{l_{2}}$only through
the point $z_{l_{1}}.$ As these are new ideas, we have explained
above the probabilities of picking various values by $X\left(z_{l}(t_{1}),\mathbb{C}_{l}\right)$.
We will slightly re-define below the probabilities and their transitions
to accommodate an easier understanding of these concepts. 

There are infinitely many options around $z_{l_{0}}$ that $X\left(z_{l}(t_{1}),\mathbb{C}_{l}\right)$
can pick during $[t_{0},t_{1}]$ each with a probability $p\left(z_{l}(t),t\right)$
for $t\in[t_{0},t_{1}].$ Let $p\left(z_{l_{0}},z_{l}(t),t\right)$
for $t\in[t_{0},t_{1}]$ that $X\left(z_{l}(t),\mathbb{C}_{l}\right)$
=$z_{l_{0}}$at $t_{0}$ has transitioned to $X\left(z_{l}(t),\mathbb{C}_{l}\right)$
=$z_{l_{1}}$at during $[t_{0},t_{1}].$ For all such probabilities
of transitions during $[t_{0},t_{1}]$, we will have 

\begin{equation}
\int_{t_{0}}^{t_{1}}p\left(z_{l_{0}},z_{l}(t),t\right)dt=1,\label{eq:pdfatzl0}
\end{equation}
where $p\left(z_{l_{0}},z_{l}(t),t\right)$ can be expressed as 

\[
p\left(z_{l_{0}},z_{l}(t),t\right)=\text{Prob}\left[X\left(z_{l}(t),\mathbb{C}_{l}\right)=z_{l_{1}}\left/X\left(z_{l}(t_{0}),\mathbb{C}_{l}\right)=z_{l_{0}}\right.\right].
\]
Let $p\left(z_{l_{0}},z_{l}(t),t\right)$ for $t\in(t_{1},t_{2}]$
be the probability that $X\left(z_{l}(t),\mathbb{C}_{l}\right)$ during
$(t_{1},t_{2}]$ will pick a number within $D(z_{l_{1}},r_{1})$$\in\mathbb{C}_{l}$
among infinitely many options. The transition probability from $z_{l_{1}}$
to $z_{l_{2}}$during $(t_{1},t_{2}]$ is 
\[
p\left(z_{l_{1}},z_{l}(t),t\right)\text{ }\left(t\in(t_{1},t_{2}]\right)=\text{Prob}[X\left(z_{l}(t),\mathbb{C}_{l}\right)=z_{l_{2}}\left/X\left(z_{l}(t_{1}),\mathbb{C}_{l}\right)=z_{l_{1}}\right.]
\]
such that

\[
\int_{t_{1}}^{t_{2}}p\left(z_{l_{1}},z_{l}(t),t\right)dt=1.
\]
Note that, 

\begin{align}
p\left(z_{l_{0}},z_{l}(t),t\right)\text{ }\left(t\in[t_{0},t_{2}]\right) & =\nonumber \\
 & p\left(z_{l_{0}},z_{l}(t),t\right)\text{ }\left(t\in[t_{0},t_{1}]\right)\times\nonumber \\
 & p\left(z_{l_{1}},z_{l}(t),t\right)\text{ }\left(t\in(t_{0},t_{2}]\right).\label{eq:two-step-transitionl0tol2}
\end{align}

A direct transition from the complex number $z_{l_{0}}$ to another
complex number $z_{l_{2}}$is not possible during $[t_{0},t_{2}]$
under the above framework. A direct transition we mean here a one-step
transition. As $X\left(z_{l}(t),\mathbb{C}_{l}\right)$ has reached
$z_{l_{1}}$during $[t_{0},t_{1}]$ and then $X\left(z_{l}(t),\mathbb{C}_{l}\right)$
starting at $z_{l_{1}}$it has taken the value $z_{l_{2}}$ during
$(t_{1},t_{2}].$ This implies, $p\left(z_{l_{0}},z_{l}(t),t\right)\text{ }\left(t\in[t_{0},t_{2}]\right)$
is not possible without having a hopping over the value $X\left(z_{l}(t),\mathbb{C}_{l}\right)=z_{l_{1}}$
during $[t_{0},t_{1}].$ If $z_{l_{2}}$is not picked by $X\left(z_{l}(t),\mathbb{C}_{l}\right)$
during $(t_{1},t_{2}]$ then it can be picked at some future time
by $X\left(z_{l}(t),\mathbb{C}_{l}\right)$ for $t>t_{2}.$ So,

\[
p\left(z_{l_{0}},z_{l}(t),t\right)\text{ }\left(t\in[t_{0},t_{2}]\right)=\begin{cases}
\begin{array}{cc}
0 & \text{(one-step transition)}\\
>0 & \text{(two or more-steps transitions)}
\end{array}\end{cases}
\]

\begin{thm}
\label{thm:A-contour-isMarkovChain}A contour formed by the set of
points generated by $X\left(z_{l}(t),\mathbb{C}_{l}\right)$ on $\mathbb{C}_{l}$
for $\mathbb{C}_{l}$ in the bundle $B_{\mathbb{R}}(\mathbb{C})$
combined with the $\mathbb{C}_{0}$ intersecting with the bundle and
satisfying (\ref{def:Distinct-complex-numbers}) will obey continuous
time Markov property.
\end{thm}

\begin{proof}
Let $\gamma_{l}$ be the contour generated by $X\left(z_{l}(t),\mathbb{C}_{l}\right)$
over the time interval $[t_{0},\infty).$ Suppose $X\left(z_{l}(t_{0}),\mathbb{C}_{l}\right)=z_{l_{0}}.$
Suppose $X\left(z_{l}(t),\mathbb{C}_{l}\right)$ has taken the value
$z_{l_{1}}$during $t\in[t_{0},t_{1}],$ and $X\left(z_{l}(t),\mathbb{C}_{l}\right)$
has taken the value $z_{l_{2}}$during $t\in(t_{1},t_{2}].$ Here
$z_{l_{1}}\in D\left(z_{l_{0}}.r_{0}\right)\subset\mathbb{C}_{l}$
for $r_{0}>0,$ and $z_{l_{2}}\in D\left(z_{l_{1}}.r_{1}\right)\subset\mathbb{C}_{l}$
for $r_{1}>0.$ By this construction, we have
\[
D\left(z_{l_{0}}.r_{0}\right)\cap D\left(z_{l_{1}}.r_{1}\right)\neq\phi\text{ (empty set)}
\]
Since (\ref{def:Distinct-complex-numbers}) holds, we have $z_{l_{1}}\neq z_{l_{0}}.$
The transition probability for $X\left(z_{l}(t),\mathbb{C}_{l}\right)$
from $z_{l_{0}}$to $z_{l_{2}}$is 

\[
p(z_{l_{0}},z_{l_{2}},t)\text{}\left(t\in[t_{0},t_{2}]\right)=p(z_{l_{0}},z_{l_{1}},t)\text{}\left(t\in[t_{0},t_{1}]\right)\text{ and }p(z_{l_{1}},z_{l_{2}},t)\text{}\left(t\in(t_{1},t_{2}]\right)
\]
\begin{align}
=\text{Prob}[X\left(z_{l}(t),\mathbb{C}_{l}\right)=z_{l_{1}}\text{}\left(t\in[t_{0},t_{1}]\right)\left/X\left(z_{l}(t_{0}),\mathbb{C}_{l}\right)=z_{l_{0}}]\times\right.\nonumber \\
\text{ }\text{Prob}[X\left(z_{l}(t),\mathbb{C}_{l}\right)=z_{l_{2}}\text{}\left(t\in(t_{1},t_{2}]\right)\left/X\left(z_{l}(t),\mathbb{C}_{l}\right)=z_{l_{1}}\text{}\left(t\in[t_{0},t_{1}]\right)]\right.\nonumber \\
\label{eq:MCin=00005Bt0,t2=00005D}
\end{align}
Through (\ref{eq:MCin=00005Bt0,t2=00005D}), we can conclude that
the random variable $X\left(z_{l}(t),\mathbb{C}_{l}\right)$ $\text{}\left(t\in[t_{0},t_{2}]\right)$
obeys Markov property during the interval $[t_{0},t_{2}].$ In (\ref{eq:MCin=00005Bt0,t2=00005D}),
the number $z_{l_{2}}$is generated within a disc around the number
$z_{l_{1}}$but not around the disc with center $z_{l_{0}}.$ A contour
is drawn from $z_{l_{0}}$ to $z_{l_{2}}$only through $z_{l_{1}}.$
In a similar way, the value of $X\left(z_{l}(t),\mathbb{C}_{l}\right)\text{}\left(t\in(t_{n-1},t_{n}]\right)$
is located in the disc $D\left(z_{l_{n-1}}.r_{n-1}\right)$ for $r_{n-1}>0$
and not on the discs $D\left(z_{l_{k}}.r_{k}\right)$ for $r_{k}>0$
and $k=0,1,...,n-2.$ That is, 
\[
\text{Prob}[X\left(z_{l}(t),\mathbb{C}_{l}\right)=z_{l_{n}}\text{}\left(t\in(t_{n-1},t_{n}]\right)\left/X\left(z_{l}(t_{0}),\mathbb{C}_{l}\right)=z_{l_{0}},X\left(z_{l}(t),\mathbb{C}_{l}\right)=z_{l_{1}}\right.
\]
\begin{align}
X\left(z_{l}(t),\mathbb{C}_{l}\right) & =z_{l_{1}}\text{}\left(t\in(t_{0},t_{1}]\right),...,...,X\left(z_{l}(t),\mathbb{C}_{l}\right)=z_{l_{n-1}}\text{}\left(t\in(t_{n-2},t_{n-1}]\right)]\nonumber \\
 & =\text{Prob}[X\left(z_{l}(t),\mathbb{C}_{l}\right)=z_{l_{n}}\text{}\left(t\in(t_{n-1},t_{n}]\right)\left/\right.\nonumber \\
 & X\left(z_{l}(t),\mathbb{C}_{l}\right)=z_{l_{n-1}}\text{}\left(t\in(t_{n-2},t_{n-1}]\right)]\label{eq:MCt0,tn}
\end{align}
\\
A contour is drawn from $z_{l_{0}}$ to $z_{l_{n}}$only connecting
the numbers (points) through $z_{l_{1}},z_{l_{2}},...,$ $z_{l_{n-1}}.$
The result in (\ref{eq:MCt0,tn}) is also true when $t_{n}\rightarrow\infty.$
Hence the contour $\gamma_{l}$ is formed using the numbers generated
by $X\left(z_{l}(t),\mathbb{C}_{l}\right)$ $\text{}\left(t\in[t_{0},\infty]\right)$
obeys properties of continuous time Markov property or obeys continuous
time Markov chain.
\end{proof}

\subsection{Behavior of $X\left(z_{l}(t),\mathbb{C}_{l}\right)$ at $\left(\mathbb{C}_{l}\cap\mathbb{C}_{0}\right).$ }

Let us understand the behavior of $X\left(z_{l}(t),\mathbb{C}_{l}\right)$
at the intersection of $\mathbb{C}_{l}\cap\mathbb{C}_{0}.$ Set $\mathbb{C}_{l}\cap\mathbb{C}_{0}$
consists of elements $\mathbb{C}_{l}\cap\mathbb{C}_{0}=\left\{ z:z\in\mathbb{C}_{l}\text{ and }z\in\mathbb{C}_{0}\right\} .$
Here $\mathbb{C}_{l}\cap\mathbb{C}_{0}\neq\phi\text{ (empty).}$ The
numbers on $\mathbb{C}_{l}\cap\mathbb{C}_{0}$ form a dense set of
numbers and they form a line. Whenever $X\left(z_{l}(t),\mathbb{C}_{l}\right)$
reaches a number, say, $z_{l_{a}}$ in the set $\mathbb{C}_{l}\cap\mathbb{C}_{0}$,
then $X\left(z_{l}(t),\mathbb{C}_{l}\right)$ can choose the next
number within the set $\mathbb{C}_{l}\cap\mathbb{C}_{0}$ or within
the plane $\mathbb{C}_{l}$ or within the plane $\mathbb{C}_{0}.$
That is, as soon as $X\left(z_{l}(t),\mathbb{C}_{l}\right)$ chooses
a number $z_{l_{a}}$after a countable number of transitions from
$z_{l_{0}},$ then the two-step randomness will help to form two discs,
namely, $D\left(z_{l_{a}},r_{a},\mathbb{C}_{l}\right)$ and $D\left(z_{l_{a}},r_{a},\mathbb{C}_{0}\right).$
Here
\begin{equation}
D\left(z_{l_{a}},r_{a},\mathbb{C}_{l}\right)\subset\mathbb{C}_{l}\text{ and }D\left(z_{l_{a}},r_{a},\mathbb{C}_{0}\right)\subset\mathbb{C}_{0}.\label{eq:TWODISCS at Intersection}
\end{equation}
One can have different radii for two discs. For both the discs we
leave the the symbol $z_{l_{a}}$ in both the discs to indicate that
the associated random variable has origins in $\mathbb{C}_{l}.$ The
next number of $X\left(z_{l}(t),\mathbb{C}_{l}\right)$, say, $z_{l_{b}}$
could fall within $D\left(z_{l_{a}},r_{a},\mathbb{C}_{l}\right)$
such that it is on line $\mathbb{C}_{l}\cap\mathbb{C}_{0}$ or it
is in the set $\mathbb{C}_{l}\cap\left(\mathbb{C}_{l}\cap\mathbb{C}_{0}\right)'$.
Alternatively, $z_{l_{b}}$ could fall within $D\left(z_{l_{a}},r_{a},\mathbb{C}_{0}\right)$
such that it is on line $\mathbb{C}_{l}\cap\mathbb{C}_{0}$ or it
is in the set $\mathbb{C}_{0}\cap\left(\mathbb{C}_{l}\cap\mathbb{C}_{0}\right)'$where
\begin{align}
\mathbb{C}_{l}\cap\left(\mathbb{C}_{l}\cap\mathbb{C}_{0}\right)' & =\left\{ z:z\in\mathbb{C}_{l}\text{ and }z\notin\left(\mathbb{C}_{l}\cap\mathbb{C}_{0}\right)\right\} \nonumber \\
\mathbb{C}_{0}\cap\left(\mathbb{C}_{l}\cap\mathbb{C}_{0}\right)' & =\left\{ z:z\in\mathbb{C}_{0}\text{ and }z\notin\left(\mathbb{C}_{l}\cap\mathbb{C}_{0}\right)\right\} .\label{eq:TEOSETS at intersection}
\end{align}

If $z_{l_{b}}\in\mathbb{C}_{0}\cap\left(\mathbb{C}_{l}\cap\mathbb{C}_{0}\right)'$,
then by joining all the numbers generated by $X\left(z_{l}(t),\mathbb{C}_{l}\right)$
from$z_{l_{0}}$ through $z_{l_{b}},$ we form a multilevel contour.
If $z_{l_{b}}\in\mathbb{C}_{l}\cap\left(\mathbb{C}_{l}\cap\mathbb{C}_{0}\right)'$,
then by joining all the numbers generated by $X\left(z_{l}(t),\mathbb{C}_{l}\right)$
from$z_{l_{0}}$ through $z_{l_{b}},$ we form a contour on $\mathbb{C}_{l}$.
Even if $z_{l_{b}}\in\mathbb{C}_{l}\cap\left(\mathbb{C}_{l}\cap\mathbb{C}_{0}\right)'$,
still $X\left(z_{l}(t),\mathbb{C}_{l}\right)$ can choose at a future
iteration a value (number) in the set $\left(\mathbb{C}_{l}\cap\mathbb{C}_{0}\right)$
and escape the plane through $\mathbb{C}_{0}.$ Again at a future
iteration the value of $X\left(z_{l}(t),\mathbb{C}_{l}\right)$ could
fall within $\mathbb{C}_{l}.$ As soon as $X\left(z_{l}(t),\mathbb{C}_{l}\right)$
leaves $\mathbb{C}_{l}$ (if in case) it reaches another plane, say,
$\mathbb{C}_{p},$ because $\left(\mathbb{C}_{p}\cap\mathbb{C}_{0}\right)\neq\phi$
and $\left(\mathbb{C}_{p}\cap\mathbb{C}_{0}\right)\subset\mathbb{C}_{p}$
for some elements of $\mathbb{C}_{0}.$ Anytime $X\left(z_{l}(t),\mathbb{C}_{l}\right)$
reaches a set of intersecting planes $\left(\mathbb{C}_{l}\cap\mathbb{C}_{0}\right)$
or $\left(\mathbb{C}_{p}\cap\mathbb{C}_{0}\right)$ or some other
similar intersecting planes, it will have the power to generate next
number in two distinct discs as described above in (). This feature
of $X\left(z_{l}(t),\mathbb{C}_{l}\right)$ helps to form multilevel
contours. This feature is summarized below:
\begin{align}
X\left(z_{l}(t),\mathbb{C}_{l}\right)=z_{l_{b}}\text{ and }z_{l_{b}}\in D\left(z_{l_{a}},r_{a},\mathbb{C}_{0}\right) & \implies\gamma_{l}\text{ is a multilevel contour}\nonumber \\
X\left(z_{l}(t),\mathbb{C}_{l}\right)=z_{l_{b}}\text{ and }z_{l_{b}}\in D\left(z_{l_{a}},r_{a},\mathbb{C}_{l}\right) & \implies\gamma_{l}\text{ is a contour on }\mathbb{C}_{l}\nonumber \\
\label{eq:intersectionRULES}
\end{align}

\begin{figure}
\includegraphics[scale=0.8]{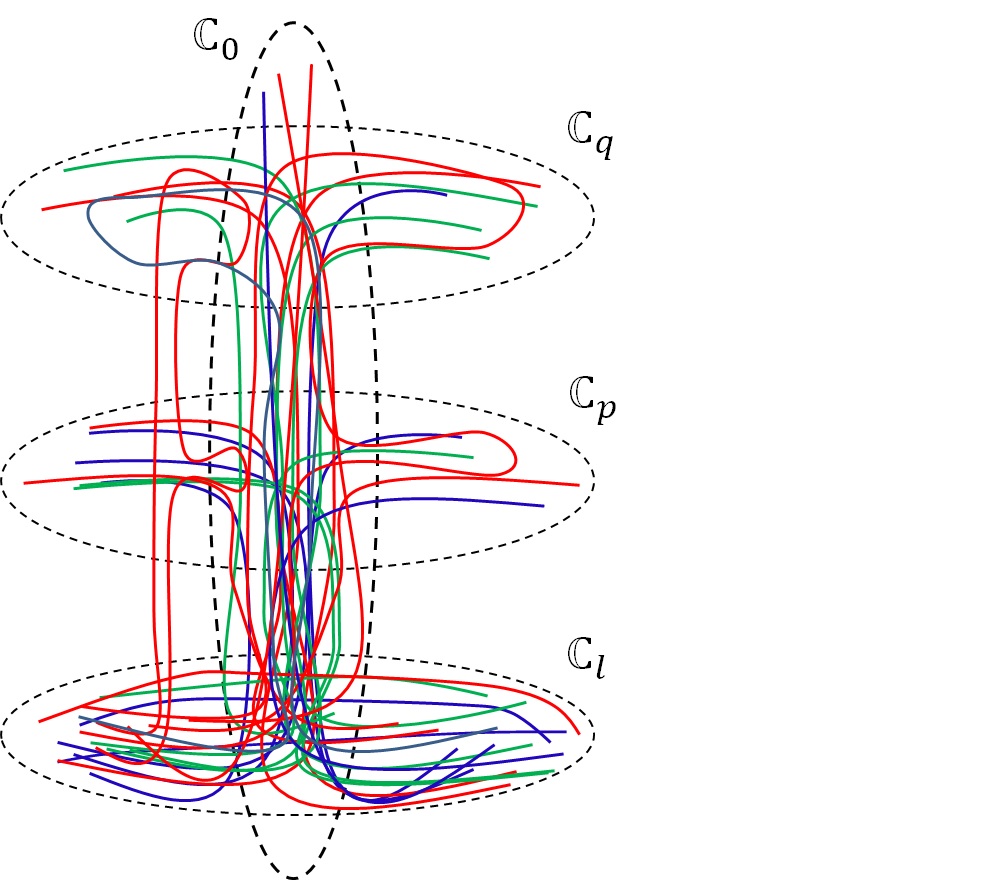}

\caption{\label{fig:Contoursonmultipleplanes}Contours spreading across one
or more planes from $\mathbb{C}_{l}$ due to random environment.}

\end{figure}

This rule (\ref{eq:intersectionRULES}) is applicable to the each
time the value of $X\left(z_{l}(t),\mathbb{C}_{l}\right)$ falls in
an intersection of planes. Once a contour attains the multilevel contour
property it will remain as a multilevel contour of that particular
$X\left(z_{l}(t),\mathbb{C}_{l}\right)$ even if the value of $X\left(z_{l}(t),\mathbb{C}_{l}\right)$
reruns and remains in $\mathbb{C}_{l}$ forever. The value of the
radius at each two-step randomness and the location of the next number
to be picked by $X\left(z_{l}(t),\mathbb{C}_{l}\right)$ decides the
time taken for a contour to become a multilevel contour (if there
is a possibility to become). See Figure \ref{fig:Contoursonmultipleplanes}.The
time interval to reach $z_{l_{a}}$ from $z_{l_{0}}$ could be at
least one, and it requires at least two time interval to reach $z_{l_{b}}$
from $z_{l_{0}}$ under the framework described above. 

Suppose $z_{l_{b}}\in D\left(z_{l_{a}},r_{a},\mathbb{C}_{0}\right)'$.
If $z_{l_{0}}$ to $z_{l_{a}}$ reaches in one time interval and $z_{l_{a}}$
to $z_{l_{b}}$ reaches in one time interval then 

\[
\left\Vert z_{l_{b}}-z_{l_{0}}\right\Vert >\left\Vert z_{l_{a}}-z_{l_{0}}\right\Vert 
\]
because 
\[
\left\Vert z_{l_{a}}-z_{l_{0}}\right\Vert +\left\Vert z_{l_{a}}-z_{l_{b}}\right\Vert =\left\Vert z_{l_{b}}-z_{l_{0}}\right\Vert 
\]
If $z_{l_{0}}$ to $z_{l_{a}}$ reaches in more than one time interval,
then the length of the contour, $L(z_{l_{o}},z_{l_{a}})$ from $z_{l_{0}}$to
$z_{l_{a}}$ is still less than $L(z_{l_{o}},z_{l_{b}})$ from $z_{l_{0}}$to
$z_{l_{b}}$ because $z_{l_{b}}$ lies in the disc $D\left(z_{l_{a}},r_{a},\mathbb{C}_{0}\right)'.$
Suppose $X\left(z_{l}(t_{0}),\mathbb{C}_{l}\right)=z_{l_{0}}$ and
$X\left(z_{l}(t),\mathbb{C}_{l}\right)$ reaches at $z_{l_{a}}$during
the $n^{th}$ time interval. Let $X\left(z_{l}(t),\mathbb{C}_{l}\right)=z_{l_{1}}$
for $z_{l_{1}}\in D\left(z_{l_{0}}.r_{0}\right)$ during the first
time interval $[t_{0},t_{1}]$ and $X\left(z_{l}(t),\mathbb{C}_{l}\right)=z_{l_{2}}$
for $z_{l_{2}}\in D\left(z_{l_{1}}.r_{1}\right)$ during the second
time interval and so on $X\left(z_{l}(t),\mathbb{C}_{l}\right)=z_{l_{a}}$
for $z_{l_{a}}\in D\left(z_{l_{n-1}}.r_{n-1}\right)$ during $n^{th}$
time interval. Suppose $X\left(z_{l}(t),\mathbb{C}_{l}\right)=z_{l_{b}}$
for $z_{l_{b}}\in D\left(z_{l_{a}}.r_{a}\right).$ Suppose $\gamma_{l}$
is described by 
\[
\gamma_{l}=z_{l}(t)\text{ for }t\in[t_{0},\infty)
\]
and the parametric representations are given by 
\[
\begin{array}{cc}
t= & \left\{ \begin{array}{c}
v_{l_{1}}(\tau)\text{ }\left(a_{0}^{l_{0}}\leq t\leq a_{0}^{l_{1}}\right)\\
v_{l_{2}}(\tau)\text{ }\left(a_{0}^{l_{1}}<t\leq a_{0}^{l_{2}}\right)\\
\vdots\\
v_{l_{n-1}}(\tau)\text{ }\left(a_{0}^{l_{n-2}}<t\leq a_{0}^{l_{n-1}}\right)\\
v_{l_{a}}(\tau)\text{ }\left(a_{0}^{l_{n-1}}<t\leq a_{0}^{l_{a}}\right)\\
v_{l_{b}}(\tau)\text{ }\left(a_{0}^{l_{a}}<t\leq a_{0}^{l_{b}}\right)\\
\\
\end{array}\right.\end{array},
\]
and the real valued functions $v_{l_{i+1}}$ mapping $(a_{0}^{l_{i}},a_{0}^{l_{i+1}}]$
for $i=1,...,n-2$ onto the intervals $[t_{0},t_{1}],$ $(t_{1},t_{2}],$
$(t_{n-2},t_{n-1}].$ The real valued function $v_{l_{0}}$ maps $[a_{0}^{l_{0}},a_{0}^{l_{1}}]$
onto the interval $([t_{0},t_{1}]$, $v_{l_{a}}$ maps $(a_{0}^{l_{n-1}},a_{0}^{l_{a}}]$
onto the interval $(t_{n-1},t_{a}]$ and the real valued function
$v_{l_{b}}$ maps $(a_{0}^{l_{a}},a_{0}^{l_{b}}]$ onto the interval
$(t_{a},t_{b}].$ Then
\begin{align}
L(z_{l_{o}},z_{l_{b}}) & =\sum_{i=1}^{n-1}\int_{a_{0}^{l_{i}}}^{a_{0}^{l_{i+1}}}\left|z\left[v_{l_{i}}(\tau)\right]\right|v'_{l_{i}}(\tau)d\tau+\int_{a_{0}^{l_{n-1}}}^{a_{0}^{l_{a}}}\left|z\left[v_{l_{a}}(\tau)\right]\right|v'_{l_{a}}(\tau)d\tau\nonumber \\
 & +\int_{a_{0}^{l_{a}}}^{a_{0}^{l_{b}}}\left|z\left[v_{l_{b}}(\tau)\right]\right|v'_{l_{b}}(\tau)d\tau\label{eq:Lz0zb>Lz0za}\\
 & >L(z_{l_{o}},z_{l_{a}})\nonumber 
\end{align}
because the first two terms of the R.H.S. of (\ref{eq:Lz0zb>Lz0za})
is $L(z_{l_{o}},z_{l_{a}}).$ Here

\begin{align*}
\int_{a_{0}^{l_{i}}}^{a_{0}^{l_{i+1}}}\left|z\left[v_{l_{i}}(\tau)\right]\right|v'_{l_{i}}(\tau)d\tau & \subset D\left(z_{l_{i}},r_{i}\right)\text{ for }i=1,2,...,n-1\\
\\
\int_{a_{0}^{l_{n-1}}}^{a_{0}^{l_{a}}}\left|z\left[v_{l_{a}}(\tau)\right]\right|v'_{l_{a}}(\tau)d\tau & \subset D\left(z_{l_{n-1}},r_{a},\mathbb{C}_{l}\right)\\
\\
\int_{a_{0}^{l_{a}}}^{a_{0}^{l_{b}}}\left|z\left[v_{l_{b}}(\tau)\right]\right|v'_{l_{b}}(\tau)d\tau & \subset D\left(z_{l_{a}},r_{a},\mathbb{C}_{0}\right)
\end{align*}
and
\[
D\left(z_{l_{0}},r_{0}\right)\cap D\left(z_{l_{1}},r_{1}\right)\cap...\cap D\left(z_{l_{n-2}},r_{n-2}\right)\neq\phi\text{ (empty)}
\]
each of these discs are non-empty and they have distinct set of numbers
on $\mathbb{C}_{l}.$ The disc $D\left(z_{l_{a}},r_{a},\mathbb{C}_{0}\right)$
has some elements outside the plane $\mathbb{C}_{l}$, and

\[
D\left(z_{l_{n-1}},r_{n-1},\mathbb{C}_{l}\right)\cap D\left(z_{l_{a}},r_{a},\mathbb{C}_{0}\right)\neq\phi\text{ (empty)}.
\]
Suppose it takes infinitely many time intervals to reach $z_{l_{a}}$
from $z_{l_{0}}$ (due to the random environment created). 

Extending the parametric representation described above, the length
of the contour from $z_{l_{0}}$ to $z_{l_{a}}$is
\[
L(z_{l_{o}},z_{l_{a}})=\sum_{i=1}^{\infty}\int_{a_{0}^{l_{i}}}^{a_{0}^{l_{\infty}}}\left|z\left[v_{l_{i}}(\tau)\right]\right|v'_{l_{i}}(\tau)d\tau<L(z_{l_{o}},z_{l_{b}})
\]
because
\[
L(z_{l_{o}},z_{l_{b}})=\sum_{i=1}^{\infty}\int_{a_{0}^{l_{i}}}^{a_{0}^{l_{\infty}}}\left|z\left[v_{l_{i}}(\tau)\right]\right|v'_{l_{i}}(\tau)d\tau+\int_{a_{0}^{l_{\infty}}}^{a_{0}^{l_{b}}}\left|z\left[v_{l_{b}}(\tau)\right]\right|v'_{l_{b}}(\tau)d\tau
\]
and
\[
\int_{a_{0}^{l_{\infty}}}^{a_{0}^{l_{b}}}\left|z\left[v_{l_{b}}(\tau)\right]\right|v'_{l_{b}}(\tau)d\tau\subset D\left(z_{l_{\infty}},r_{\infty},\mathbb{C}_{0}\right)
\]
\begin{equation}
\bigcap_{i=0}^{\infty}D\left(z_{l_{i}},r_{i}\right)\neq\phi\text{ (empty)}.\label{eq:INFINITE-DISCS}
\end{equation}

\begin{thm}
\label{thm:Nested discs}Suppose it takes infinitely many time intervals
for $X\left(z_{l}(t),\mathbb{C}_{l}\right)$ to reach $z_{l_{a}}$
from $z_{l_{0}}$ The infinitely many discs (uncountable) created
while reaching $z_{l_{a}}$from $z_{l_{0}}$ could be nested under
the two-step random environment created by $X\left(z_{l}(t),\mathbb{C}_{l}\right)$
and $\bigcap_{i=0}^{\infty}D\left(z_{l_{i}},r_{i}\right)\neq\phi.$
\end{thm}

\begin{proof}
Suppose $X\left(z_{l}(t_{0}),\mathbb{C}_{l}\right)=z_{l_{0}}.$ Let
$D\left(z_{l_{0}},r_{0}\right)$ formed out of two-step randomness
and $z_{l_{1}}\in D\left(z_{l_{0}},r_{0}\right).$ Suppose $r_{1}$
is generated randomly such that $D\left(z_{l_{1}},r_{1}\right)$$\subset D\left(z_{l_{0}},r_{0}\right)$.
Further, let $r_{i}$ is generated randomly such that 
\[
D\left(z_{l_{i}},r_{i}\right)\subset D\left(z_{l_{i-1}},r_{i-1}\right)\text{ for }i=2,3,...,\infty
\]
and
\[
z_{l_{i}}\in D\left(z_{l_{i-1}},r_{i-1}\right)\text{ for }i=1,2,...,\infty
\]
See Figure \ref{fig:Nested-discs-and}. Given that $z_{l_{a}}$has
been generated by $X\left(z_{l}(t),\mathbb{C}_{l}\right)$, it must
be in one of the infinitely many discs formed described earlier before
the Theorem \ref{thm:Nested discs}. Moreover, $z_{l_{a}}\in\mathbb{C}_{l}\cap\mathbb{C}_{0}.$
So $\bigcap_{i=0}^{\infty}D\left(z_{l_{i}},r_{i}\right)\neq\phi.$

\begin{figure}
\includegraphics[scale=0.45]{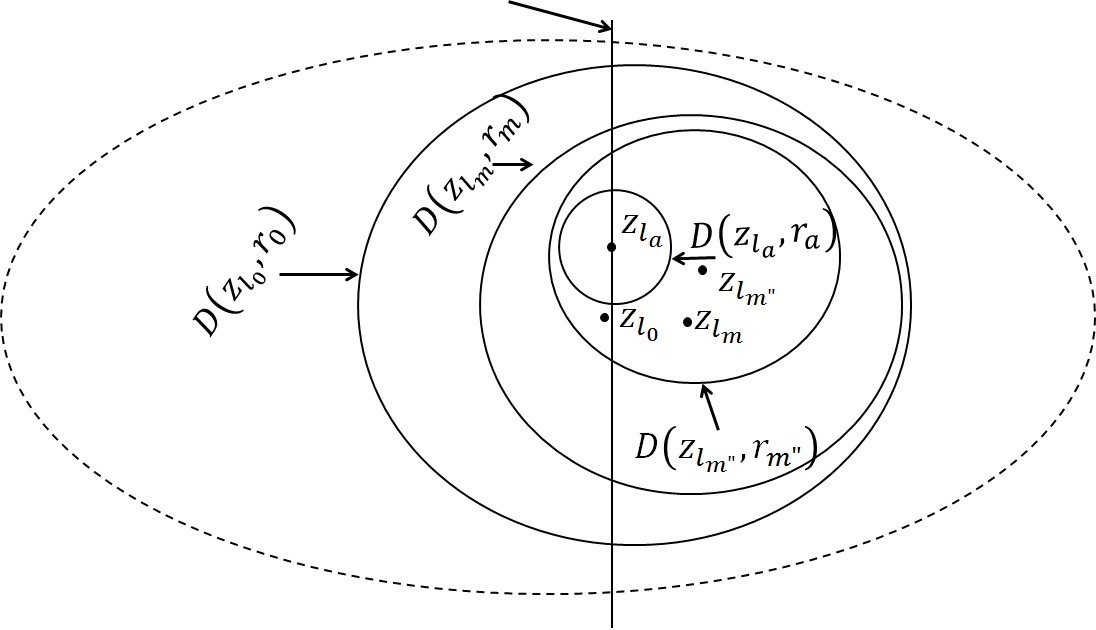}

\caption{\label{fig:Nested-discs-and}Nested discs and points on the intersection
of planes.}

\end{figure}
\end{proof}
\begin{rem}
As a consequence of Theorem \ref{thm:Nested discs}, we will see that
\[
\sum_{i=1}^{\infty}\int_{a_{0}^{l_{i}}}^{a_{0}^{l_{\infty}}}\left|z\left[v_{l_{i}}(\tau)\right]\right|v'_{l_{i}}(\tau)d\tau\subset D\left(z_{l_{0}},r_{0}\right).
\]
\end{rem}

The random fluctuations and transitions explained in this sub-section
could arise infinitely many times during $[t_{0},\infty).$ The random
environment created on $\mathbb{C}_{l}$ and its behavior at the intersection
of $\mathbb{C}_{l}$ and $\mathbb{C}_{0}$ is crucial for the creation
of multilevel contours. A contour $\gamma_{l}$ could remain forever
in $\mathbb{C}_{l}$ or could go beyond $\mathbb{C}_{l}$ and reach
other planes. Contour $\gamma_{l}$ also has the potential to reach
infinitely many planes and also has the potential to return to each
and every plane infinitely many times. 
\begin{thm}
There exists a unique contour $\gamma_{l}$ for each initial number
chosen on $\mathbb{C}_{l}$ for $\mathbb{C}_{l}$ in the bundle $B_{\mathbb{R}}(\mathbb{C})$
combined with the $\mathbb{C}_{0}$ intersecting with the bundle and
satisfying the property (\ref{def:Distinct-complex-numbers}).
\end{thm}

\begin{proof}
Let $\gamma_{l}$ be a contour described by $z_{l}(t)$ for $t\in[t_{0},\infty)$
that has a starting point on the plane $\mathbb{C}_{l}$ with $z_{l}(t_{0})=z_{l_{0}}$
(say). The values of $z_{l}(t)$ can reach other planes because the
bundle $B_{\mathbb{R}}(\mathbb{C})$ in which $\mathbb{C}_{l}$ lies
is intersecting with the $\mathbb{C}_{0}.$ Due to the property (\ref{def:Distinct-complex-numbers}),
$X\left(z_{l}(t),\mathbb{C}_{l}\right)$ keeps on generating new numbers
during $[t_{0},\infty)$ and $z_{l_{a}}\neq z_{l_{b}}$ for $a\neq b$
and $a,b\in\mathbb{R}.$ At some of time $\gamma_{l}$ could become
a multilevel contour if the value generated by $X\left(z_{l}(t),\mathbb{C}_{l}\right)$
reaches another plane, say $\mathbb{C}_{p}$ for $\mathbb{C}_{p}\neq\mathbb{C}_{l}$.
Reaching of another plane is possible due to the positioning of $\mathbb{C}_{0}.$
Or $\gamma_{l}$ could remain as a contour on $\mathbb{C}_{l}$ for
$t\in[t_{0},\infty).$ 

In Theorem (\ref{thm:A-contour-isMarkovChain}), we saw that a contour
drawn from an initial number $z_{l_{0}}$ reaches $z_{l_{n}}$only
through distinct numbers $z_{l_{1}},$ $z_{l_{2}},$ $...,z_{l_{n-1}}.$
The numbers $z_{l_{1}},$ $z_{l_{2}},$ $...,z_{l_{n-1}}$ were generated
by $X\left(z_{l}(t),\mathbb{C}_{l}\right)$ . Alternatively, suppose
the initial value chosen by $X\left(z_{l}(t),\mathbb{C}_{l}\right)$
is say, $z'_{l_{0}}$ for $z'_{l_{0}}\neq z_{l_{0}},$then even if
the numbers $z_{l_{1}},$ $z_{l_{2}},$ $...,z_{l_{n-1}}$ are the
same through which a contour up to $z_{l_{n}}$is drawn (due to the
randomness in selection of complex numbers by $X\left(z_{l}(t),\mathbb{C}_{l}\right)$),
a contour drawn from $z_{l_{0}}$to $z_{l_{n}}$ would be different
than $z'_{l_{0}}.$ This argument is valid for $\gamma_{l}$ a contour
on $\mathbb{C}_{l}$ or a multilevel contour. Hence, as $t_{n}\rightarrow\infty$
the numbers generated by $X\left(z_{l}(t),\mathbb{C}_{l}\right)$
would be distinct due to the property (\ref{def:Distinct-complex-numbers}).
Therefore $\gamma_{l}$ is unique. 
\end{proof}
Each point (number) on $\mathbb{C}_{l}$ has potentially capable to
produce a contour which could remain forever in $\mathbb{C}_{l}$
or could become a multilevel contour by crossing $\mathbb{C}_{l}$
through $\mathbb{C}_{0}.$$X\left(z_{l}(t),\mathbb{C}_{l}\right)$
has the power to choose the first number of $\mathbb{C}_{l}$and generate
rest of the numbers randomly. Due to the property (\ref{def:Distinct-complex-numbers}),
$X\left(z_{l}(t),\mathbb{C}_{l}\right)$ looses the power to choose
a number that was already been chosen earlier. That is

\[
p(z_{l_{a}},z_{l_{a}},t)\text{ }(t\in[t_{0},\infty))=0
\]
where $p(z_{l_{a}},z_{l_{a}},t)$ is the probability of transition
to the same number is zero. Suppose $z_{l_{0}},$ $z_{l_{1}}$, $z_{l_{2}},...,z_{l_{n}},...$
are the set of numbers generated by $X\left(z_{l}(t),\mathbb{C}_{l}\right)$
over the time. Let $p(z_{l_{0}},z_{l_{1}},t)^{(1)}$ represent the
probability of reaching $z_{l_{1}}$from $z_{l_{0}}$in $1-$time
interval $[t_{0},t_{1}]$, $p(z_{l_{0}},z_{l_{1}},t)^{(2)}$ represent
the probability of reaching $z_{l_{2}}$from $z_{l_{0}}$in $2-$time
intervals $[t_{0},t_{1}]$, $(t_{1},t_{2}]$, and so on. We note that
$z_{l_{2}}$cannot be reached directly from $z_{l_{0}}$ in $1-$time
interval $[t_{0},t_{1}]$ or $[t_{0},t_{2}]$ as mentioned previously
in the article. So 
\begin{equation}
p(z_{l_{0}},z_{l_{2}},t)^{(2)}>0\text{ and }p(z_{l_{0}},z_{l_{2}},t)^{(1)}=0,\label{eq:highertransitionszero}
\end{equation}
and
\begin{equation}
p(z_{l_{1}},z_{l_{3}},t)^{(2)}>0\text{ and }p(z_{l_{1}},z_{l_{3}},t)^{(m)}=0,\text{ for }m=1,3,4,5,...\label{eq:higherordertransitions}
\end{equation}
The $n-$time intervals transition probabilities between any other
two distinct complex numbers on $\mathbb{C}_{l}$ can be expressed
as in (\ref{eq:highertransitionszero}) and (\ref{eq:higherordertransitions}).
\begin{rem}
The notation $X\left(z_{l}(t),\mathbb{C}_{l}\right)$ $t\in[t_{0},\infty)$
is used even if $X\left(z_{l}(t),\mathbb{C}_{l}\right)$ starts generating
numbers from different plane after crossing through $\mathbb{C}_{0}$
after the initial value $z_{l_{0}}$ was chosen in $\mathbb{C}_{l}.$
Such notation will help identify the origin plane of $X\left(z_{l}(t),\mathbb{C}_{l}\right).$ 
\end{rem}

\begin{thm}
Given $X\left(z_{l}(t_{0}),\mathbb{C}_{l}\right)=z_{l_{0}}.$ Suppose
$p(z_{l_{a}},z_{l_{b}},t)^{(n)}$ represent $n-$time intervals transition
probabilities from $z_{l_{a}}$to $z_{l_{b}}$ where 
\begin{align*}
X\left(z_{l}(t),\mathbb{C}_{l}\right)=z_{l_{a}} & \text{ for }t\in(t_{a-1},t_{a}]\,\text{and}\\
X\left(z_{l}(t),\mathbb{C}_{l}\right)=z_{l_{b}} & \text{ for }t\in(t_{b-1},t_{b}],
\end{align*}
for $t_{b}>t_{a}$, then 
\[
p(z_{l_{a}},z_{l_{b}},t)^{(n)}=\left\{ \begin{array}{c}
>0\text{ if }n=t_{b}-t_{a\text{ }}\text{ sequential time intervals}\\
=0\text{ if }n\neq t_{b}-t_{a}\text{ sequential time intervals}
\end{array}\right..
\]
\end{thm}

\begin{proof}
The $n-$time interval transition probability $p(z_{l_{a}},z_{l_{b}},t)^{(n)}$
is written as 
\begin{align}
p(z_{l}(t_{a}),z_{l}(t_{b}),t)^{(n)}\text{ }\left(t\in(t_{a},t_{b}]\right) & =p(z_{l_{a}},z_{l_{a_{1}}},t)^{(1)}\text{ }\left(t\in(t_{a},t_{a_{1}}]\right).\nonumber \\
\nonumber \\
 & p(z_{l_{a_{1}}},z_{l_{a_{2}}},t)^{(1)}\text{ }\left(t\in(t_{a_{1}},t_{a_{2}}]\right).\nonumber \\
\nonumber \\
 & ....p(z_{l_{a_{n-1}}},z_{l_{b}},t)^{(1)}\text{ }\left(t\in(t_{a_{n-1}},t_{b}]\right),\label{eq:n-stepProbinTh9}
\end{align}
for $t_{a}<t_{a_{1}}<...t_{a_{n-1}}<t_{b}.$ In (\ref{eq:n-stepProbinTh9}),
$z_{l_{a_{1}}}$ is generated by $X\left(z_{l}(t),\mathbb{C}_{l}\right)$
during $(t_{a},t_{a_{1}}]$ from the set of numbers of the disc $D(z_{l_{a}},r_{a})$,
$z_{l_{a_{2}}}$ is generated by $X\left(z_{l}(t),\mathbb{C}_{l}\right)$
during $(t_{a_{1}},t_{a_{2}}]$ from the set of numbers of the disc
$D(z_{l_{a_{1}}},r_{a_{1}})$, and so on, $z_{l_{b}}$ is generated
by $X\left(z_{l}(t),\mathbb{C}_{l}\right)$ during $(t_{a_{n-1}},t_{b}]$
from the set of numbers of the disc $D(z_{l_{a_{n-1}}},r_{a_{n-1}}).$
The numbers $z_{l_{a_{1}}},$ $z_{l_{a_{2}}}$, $...,$ $z_{l_{a_{n-1}}}$,
$z_{l_{b}}$were sequentially generated from the sets of distinct
discs 
\begin{equation}
\left\{ D(z_{l_{a}},r_{a}),D(z_{l_{a_{1}}},r_{a_{1}}),...,D(z_{l_{a_{n-1}}},r_{a_{n-1}})\right\} \label{eq:sequential DISCS}
\end{equation}
within the sequential time intervals
\begin{equation}
\left\{ (t_{a},t_{a_{1}}],(t_{a_{1}},t_{a_{2}}],...,(t_{a_{n-1}},t_{b}]\right\} \label{eq:sequencial-TIMESINTERVALS}
\end{equation}
Due to the sequential nature of (\ref{eq:sequential DISCS}) and (\ref{eq:sequencial-TIMESINTERVALS}),
we will have
\begin{align}
p(z_{l_{a}},z_{l_{a_{1}}},t)^{(1)}\text{ }\left(t\in(t_{a},t_{a_{1}}]\right) & >0\nonumber \\
p(z_{l_{a_{1}}},z_{l_{a_{2}}},t)^{(1)}\text{ }\left(t\in(t_{a_{1}},t_{a_{2}}]\right) & >0\nonumber \\
\vdots\nonumber \\
p(z_{l_{a_{n-1}}},z_{l_{b}},t)^{(1)}\text{ }\left(t\in(t_{a_{n-1}},t_{b}]\right) & >0.\label{eq:oneSTEP probabilities}
\end{align}
From (\ref{eq:oneSTEP probabilities}) we conclude that
\begin{equation}
p(z_{l}(t_{a}),z_{l}(t_{b}),t)^{(n)}\text{ }\left(t\in(t_{a},t_{b}]\right)>0.\label{eq:mainstatementTH9}
\end{equation}
We also note that,
\begin{align}
t_{a_{1}-t_{a}} & =1-\text{ step time interval}\nonumber \\
t_{a_{2}-t_{a_{1}}} & =1-\text{ step time interval}\nonumber \\
\vdots\nonumber \\
t_{b}-t_{a_{n-1}} & =1-\text{ step time interval}\label{eq:1-stepINTERVALS}
\end{align}
Summing up the terms of the L.H.S. of (\ref{eq:1-stepINTERVALS})
and equating it to sum of the quantities of the R.H.S. of (\ref{eq:1-stepINTERVALS}),
we see that 
\begin{equation}
t_{b}-t_{a}=n-\text{ sequential time intervals}.\label{eq:tb-ta=00003Dn}
\end{equation}
 From (\ref{eq:mainstatementTH9}) and (\ref{eq:tb-ta=00003Dn}) we
conclude that $p(z_{l_{a}},z_{l_{b}},t)^{(n)}>0$ when $t_{b}-t_{a}=n-\text{ sequential time intervals}.$

Suppose $t_{b}-t_{a}\neq m-\text{ sequential time intervals}.$ This
implies $z_{l_{b}}$is generated either through less than $m-$ sequential
time intervals after choosing $z_{l_{a}}$ by $X\left(z_{l}(t),\mathbb{C}_{l}\right)$
or $z_{l_{b}}$is generated through more than $m-$ sequential time
intervals after choosing $z_{l_{a}}$ by $X\left(z_{l}(t),\mathbb{C}_{l}\right).$
If $m=1,$then

\[
t_{b}-t_{a}\neq1,\text{which implies,}p(z_{l_{a}},z_{l_{b}},t)^{(1)}=0
\]
if $m=2,$then

\[
t_{b}-t_{a}\neq2,\text{which implies,}p(z_{l_{a}},z_{l_{b}},t)^{(2)}=0
\]
and so on if $m=n,$then 
\[
t_{b}-t_{a}\neq n,\text{which implies,}p(z_{l_{a}},z_{l_{b}},t)^{(n)}=0.
\]
This concludes that $p(z_{l_{a}},z_{l_{b}},t)^{(n)}=0$ if $n\neq t_{b}-t_{a}\text{ sequential time intervals}.$

Each number in $\mathbb{C}_{l}$ that was chosen initially by $X\left(z_{l}(t),\mathbb{C}_{l}\right)$
and further points chosen by $X\left(z_{l}(t),\mathbb{C}_{l}\right)$
over the time $[t_{0},\infty)$ forms the state space of $X\left(z_{l}(t),\mathbb{C}_{l}\right)$
and is given by
\begin{align*}
A_{X}(\mathbb{C}_{l})= & \left\{ z:z=z_{l_{0}}\text{ at }X\left(z_{l}(t_{0}),\mathbb{C}_{l}\right)\text{ and \ensuremath{X\left(z_{l}(t),\mathbb{C}_{l}\right)\left(t_{0},\infty\right)=z} }\right.\\
 & \text{for z\ensuremath{\in B_{\mathbb{R}}}(\ensuremath{\mathbb{C}}).}\\
 & \text{Here the order of choosing }z\text{ is preserved and }z_{l_{a}}\neq z_{l_{b}}\\
 & \left.\text{ if }a\neq b,\text{ }a,b>t_{0}\right\} 
\end{align*}
As the time progresses it will keep on reaching new numbers in $B_{\mathbb{R}}(\mathbb{C}).$
All the elements of the set $A_{X}(\mathbb{C}_{l})$ need not be in
$\mathbb{C}_{l}$ due to behavior of $X\left(z_{l}(t_{0}),\mathbb{C}_{l}\right)$
at $\mathbb{C}_{l}\cap\mathbb{C}_{0}$ explained previously. The points
of $A_{X}(\mathbb{C}_{l})$ could be forever in $\mathbb{C}_{l}$
or they could spread across one or more planes of $B_{\mathbb{R}}(\mathbb{C}).$
Each element in $A_{X}(\mathbb{C}_{l})$ is called a state of the
process $\left\{ X\left(z_{l}(t),\mathbb{C}_{l}\right)\right\} _{t\geq t_{0}}.$
If there is no certainty to return to a state after leave that state
by the random variable, we call it a transient state. In this case,
once $X\left(z_{l}(t),\mathbb{C}_{l}\right)$ chooses a state then
it won't be able to return to that state forever. Since there is an
equal probability for $X\left(z_{l}(t),\mathbb{C}_{l}\right)$ to
choose any complex number on $\mathbb{C}_{l}$, so there is a possibility
to form contours of infinite lengths starting from each point on plane
$\mathbb{C}_{l}.$ However, these infinite lengths of contours need
not be identical. Suppose we define another process $\left\{ Y\left(z_{l}(t),\mathbb{C}_{l}\right)\right\} _{t\geq t_{0}}$
with a two-step randomness that we used for generating the discs and
radii then we consider these two processes $\left\{ X\left(z_{l}(t),\mathbb{C}_{l}\right)\right\} _{t\geq t_{0}}$
and $\left\{ X\left(z_{l}(t),\mathbb{C}_{l}\right)\right\} _{t\geq t_{0}}$
are not disjoint to each other, that is, both these processes have
some probability to choose identical values in the same order or in
different order during $[t_{0},\infty).$ 
\end{proof}
\begin{thm}
Two contours formed during $[t_{0},\infty)$ need not be identical
but their lengths could be identical. 
\end{thm}

\begin{proof}
Let $\gamma_{l}(X)$ and $\gamma_{l}(Y)$ be two contours formed out
of the points created by two processes $\left\{ X\left(z_{l}(t),\mathbb{C}_{l}\right)\right\} _{t\geq t_{0}}$
and $\left\{ Y\left(z_{l}(t),\mathbb{C}_{l}\right)\right\} _{t\geq t_{0}}$
with $\left\{ X\left(z_{l}(t_{0}),\mathbb{C}_{l}\right)\right\} =z_{l_{0}}(X)$
and $\left\{ Y\left(z_{l}(t_{0}),\mathbb{C}_{l}\right)\right\} =z_{l_{0}}(Y).$
Let $\gamma_{l}(X)$ be described by $z_{X}(t)$ $[t_{0},\infty)$
and $\gamma_{l}(Y)$ be described by $z_{Y}(t)$ $[t_{0},\infty)$.
The two state spaces corresponding to the two processes are 
\begin{align*}
A_{X}(\mathbb{C}_{l}) & =\left\{ z:z=z_{l_{0}}(X)\text{ at }X\left(z_{l}(t_{0}),\mathbb{C}_{l}\right)\text{ and \ensuremath{X\left(z_{l}(t),\mathbb{C}_{l}\right)\left(t_{0},\infty\right)=z} }\right.\\
 & \text{for }z\in B_{\mathbb{R}}(\mathbb{C})\\
 & \text{with some order of choosing }z\text{ and }z_{l_{a}}(X)\neq z_{l_{b}}(X)\\
 & \left.\text{ if }a\neq b,\text{ }a,b>t_{0}\right\} 
\end{align*}

\begin{align*}
A_{Y}(\mathbb{C}_{l}) & =\left\{ z:z=z_{l_{0}}(Y)\text{ at }Y\left(z_{l}(t_{0}),\mathbb{C}_{l}\right)\text{ and \ensuremath{Y\left(z_{l}(t),\mathbb{C}_{l}\right)\left(t_{0},\infty\right)=z} }\right.\\
 & \text{for z\ensuremath{\in B_{\mathbb{R}}}(\ensuremath{\mathbb{C}}).}\\
 & \text{with some order of choosing }z\text{ and }z_{l_{a}}(Y)\neq z_{l_{b}}(Y)\\
 & \left.\text{ if }a\neq b,\text{ }a,b>t_{0}\right\} 
\end{align*}
Note, $\gamma_{l}(X)$ is identical to $\gamma_{l}(Y)$ if, and only
if, $z_{l_{0}}(X)=z_{l_{0}}(Y)$ and all other $z$ values generated
out of infinite iterations of two processes are identical, i.e. $A_{X}(\mathbb{C}_{l})=A_{Y}(\mathbb{C}_{l})$.
Since $A_{X}(\mathbb{C}_{l})$ and $A_{Y}(\mathbb{C}_{l})$ are not
disjoint, there is a possibility that $\left\{ X\left(z_{l}(t),\mathbb{C}_{l}\right)\right\} _{t\geq t_{0}}$
and $\left\{ Y\left(z_{l}(t),\mathbb{C}_{l}\right)\right\} _{t\geq t_{0}}$
may choose same numbers during $[t_{0},\infty).$ If $A_{X}(\mathbb{C}_{l})\neq A_{Y}(\mathbb{C}_{l})$,
then anyway $\gamma_{l}(X)$ is not identical to $\gamma_{l}(Y).$
Given that $A_{X}(\mathbb{C}_{l})$ and $A_{Y}(\mathbb{C}_{l})$ are
available, let $L\left[(z_{l_{0}}(X),z_{l_{a}}(X)\right]$ be the
length of the contour from $z_{l_{0}}(X)$ to $z_{l_{a}}(X)$ and
$L\left[(z_{l_{a}}(X),z_{l_{b}}(X)\right]$ be the length of the contour
from $z_{l_{a}}(X)$ to $z_{l_{b}}(X)$ such that the length of the
contour from $z_{l_{0}}(X)$ to $z_{l_{a}}(X)$ is computed as

\begin{equation}
L\left[(z_{l_{0}}(X),z_{l_{b}}(X)\right]=L\left[(z_{l_{0}}(X),z_{l_{a}}(X)\right]+L\left[(z_{l_{a}}(X),z_{l_{b}}(X)\right].\label{eq:Lz0za(X)}
\end{equation}
Let $L\left[(z_{l_{0}}(Y),z_{l_{a}}(Y)\right]$ be the length of the
contour from $z_{l_{0}}(Y)$ to $z_{l_{a}}(Y)$ and $L\left[(z_{l_{a}}(Y),z_{l_{b}}(Y)\right]$
be the length of the contour from $z_{l_{a}}(Y)$ to $z_{l_{b}}(Y)$
such that the length of the contour from $z_{l_{0}}(Y)$ to $z_{l_{a}}(Y)$
is computed as

\begin{equation}
L\left[(z_{l_{0}}(Y),z_{l_{b}}(Y)\right]=L\left[(z_{l_{0}}(Y),z_{l_{a}}(Y)\right]+L\left[(z_{l_{a}}(Y),z_{l_{b}}(Y)\right].\label{eq:Lz0zbY}
\end{equation}

In (\ref{eq:Lz0za(X)}) and (\ref{eq:Lz0zbY}), it is assumed that
$z_{l_{0}}(X)\neq z_{l_{0}}(Y)$, $z_{l_{a}}(X)\neq z_{l_{a}}(X)$,
and $z_{l_{b}}(X)\neq z_{l_{b}}(Y).$ Suppose

\[
L\left[(z_{l_{0}}(X),z_{l_{a}}(X)\right]\neq L\left[(z_{l_{0}}(Y),z_{l_{a}}(Y)\right]\text{ and }
\]
\[
L\left[(z_{l_{0}}(X),z_{l_{a}}(X)\right]<L\left[(z_{l_{0}}(Y),z_{l_{a}}(Y)\right]
\]
and
\[
L\left[(z_{l_{a}}(X),z_{l_{b}}(X)\right]\neq L\left[(z_{l_{a}}(Y),z_{l_{b}}(Y)\right]\text{ and }
\]
\[
L\left[(z_{l_{a}}(X),z_{l_{b}}(X)\right]>L\left[(z_{l_{a}}(Y),z_{l_{b}}(Y)\right]
\]
such that
\begin{align}
L\left[(z_{l_{0}}(Y),z_{l_{a}}(Y)\right]-L\left[(z_{l_{0}}(X),z_{l_{a}}(X)\right] & =\nonumber \\
 & L\left[(z_{l_{a}}(X),z_{l_{b}}(X)\right]-L\left[(z_{l_{a}}(Y),z_{l_{b}}(Y)\right]\label{eq:excessL=00003DlessL}
\end{align}

By (\ref{eq:excessL=00003DlessL}), we conclude that $L\left[(z_{l_{0}}(X),z_{l_{b}}(X)\right]=L\left[(z_{l_{0}}(Y),z_{l_{b}}(Y)\right].$
Since $t_{a}$ and $t_{b}$ are arbitrary, one can extend the result
to other contour distances. 
\end{proof}
\begin{thm}
Two state spaces $A_{X}(\mathbb{C}_{l})$ and $A_{Y}(\mathbb{C}_{l})$
are identical need not imply the corresponding contours are identical. 
\end{thm}

\begin{proof}
Given that the two state spaces $A_{X}(\mathbb{C}_{l})$ and $A_{Y}(\mathbb{C}_{l})$
are identical. This implies the states in $A_{X}(\mathbb{C}_{l})$
and $A_{Y}(\mathbb{C}_{l})$ are the same. Suppose the order of the
states generated by $\left\{ X\left(z_{l}(t),\mathbb{C}_{l}\right)\right\} _{t\geq t_{0}}$
and $\left\{ Y\left(z_{l}(t),\mathbb{C}_{l}\right)\right\} _{t\geq t_{0}}$
are the same, then the two contours $\gamma_{l}(X)$ and $\gamma_{l}(Y)$
are identical. 

Suppose $z_{l_{0}}(X)=z_{l_{0}}(Y)$ but the randomness has resulted
in distinct order of the states in $A_{X}(\mathbb{C}_{l})$ and $A_{Y}(\mathbb{C}_{l})$
such that
\[
z_{l_{a}}(X)=z_{l_{b}}(Y)\text{ and }z_{l_{b}}(X)=z_{l_{a}}(Y)
\]
 for some arbitrary $t_{a}$ and $t_{b}.$ This implies $\gamma_{l}(X)$
be no more identical to $\gamma_{l}(Y).$ 
\end{proof}
The length of a contour is less sensitive for fluctuations due to
randomness than the contour itself. When $z_{l_{0}}(X)\neq z_{l_{0}}(Y),$
then $\gamma_{l}(X)\neq\gamma_{l}(Y)$ but in the long run, say for
some $t$ $(T,\infty)$, the lengths of $\gamma_{l}(X)$ and $\gamma_{l}(Y)$
could be the same. Once $A_{X}(\mathbb{C}_{l})$ and $A_{Y}(\mathbb{C}_{l})$
are formed by the two non-disjoint processes $\left\{ X\left(z_{l}(t),\mathbb{C}_{l}\right)\right\} _{t\geq t_{0}}$
and $\left\{ Y\left(z_{l}(t),\mathbb{C}_{l}\right)\right\} _{t\geq t_{0}}$,
then
\[
p\left(z_{l_{0}}(X),z_{l_{0}}(X),t\right)^{(n)}\left(t\in[t_{0},\infty)\right)=0\text{ for }n=1,2,...
\]
\[
p\left(z_{l_{0}}(Y),z_{l_{0}}(Y),t\right)^{(n)}\left(t\in[t_{0},\infty)\right)=0\text{ for }n=1,2,...
\]
and

\[
p\left(z_{l_{a}}(X),z_{l_{a}}(X),t\right)^{(n)}\left(t\in[t_{0},\infty)\right)=0\text{ for }n=1,2,...
\]
for an arbitrary $z_{l_{a}}(X)$ chosen by $X\left(z_{l}(t),\mathbb{C}_{l}\right)$,
\[
p\left(z_{l_{a}}(Y),z_{l_{a}}(Y),t\right)^{(n)}\left(t\in[t_{0},\infty)\right)=0\text{ for }n=1,2,...
\]
for an arbitrary $z_{l_{a}}(Y)$ chosen by $Y\left(z_{l}(t),\mathbb{C}_{l}\right)$.
These imply,

\[
\sum_{n=1}^{\infty}p\left(z_{l_{a}}(X),z_{l_{a}}(X),t\right)^{(n)}\left(t\in[t_{0},\infty)\right)=0
\]
\[
\sum_{n=1}^{\infty}p\left(z_{l_{a}}(Y),z_{l_{a}}(Y),t\right)^{(n)}\left(t\in[t_{0},\infty)\right)=0.
\]
Suppose the real valued function $v_{l}(X,\tau)$ maps $[a_{0}^{l},\infty)$
onto the interval $[t_{0},\infty)$, then
\begin{align*}
L(\gamma_{l},X) & =\int_{a_{0}^{l}}^{\infty}\left|z\left[v_{l}(X,\tau)\right]\right|v'_{l}(X,\tau)d\tau\\
 & =\int_{a_{0}^{l}}^{\infty}\left|z\left[v_{l}(Y,\tau)\right]\right|v'_{l}(Y,\tau)d\tau\\
 & =L(\gamma_{l},Y)
\end{align*}
where~ the real valued function $v_{l}(Y,\tau)$ maps $[a_{0}^{l},\infty)$
onto the interval $[t_{0},\infty)$. Only looking at the integral
expressions used for $L(\gamma_{l},X)$ or $L(\gamma_{l},Y)$, we
are unable to tell about whether a contour $\gamma_{l}(X)$ has traveled
to any other planes beyond $\mathbb{C}_{l}.$ The symbol $l$ in $\gamma_{l}(X)$
stands for the plane from which this contour has originated and $X$
stands $X\left(z_{l}(t),\mathbb{C}_{l}\right)$ indicates the random
variable responsible for generating data required to form $\gamma_{l}(X).$
Suppose we consider infinitely many random variables of the type $X\left(z_{l}(t),\mathbb{C}_{l}\right)$
to satisfy two conditions, $(i).$ Each of these works non-disjointly
such that they may choose an initial value that was chosen by a different
random variable, and $(ii).$ Each of these random variables chooses
an initial value that is distinct from others such that the number
of initial values is again the number of random variables. Let $\alpha$
be the number of distinct initial values satisfying the condition
$(i)$ such that $\alpha$ is less than the cardinality of $\mathbb{C}_{l}$
and let $X'$ be the index random variable. Then the total lengths
of all the contours originated by all the random variables of condition
$(i)$ is 
\begin{equation}
L(\gamma_{l},X',\alpha)=\sum_{\alpha}\int_{a_{0}^{l}}^{\infty}\left|z\left[v_{l}(X',\alpha,\tau)\right]\right|v'_{l}(X',\alpha,\tau)d\tau d\alpha\label{eq:totalLalpha}
\end{equation}
Let $\beta$ be the distinct initial values in the condition $(ii)$
due to distinct random variables within the condition $(ii)$ with
a then the total length of all the contours generated due to condition
$(ii)$ is
\begin{equation}
L(\gamma_{l},X',\beta)=\sum_{\beta}\int_{a_{0}^{l}}^{\infty}\left|z\left[v_{l}(X',\beta,\tau)\right]\right|v'_{l}(X',\alpha,\tau)d\tau d\beta.\label{eq:totalLbeta}
\end{equation}
There is no comparative measure between (\ref{eq:totalLalpha}) and
(\ref{eq:totalLbeta}), but 

\begin{align*}
L(\gamma_{l},X',\alpha) & \subset B_{\mathbb{R}}(\mathbb{C}),\\
L(\gamma_{l},X',\beta) & \subset B_{\mathbb{R}}(\mathbb{C}).
\end{align*}
and

\begin{align*}
L(\gamma_{l},X',\alpha) & \subset\bigcup_{\alpha}\bigcup_{z_{l_{a}}\in\mathbb{C}_{l}}D(z_{l_{a}},r_{a},X',\alpha)\\
\\
L(\gamma_{l},X',\beta) & \subset\bigcup_{\beta}\bigcup_{z_{l_{a}}\in\mathbb{C}_{l}}D(z_{l_{a}},r_{a},X',\beta)
\end{align*}
where $D(z_{l_{a}},r_{a},X',\alpha)$ and $D(z_{l_{a}},r_{a},X',\alpha)$
represent discs generated through two-step randomness for conditions
$(i)$ and $(ii).$The procedure for generating these discs remains
the same as described previously. Currently, we have not considered
the spaces created due to overlapping contours by these infinitely
many contours generated due to conditions $(i)$ and $(ii).$ However,
\[
\bigcup_{\alpha}\bigcup_{z_{l_{a}}\in\mathbb{C}_{l}}D(z_{l_{a}},r_{a},X',\alpha)\neq\bigcup_{\beta}\bigcup_{z_{l_{a}}\in\mathbb{C}_{l}}D(z_{l_{a}},r_{a},X',\beta)
\]
if the origins of contours created in $(i)$ and $(ii)$ are different. 
\begin{thm}
Suppose infinitely many random variables of the type $X\left(z_{l}(t),\mathbb{C}_{l}\right)$
are available whose cardinality is same as that of $\mathbb{C}_{l}$
and two different conditions $(i)$ and $(ii)$ above are given. The
union of the sets of discs formed under these two conditions could
be different or the same. 
\end{thm}

\begin{proof}
Let us consider infinitely many random variables within the condition
$(i).$ Let the arbitrary variable be $X'\left(z_{l}(t),\mathbb{C}_{l},\alpha\right)\text{}\left(t\in[t_{0},\infty\right)$
for $\alpha$ as in condition $(i)$ described above. Let $X'\left(z_{l}(t),\mathbb{C}_{l},\alpha\right)=z_{l_{0}}(t_{0},X',\alpha).$
The set of discs formed due to each $X'\left(z_{l}(t),\mathbb{C}_{l},\alpha\right)$
are infinite. Each point on the plane could be the origin of a contour
on $\mathbb{C}_{l}.$ This imply, there is a possibility that 
\begin{equation}
\bigcup_{z_{l_{a}}\in\mathbb{C}_{l}}D(z_{l_{0}},r_{0},X',\alpha)=\bigcup_{z_{l_{a}}\in\mathbb{C}_{l}}D(z_{l_{0}},r_{0},X',\beta)\label{eq:Unionof twodiscs}
\end{equation}
Note that $r_{0}$ is associated with two step randomness of each
$X'\left(z_{l}(t),\mathbb{C}_{l},\alpha\right)$ and $r_{0}$values
are generated separately for conditions $(i)$ and $(ii).$ So, $r_{0}$
in the L.H.S. of (\ref{eq:Unionof twodiscs}) need not be equal to
the $r_{0}$ of R.H.S of (\ref{eq:Unionof twodiscs}). When randomly
each $X'\left(z_{l}(t),\mathbb{C}_{l},\alpha\right)$ chooses different
origins and $r_{0}$s of each $z_{l_{0}}(t_{0},X',\alpha)$ corresponding
to each $X'\left(z_{l}(t),\mathbb{C}_{l},\alpha\right)$ are identical
then (\ref{eq:Unionof twodiscs}) holds. If these $r_{0}s$ are different
then (\ref{eq:Unionof twodiscs}) does not hold. Once (\ref{eq:Unionof twodiscs})
holds, then suppose the $r_{a}$s for $z_{l_{a}}(t_{a},X',\alpha)$
corresponding to each $X'\left(z_{l}(t),\mathbb{C}_{l},\alpha\right)$
are identical to the $r_{a}$s for $z_{l_{a}}(t_{a},X',\beta)$ corresponding
to each $X'\left(z_{l}(t),\mathbb{C}_{l},\beta\right)$ for each of
the infinitely many time intervals, then 

\[
\bigcup_{\alpha}\bigcup_{z_{l_{a}}\in\mathbb{C}_{l}}D(z_{l_{a}},r_{a},X',\alpha)=\bigcup_{\beta}\bigcup_{z_{l_{a}}\in\mathbb{C}_{l}}D(z_{l_{a}},r_{a},X',\beta),
\]
else, if at least one such $r_{a}$ that was chosen randomly is different
in conditions $(i)$ and $(ii)$ then 
\[
\bigcup_{\alpha}\bigcup_{z_{l_{a}}\in\mathbb{C}_{l}}D(z_{l_{a}},r_{a},X',\alpha)\neq\bigcup_{\beta}\bigcup_{z_{l_{a}}\in\mathbb{C}_{l}}D(z_{l_{a}},r_{a},X',\beta).
\]
\end{proof}
We have

\[
\bigcup_{\beta}\bigcup_{z_{l_{a}}\in\mathbb{C}_{l}}D(z_{l_{a}},r_{a},X',\beta)=\mathbb{C}_{l}
\]
because $z_{l_{a}}\in\bigcup_{\beta}\bigcup_{z_{l_{a}}\in\mathbb{C}_{l}}D(z_{l_{a}},r_{a},X',\beta)$
implies $z_{l_{a}}\in\mathbb{C}_{l}$ for each $z_{l_{a}}$ and $z_{l_{a}}\in\mathbb{C}_{l}$
implies $z_{l_{a}}\in\bigcup_{\beta}\bigcup_{z_{l_{a}}\in\mathbb{C}_{l}}D(z_{l_{a}},r_{a},X',\beta).$
For condition $(ii)$ within every disc, there are infinitely many
points of other contours. Whereas, such an assertion is not possible
for the discs generated under the condition $(ii).$ There is no chance
to form an isolated disc under the two step randomness procedure and
Markov property derived earlier still holds for the discs formed under
these two conditions. For a general description of continuous-time
Markov property, refer, for example to \cite{good,bRBhat,ChenMa,Gani,GoswamiRao}. 
\begin{rem}
Under random environment the possibility for having identical $r_{a}$
values in each iteration of $X'\left(z_{l}(t),\mathbb{C}_{l},\alpha\right)$
for infinitely many time interval is very small. So the chances for
below equality can be treated as a rare event:

\begin{align*}
D(z_{l_{a}},r_{a},X',\alpha) & =D(z_{l_{a}},r_{a},X',\beta)\text{ for }t\in(t_{a},t_{a'}]\\
D(z_{l_{a'}},r_{a'},X',\alpha) & =D(z_{l_{a'}},r_{a'},X',\beta)\text{ for }t\in(t_{a'},t_{a''}]\\
\vdots & \vdots\\
D(z_{l_{b}},r_{b},X',\alpha) & =D(z_{l_{b}},r_{b},X',\beta)\text{ for }t\in(t_{b},t_{b}]\\
\vdots & \vdots
\end{align*}
\end{rem}

$ $
\begin{rem}
When we relax the assumption in (\ref{eq:zloNEQzl1}) and (\ref{def:Distinct-complex-numbers})
for the possibility to choose the same state by $X\left(z_{l}(t),\mathbb{C}_{l}\right)$
after that state has been chosen earlier by $X\left(z_{l}(t),\mathbb{C}_{l}\right)$,
then each state in $A_{X}(\mathbb{C}_{l})$ becomes recurrent. For
a recurrent state, the probability to return to a state is certain
even if it takes a very large number of time intervals. We can draw
many contours like $\gamma_{l}(X)$, $\gamma_{l}(Y)$, etc., Each
contour will have its starting point or the origin depending upon
the initial value chosen by the random variable responsible to generate
the data required. A thick forest of contours can be formed from infinitely
many random variables. A family of infinitely many random variables
of type $X\left(z_{l}(t),\mathbb{C}_{l}\right)$ could form a forest
of contours. Let this family be $\mathcal{F}_{l}$ 
\[
\mathcal{F}_{l}=\left\{ X\left(z_{l}(t),\mathbb{C}_{l}\right),Y\left(z_{l}(t),\mathbb{C}_{l}\right),...\right\} 
\]
for the set of random variables defined on $\mathbb{C}_{l}$. Let
$\mathcal{F}_{l}$ satisfies (\ref{eq:zloNEQzl1}) and (\ref{def:Distinct-complex-numbers}).
Each element within $\mathcal{F}_{l}$ will have infinitely many points
called the state spaces. Each state space will have infinitely many
states. A family $\mathcal{F}_{l}$ is recurrent if all the elements
of $\mathcal{F}_{l}$ are recurrent and if not $\mathcal{F}_{l}$
is called transient. A transient family $\mathcal{F}_{l}$ would have
a higher possibility to form a relatively quicker dense forest of
contours than a recurrent family. These dense forests of contours
could spread over one or more planes in bundle $B_{\mathbb{R}}(\mathbb{C})$. 
\end{rem}

\subsection{Loss of Spaces in Bundle $B_{\mathbb{R}}(\mathbb{C})$}

Suppose we continue our investigations of the behavior of $X\left(z_{l}(t),\mathbb{C}_{l}\right)$
with the properties of distinct complex numbers as described in Section
3.1. Let $\gamma_{l}(X,t)$ be the contour formed out of the set of
points $z_{l}(t)$ $(t\in[t_{0},\infty))$ sequentially chosen as
per two step randomness of $\left\{ X\left(z_{l}(t),\mathbb{C}_{l}\right)\right\} _{t\geq t_{0}}$
and has origin in the plane $\mathbb{C}_{l}$. Suppose the space created
by $\gamma_{l}(X,t)$ $(t\in[t_{0},\infty)$ is removed from $B_{\mathbb{R}}(\mathbb{C}).$
Let $[\gamma_{l}(X,t)]^{c}$ be the space of all points of the $B_{\mathbb{R}}(\mathbb{C})$
minus the points of the contour. Let us assume that $\gamma_{l}(X,t)$
is formed out of the distinct complex numbers described in the previous
section. That is,

\[
[\gamma_{l}(X,t)]^{c}=B_{\mathbb{R}}(\mathbb{C})\backslash\gamma_{l}(X,t)=\left\{ z_{l}:z_{l}\in B_{\mathbb{R}}(\mathbb{C})\text{ and }z_{l}\notin\gamma_{l}(X,t)\right\} .
\]
When we introduce another random variable $Y\left(z_{l}(t),\mathbb{C}_{l}\right),$
we assume that the space lost due to $\gamma_{l}(X,t)$ is not available
for $Y\left(z_{l}(t),\mathbb{C}_{l}\right).$ That is, $Y\left(z_{l}(t),\mathbb{C}_{l}\right)$
can choose numbers out of $[\gamma_{l}(X,t)]^{c}.$ The two step randomness
is flexible to choose a radius and the next number within a disc until
a number is found. This implies $\gamma_{l}(Y,t)$ can be formed continuously
without any obstructions from the available numbers of the bundle.
The timing between introducing process $X\left(z_{l}(t),\mathbb{C}_{l}\right)$
and removing a set of data created by $\gamma_{l}(X,t)$ up to a certain
time $t$ would form a removal process. The elements or numbers that
the process $X\left(z_{l}(t),\mathbb{C}_{l}\right)$ occupies over
a time interval will be nothing new and they are part of $B_{\mathbb{R}}(\mathbb{C}).$
Due to the removal of the data occupied by the contour $\gamma_{l}(X,t)$
during the time interval, say $[t_{0},t_{d_{1}}]$, the bundle has
lost elements from it. Since $\left\{ X\left(z_{l}(t),\mathbb{C}_{l}\right)\right\} _{t\geq t_{0}}$
is a continuous process, the contour $\gamma_{l}(Y,t)$ still be forming
after removal process has started at $t_{d_{1}}.$ Suppose $\left\{ X\left(z_{l}(t),\mathbb{C}_{l}\right)\right\} _{t\geq t_{0}}$
has generated discs for a long time intervals up to $t_{b}$ by the
time removal process has started. Here $t_{b}>t_{d_{1}}.$ All the
elements of the contour $\gamma_{l}(X,t)$ that was formed during
$t_{0},t_{d_{1}}]$ are nothing but the elements on the c0ontout until
$t_{d_{1}}$ whose length is 
\[
\int_{a_{0}^{l_{0}}}^{a_{0}^{l_{d_{1}}}}\left|z[v_{l_{d_{1}}}(\tau)]\right|v'_{l_{d_{1}}}(\tau)d\tau
\]
for a real-valued function $v_{l_{d_{1}}}$ mapping $[a_{0}^{l_{0}},a_{0}^{l_{d_{1}}}]$
onto the interval $[t_{0},t_{d_{1}}].$ The set of elements on this
contour are the set $z_{l}(t)$ $(t\in[t_{0},t_{d_{1}}])$ and denoted
by $\gamma_{l}(X,t)$ $(t\in[t_{0},t_{d_{1}}])$. The remaining elements
in the bundle are 
\begin{align*}
B_{\mathbb{R}}(\mathbb{C}) & -\phi_{1}(X)B_{\mathbb{R}}(\mathbb{C})\\
\text{where } & \phi_{1}(X)=\frac{\gamma_{l}(X,t)(t\in[t_{0},t_{d_{1}}])}{B_{\mathbb{R}}(\mathbb{C})}
\end{align*}
and differential equation describing the dynamics is 
\[
\frac{dB_{\mathbb{R}}(\mathbb{C})}{dt}=B_{\mathbb{R}}(\mathbb{C})-\phi_{1}(X)B_{\mathbb{R}}(\mathbb{C})\text{ for (t\ensuremath{\in}[\ensuremath{t_{0}},\ensuremath{t_{d_{1}}}])}.
\]
Suppose we remove the elements of the contour $\gamma_{l}(X,t)$ that
was formed during $(t_{d_{1}},t_{d_{2}}]$ from $B_{\mathbb{R}}(\mathbb{C})$
such that $t_{d_{2}}-t_{d_{1}}=t_{d_{1}}-t_{0}.$ The rate of change
in the bundle $B_{\mathbb{R}}(\mathbb{C})\backslash\gamma_{l}(X,t)$
$(t\in[t_{0},t_{d_{1}}]$ during $(t_{d_{1}},t_{d_{2}}]$ is 
\[
\frac{dB_{\mathbb{R}}(\mathbb{C})}{dt}\mid_{t=t_{d_{1}}}=B_{\mathbb{R}}(\mathbb{C})\mid_{t=t_{d_{1}}}-\phi_{2}(X)B_{\mathbb{R}}(\mathbb{C})\mid_{t=t_{d_{1}}}\text{ for (t\ensuremath{\in}[\ensuremath{t_{d_{1}}},\ensuremath{t_{d_{2}}}])}
\]
where $B_{\mathbb{R}}(\mathbb{C})\mid_{t=t_{d_{1}}}$ indicates the
space of $B_{\mathbb{R}}(\mathbb{C})$ that was available at $t=t_{d_{1}}$and
\[
\phi_{2}(X)=\frac{\gamma_{l}(X,t)(t\in[t_{d_{1}},t_{d_{2}}])}{B_{\mathbb{R}}(\mathbb{C})\mid_{t=t_{d_{1}}}}.
\]
Suppose $t_{d_{1}}$ is randomly chosen, and the rest of all the time
intervals are fixed to maintain the interval lengths equal to $t_{d_{1}}-t_{0}.$
The time intervals have constant length but the piecewise contour
lengths in these intervals need not be identical because the contour
formation is dependent on two step randomness and corresponding discs
formations. The process of removal continues after $t_{d_{1}}$and
$\phi$ be the rate of removal of elements from $B_{\mathbb{R}}(\mathbb{C})$,
then this can be expressed with a differential equation
\begin{equation}
\frac{dB_{\mathbb{R}}(\mathbb{C})}{dt}=B_{\mathbb{R}}(\mathbb{C})-\phi(X)B_{\mathbb{R}}(\mathbb{C}).\label{eq:diff-eq-phi}
\end{equation}
A constant rate of removal of elements is difficult to imagine because
within the each time intervals 
\[
\left\{ [t_{0},t_{d_{1}}],(t_{d_{1}},t_{d_{2}}],...\right\} 
\]
the number of elements to be removed depends on the lengths of contour
formed during these intervals. These contour lengths are 
\begin{align}
\int_{a_{0}^{l_{d_{1}}}}^{a_{0}^{l_{d_{2}}}}\left|z[v_{l_{d_{2}}}(\tau)]\right|v'_{l_{d_{2}}}(\tau)d\tau,\nonumber \\
\int_{a_{0}^{l_{d_{2}}}}^{a_{0}^{l_{d_{3}}}}\left|z[v_{l_{d_{3}}}(\tau)]\right|v'_{l_{d_{3}}}(\tau)d\tau,\label{eq:lengthsofcontours removed}\\
\vdots\nonumber 
\end{align}
We know that the two step randomness creates discs at each iteration
and the space occupied by these discs on $B_{\mathbb{R}}(\mathbb{C})$
need not be identical. That means the lengths of contours formed during
$\left\{ [t_{0},t_{d_{1}}],(t_{d_{1}},t_{d_{2}}],...\right\} $ need
not be identical. The quantity $\phi$ can only be retrospectively
estimated from the data on the sets of elements created by the piecewise
contours within the intervals $\left\{ [t_{0},t_{d_{1}}],(t_{d_{1}},t_{d_{2}}],...\right\} .$
So a better way to express the dynamics due to removal of elements
from $B_{\mathbb{R}}(\mathbb{C})$ due to the removal of piecewise
contours is

\begin{equation}
\frac{dB_{\mathbb{R}}(\mathbb{C})}{dt}=B_{\mathbb{R}}(\mathbb{C})-\phi(X,t)B_{\mathbb{R}}(\mathbb{C}),\label{eq:diff-eq(timedependent}
\end{equation}
where $\phi(t)$ can be approximated by 
\[
\phi(X,t)=\frac{\gamma_{l}(X,t)(t\in[t_{d},t_{d''}])}{B_{\mathbb{R}}(\mathbb{C})}.
\]
Over time (\ref{eq:diff-eq(timedependent}) will produce the dynamics
within bundle $B_{\mathbb{R}}(\mathbb{C}).$ The total elements inside
$B_{\mathbb{R}}(\mathbb{C})$ keep on decreasing due to the removal
of piecewise contours (can be treated as a death rate of data on piecewise
contours). The questions that remain to understand here are if the
rate of removal of contours is faster than the formation of the contours
(a possibility exists), then, does the removal rate becomes an instantaneous
rate? What if the contour $\gamma_{l}(X,t)$ is forming continuously
such that it is spreading into infinitely many planes of $B_{\mathbb{R}}(\mathbb{C}),$
and we start removing the space created by $\gamma_{l}(X,t)$ then
how the dynamics of $B_{\mathbb{R}}(\mathbb{C})$ look like? 

The rate of removal of $\gamma_{l}(X,t)$ in an interval will be zero
if no contour data is available for that interval. The removal of
contour data resumes as soon as the contour data becomes available.
This also implies the removal process could be temporarily discontinued.
By the set-up of the time intervals that are used for removing contours,
the removal rate of contours might be higher than the formation rates
or vice versa, or they both might be identical. First, an interval
of time is decided and within this interval whatever the contour lies
that set of points (numbers) will be removed. If within that chosen
time interval no contour data is available then the removal process
halts temporarily. The removal process resumes once data on contours
becomes available. It is difficult to model a form for $\phi(t)$
because it is dependent on the time interval that was used to remove
for and the length of contour that was formed by the process $\left\{ X\left(z_{l}(t),\mathbb{C}_{l}\right)\right\} _{t\geq t_{0}}$
through the two step randomness. The lengths that will be removed
during these intervals are shown in (\ref{eq:lengthsofcontours removed}).
At a given $t_{M}>t_{0}$ the length of $\gamma_{l}(X,t_{M})(t_{M}\in(t_{0},\infty))$
formed until $t_{M}$ could be larger than the sum of these above
intervals or could be equal, that is 
\[
\]

\begin{align}
\gamma_{l}(X,t_{M})(t_{M}\in(t_{0},\infty))\left\{ \begin{array}{c}
=\int_{a_{0}^{l_{0}}}^{a_{0}^{l_{d_{1}}}}\left|z[v_{l_{d_{1}}}(\tau)]\right|v'_{l_{d_{1}}}(\tau)d\tau+\\
\sum_{i=1}^{\infty}\int_{a_{0}^{l_{d_{i}}}}^{a_{0}^{l_{d_{i+2}}}}\left|z[v_{l_{d_{i+1}}}(\tau)]\right|v'_{l_{d_{i+1}}}(\tau)d\tau\\
\text{( whenever }\phi(t_{M})=0)\\
\\
>\int_{a_{0}^{l_{0}}}^{a_{0}^{l_{d_{1}}}}\left|z[v_{l_{d_{1}}}(\tau)]\right|v'_{l_{d_{1}}}(\tau)d\tau+\\
\sum_{i=1}^{\infty}\int_{a_{0}^{l_{d_{i}}}}^{a_{0}^{l_{d_{i+2}}}}\left|z[v_{l_{d_{i+1}}}(\tau)]\right|v'_{l_{d_{i+1}}}(\tau)d\tau\\
\text{(otherwise)}\\
\text{ }
\end{array}\right.\nonumber \\
\label{eq:Gamma=00003Dand>equation}
\end{align}
Here $\phi(t_{M})=0$ indicates there is no contour data available
that is to be removed from $B_{\mathbb{R}}(\mathbb{C}).$ The event
\begin{align*}
\gamma_{l}(X,t_{M})(t_{M} & \in(t_{0},\infty))<\int_{a_{0}^{l_{0}}}^{a_{0}^{l_{d_{1}}}}\left|z[v_{l_{d_{1}}}(\tau)]\right|v'_{l_{d_{1}}}(\tau)d\tau+\\
 & \sum_{i=1}^{\infty}\int_{a_{0}^{l_{d_{i}}}}^{a_{0}^{l_{d_{i+2}}}}\left|z[v_{l_{d_{i+1}}}(\tau)]\right|v'_{l_{d_{i+1}}}(\tau)d\tau
\end{align*}
is impossible. Whenever 
\begin{align}
\gamma_{l}(X,t)(t_{M} & \in(t_{0},\infty))=\int_{a_{0}^{l_{0}}}^{a_{0}^{l_{d_{1}}}}\left|z[v_{l_{d_{1}}}(\tau)]\right|v'_{l_{d_{1}}}(\tau)d\tau+\nonumber \\
 & \sum_{i=1}^{\infty}\int_{a_{0}^{l_{d_{i}}}}^{a_{0}^{l_{d_{i+2}}}}\left|z[v_{l_{d_{i+1}}}(\tau)]\right|v'_{l_{d_{i+1}}}(\tau)d\tau\label{eq:contourt=00003Dtb}
\end{align}
at $t=t_{b}$ (say), then during $(t_{b},t_{b'}]$ for $(t_{b},t_{b'}]=(t_{0},t_{d_{1}}],$
the amount of data removed could be equal to the amount of $\gamma_{l}(X,t)$
that is available during $(t_{d_{b}},t_{d_{b'}}].$ Also when (\ref{eq:contourt=00003Dtb})
is true at $t=t_{b}$, then
\begin{equation}
\frac{dB_{\mathbb{R}}(\mathbb{C})}{dt}=0.\label{eq:dBR=00003D0}
\end{equation}
Satisfying (\ref{eq:dBR=00003D0}) does not indicate a removal process
has attained stationary solution or a steady state solution. As noted
earlier after attaining (\ref{eq:dBR=00003D0}) at some $t>t_{0}$,
the rate of removal continues soon after the formation of a new piece
of contour in $\gamma_{l}(X,t).$ 
\begin{thm}
\label{thm:NoStabilityX}The differential equation describing the
removal process
\[
\frac{dB_{\mathbb{R}}(\mathbb{C})}{dt}=B_{\mathbb{R}}(\mathbb{C})-\phi(X,t)B_{\mathbb{R}}(\mathbb{C})
\]
never attains global stability. 
\end{thm}

\begin{proof}
The removal process of the data generated by $\gamma_{l}(X,t)$ continues
even after $\frac{dB_{\mathbb{R}}(\mathbb{C})}{dt}=0$ at $t$ for
$t\in(t_{0},\infty)$. The amount of $\phi(X,t)$ after $\frac{dB_{\mathbb{R}}(\mathbb{C})}{dt}=0$
depends on the availability of the length of $\gamma_{l}(X,t)$ just
after attaining $\frac{dB_{\mathbb{R}}(\mathbb{C})}{dt}=0$ and it
could be smaller than the set of data points generated by the piece
of contour formed after 
\begin{align*}
\gamma_{l}(X,t)(t_{M} & \in(t_{0},\infty))=\int_{a_{0}^{l_{0}}}^{a_{0}^{l_{d_{1}}}}\left|z[v_{l_{d_{1}}}(\tau)]\right|v'_{l_{d_{1}}}(\tau)d\tau+\\
 & \sum_{i=1}^{\infty}\int_{a_{0}^{l_{d_{i}}}}^{a_{0}^{l_{d_{i+2}}}}\left|z[v_{l_{d_{i+1}}}(\tau)]\right|v'_{l_{d_{i+1}}}(\tau)d\tau
\end{align*}
or equal to the piece of the contour formed. We could never attain
a situation of 
\[
\int_{a_{0}^{l_{db}}}^{a_{0}^{l_{d_{b'}}}}\left|z[v_{l_{d_{b'}}}(\tau)]\right|v'_{l_{d_{b'}}}(\tau)d\tau<\epsilon
\]
for $\epsilon>0$ chosen. Hence, the rate of removal of the space
of data in $B_{\mathbb{R}}(\mathbb{C})$ can never attain stability
as long as the contour formation process $\left\{ X\left(z_{l}(t),\mathbb{C}_{l}\right)\right\} _{t\geq t_{0}}$
continues. 
\end{proof}
\begin{thm}
Suppose $R_{M}$ be an upper bound such that the length of the contour
removed
\begin{equation}
\int_{a_{0}^{l_{db}}}^{a_{0}^{l_{d_{b'}}}}\left|z[v_{l_{d_{b'}}}(\tau)]\right|v'_{l_{d_{b'}}}(\tau)d\tau\leq R_{M}\label{eq:length<RM}
\end{equation}
for an arbitrary interval $(t_{d_{b}},t_{d_{b'}}].$ Then such an
$R_{M}$ does not exist for all the intervals of the type $(t_{d_{b}},t_{d_{b'}}].$
\end{thm}

\begin{proof}
Suppose the quantity $R_{M}$ exists for all the intervals of the
type $t_{d_{b}},t_{d_{b'}}]$ such that (\ref{eq:length<RM}) is true.
This implies for any given arbitrary interval $(t_{d_{c}},t_{d_{c'}}]$
where $(t_{d_{c}},t_{d_{c'}}]$ occurred prior to $(t_{d_{b}},t_{d_{b'}}]$
or has occurred after $(t_{d_{b}},t_{d_{b'}}]$, but the length of
the contour whose data to be removed does not exceed $R_{M}.$ Such
an assertion is true only if $R_{M}\rightarrow\infty,$ and not for
a finite $R_{M}$ because the piece of the contour $\gamma_{l}(X,t)$
whose data to be removed depends on the length of the contour that
is available. This implies, there is no upper limit for the length
of the contour to be formed. This contradicts that $R_{M}$ can be
attained such that (\ref{eq:length<RM}) holds. The set of the data
created by the length 
\[
\int_{a_{0}^{l_{0}}}^{a_{0}^{l_{d_{1}}}}\left|z[v_{l_{d_{1}}}(\tau)]\right|v'_{l_{d_{1}}}(\tau)d\tau+\sum_{i=1}^{\infty}\int_{a_{0}^{l_{d_{i}}}}^{a_{0}^{l_{d_{i+2}}}}\left|z[v_{l_{d_{i+1}}}(\tau)]\right|v'_{l_{d_{i+1}}}(\tau)d\tau
\]
could reach $\gamma_{l}(X,t)$ from left (or from below) once or more
than once. Contour formation and corresponding removal process once
initiated will continue forever. 
\end{proof}
One can also use a different strategy to remove a space of data points
formed by $\gamma_{l}(X,t)$ for $t\in[t_{0},t_{b}].$ Suppose we
assume $\phi(X,t)$ follows a certain parametric form to decide the
number of elements of the set $z_{l}(t)$ on $\gamma_{l}(X,t)$ over
various time intervals. say $\left\{ [t_{0},t_{d_{1}}],(t_{d_{1}},t_{d_{2}}],...\right\} .$
We will know the length of $\gamma_{l}(X,t)$ at $t_{b},$ which is
\begin{equation}
\int_{a_{0}^{l_{0}}}^{a_{0}^{l_{d_{b}}}}\left|z[v_{l_{d_{b}}}(\tau)]\right|v'_{l_{d_{b}}}(\tau)d\tau.\label{eq:lengthofGammaattb}
\end{equation}
So we choose the removal rate of the set of data created on this contour
up to $t_{b}$ such that $\phi(t)$ at each interval $\left\{ [t_{0},t_{d_{1}}],(t_{d_{1}},t_{d_{2}}],...\right\} $
is less than the corresponding pieces of the contours formed during
$\left\{ [t_{b},t_{b'}],(t_{b'},t_{b''}],...\right\} .$ Note that
these two sets of intervals need not have same interval lengths. The
intervals to form $\gamma_{l}(X,t)$ are emerged out of two step randomness.
At $t_{b}$, we first form 
\begin{equation}
D\left(z_{l_{b}},r_{b}\right)\label{eq:DISCatzlb}
\end{equation}
using $r_{b}$ and $\left\{ X\left(z_{l}(t),\mathbb{C}_{l}\right)\right\} _{t\geq t_{0}}$
chooses $z_{l_{b'}}$ from (\ref{eq:DISCatzlb}). If we choose $\phi(X,t)$
such that 
\begin{equation}
\phi(X,t)=\psi(X,t)\int_{a_{0}^{l_{db}}}^{a_{0}^{l_{d_{b'}}}}\left|z[v_{l_{d_{b'}}}(\tau)]\right|v'_{l_{d_{b'}}}(\tau)d\tau\text{ \ensuremath{\left(t\in(t_{b},t_{b'}]\right)}}\label{eq:phiinintegral}
\end{equation}
for $0<\psi(X,t)<1.$ The dynamics in the bundle would be
\begin{align*}
\frac{dB_{\mathbb{R}}(\mathbb{C})}{dt} & =B_{\mathbb{R}}(\mathbb{C})-\left[\psi(X,t)\int_{a_{0}^{l_{d_{b}}}}^{a_{0}^{l_{d_{b'}}}}\left|z[v_{l_{d_{b'}}}(\tau)]\right|v'_{l_{d_{b'}}}(\tau)d\tau\text{ \ensuremath{\left(t\in(t_{b},t_{b'}]\right)}}\right]\times\\
 & B_{\mathbb{R}}(\mathbb{C})\text{ }\left(t\in(t_{b},t_{b'}]\right),
\end{align*}
and for each of the interval, we can choose an $\psi(t)$ or it could
be a constant value in $(0,1).$ We assure through (\ref{eq:phiinintegral})
that the set of numbers on $\gamma_{l}(X,t)$ removed during $(t_{b},t_{b'}]$
are less than the set of numbers formed on contour during $(t_{b},t_{b'}].$
This way the data of the contour $\gamma_{l}(X,t)$ remaining unremoved
are at least the set of data points that required to draw the distance
(\ref{eq:lengthofGammaattb}). Similarly, the dynamics in bundle due
to the removal of the set of numbers (data points) removed during
$(t_{b'},t_{b''}]$ is

\begin{align*}
\frac{dB_{\mathbb{R}}(\mathbb{C})}{dt} & =B_{\mathbb{R}}(\mathbb{C})-\left[\psi(X,t)\int_{a_{0}^{l_{d_{b'}}}}^{a_{0}^{l_{d_{b''}}}}\left|z[v_{l_{d_{b''}}}(\tau)]\right|v'_{l_{d_{b''}}}(\tau)d\tau\text{ \ensuremath{\left(t\in(t_{b'},t_{b''}]\right)}}\right]\times\\
 & B_{\mathbb{R}}(\mathbb{C})\text{ }\left(t\in(t_{b'},t_{b''}]\right).
\end{align*}
Through the strategy explained here, the removal of the space over
the long period of time can be approximated by 
\begin{align}
\frac{dB_{\mathbb{R}}(\mathbb{C})}{dt} & =B_{\mathbb{R}}(\mathbb{C})-\left[\psi(X,t)\int_{a_{0}^{l_{0}}}^{a_{0}^{l_{\infty}}}\left|z[v_{l_{\infty}}(\tau)]\right|v'_{l_{\infty}}(\tau)d\tau\text{ \ensuremath{\left(t\in(t_{b},t_{\infty}]\right)}}\right]\times\nonumber \\
 & B_{\mathbb{R}}(\mathbb{C})\text{ }\left(t\in([t_{0},\infty]\right)\nonumber \\
\label{eq:Db/dtforinfinotelength}
\end{align}
The differential equation (\ref{eq:Db/dtforinfinotelength}) gives
an approximation of overall dynamics generated in various intervals
$\left\{ [t_{b},t_{b'}],(t_{b'},t_{b''}],...\right\} .$ The amount
of space removed would never be able to reach a situation where $\phi(t)=0$
in these differential equations because the data points due to the
length of the contour in (\ref{eq:lengthofGammaattb}) will be still
in excess. Through the differential equation (\ref{eq:Db/dtforinfinotelength})
we made sure that $\gamma_{l}(X,t_{M})(t_{M}\in(t_{0},\infty))$ of
(\ref{eq:Gamma=00003Dand>equation}) satisfies 
\begin{align}
\gamma_{l}(X,t_{M})(t_{M} & \in(t_{0},\infty))>\int_{a_{0}^{l_{0}}}^{a_{0}^{l_{d_{1}}}}\left|z[v_{l_{d_{1}}}(\tau)]\right|v'_{l_{d_{1}}}(\tau)d\tau+\nonumber \\
 & \sum_{i=1}^{\infty}\int_{a_{0}^{l_{d_{i}}}}^{a_{0}^{l_{d_{i+2}}}}\left|z[v_{l_{d_{i+1}}}(\tau)]\right|v'_{l_{d_{i+1}}}(\tau)d\tau\label{eq:gammal>twosums}
\end{align}
if the removal process follows $\phi(X,t)$ of (\ref{eq:phiinintegral}).
There are no specific advantages of if (\ref{eq:contourt=00003Dtb})
holds unless we are having any difficulties with discontinuity of
the removal process. 

Suppose we wanted to introduce infinitely many random processes to
generate contours as in (\ref{eq:totalLalpha}) and (\ref{eq:totalLbeta})
and then initiate corresponding removal processes. The space of the
data lost in $B_{\mathbb{R}}(\mathbb{C})$ over a period of time intervals
and constructions of such sets would involve careful considerations
of contour formation and removal processes. For the sake of understanding
the dynamics in $B_{\mathbb{R}}(\mathbb{C})$ due to these multiple
contour formation and removal processes, let us introduce a second
process $\left\{ Y\left(z_{l}(t),\mathbb{C}_{l}\right)\right\} _{t\geq t_{c}}$
for $t_{c}>t_{b}.$ recollect that when $\left\{ X\left(z_{l}(t),\mathbb{C}_{l}\right)\right\} _{t\geq t_{0}}$
has reached $t=t_{b}$ we have introduced the removal process of $\gamma_{l}(X,t).$
This implies, $\left\{ Y\left(z_{l}(t),\mathbb{C}_{l}\right)\right\} _{t\geq t_{c}}$
is introduced $t_{c}-t_{b}$ time units after a removal process of
$\gamma_{l}(X,t)$ was initiated, and $t_{c}$ time units after the
process $\left\{ X\left(z_{l}(t),\mathbb{C}_{l}\right)\right\} _{t\geq t_{0}}$
was introduced in bundle $B_{\mathbb{R}}(\mathbb{C}).$ At the time
of introduction of $\left\{ Y\left(z_{l}(t),\mathbb{C}_{l}\right)\right\} _{t\geq t_{c}}$
the bundle $B_{\mathbb{R}}(\mathbb{C})$ has lost some set of data
points due to the removal process of $\gamma_{l}(X,t)$ at $t=t_{b}.$
The two step randomness of $\left\{ Y\left(z_{l}(t),\mathbb{C}_{l}\right)\right\} $
will choose a number $z_{l}(Y,t)$ to $t=t_{c}$ on the plane $\mathbb{C}_{l}.$
Let us call this initial value of the new contour $z_{l_{0}}(Y)$
and the contour be $\gamma_{l}(Y,t).$ Since the space of the contour
$\gamma_{l}(Y,t)$ for the period $t_{0}$ to $t_{c}$ was removed,
the number $z_{l_{0}}(Y)$ will be point of $\mathbb{C}_{l}$ such
that it is within the set 
\[
\left\{ \mathbb{C}_{l}\backslash z_{l}(X,t)\text{ \ensuremath{\left(t\in[t_{0},t_{c}]\right)}}\right\} ,
\]
where
\[
\mathbb{C}_{l}\backslash z_{l}(X,t)\text{ \ensuremath{\left(t\in[t_{0},t_{c}]\right)}}=\left\{ z_{l}:z_{l}\in\mathbb{C}_{l}\right\} \text{ and }z_{l}\notin z_{l}(X,t)\text{ for \ensuremath{t\in[t_{0},t_{c}]}}.
\]
\begin{figure}
\includegraphics[scale=0.5]{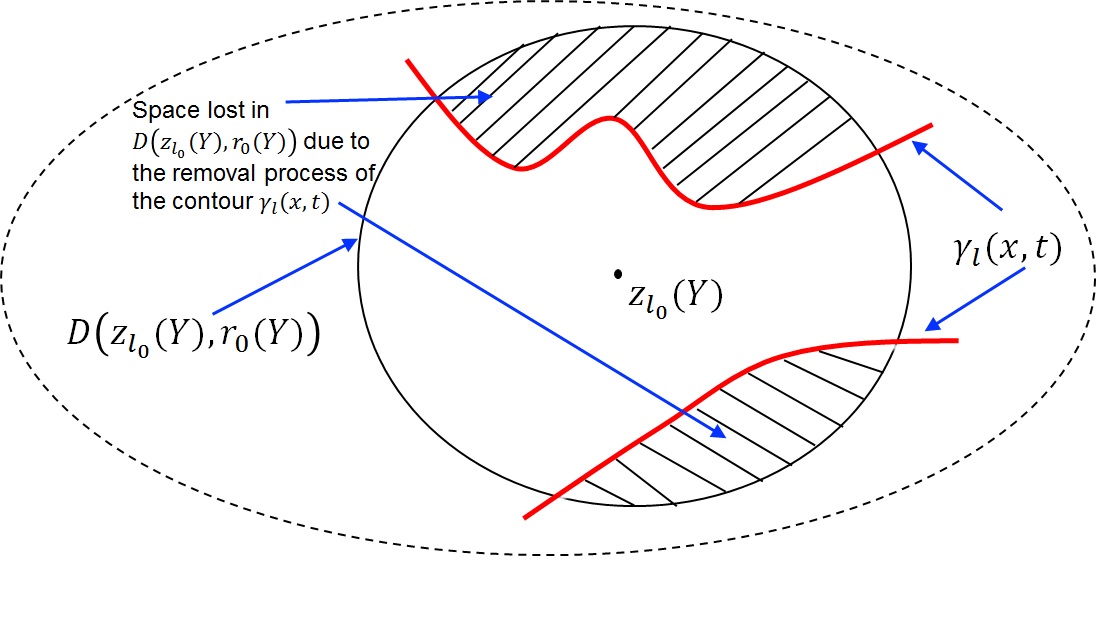}\caption{\label{fig:Non-availability-of-thespace}Non-availability of the space
for $\left\{ Y\left(z_{l}(t),\mathbb{C}_{l}\right)\right\} $ in a
disc after the time $t_{c}$ due to removal of a piece of contour
$\gamma_{l}(x,t)$. The shaded region cannot be reached from $z_{l_{0}}$
within a disc $D\left(z_{l_{0}}(Y),r_{0}(Y)\right)$. }
\end{figure}

Once $z_{l_{0}}(Y)$ is chosen by $Y\left(z_{l}(t),\mathbb{C}_{l}\right)$,
then a disc 
\begin{equation}
D\left(z_{l_{0}}(Y),r_{0}(Y)\right)\subset\mathbb{C}_{l}\label{eq:newinitialdiscofY}
\end{equation}
 will be formed. As the second step $Y\left(z_{l}(t),\mathbb{C}_{l}\right)$
will pick a point within the disc (\ref{eq:newinitialdiscofY}) such
that the points of this disc are all part of the available space in
$\mathbb{C}_{l}$ that was available after points (numbers) due to
the removal process of the contour $\gamma_{l}(X,t)$ for $t=t_{0}$
to $t=t_{c}$ is implemented. Once $z_{l_{0}}(Y)$ is chosen and the
disc (\ref{eq:newinitialdiscofY}) is formed, there may be situation
that the entire space of this disc is not available for choosing $z_{l_{1}}$such
that a contour or a piecewise contour can not be formed by joining
$z_{l_{0}}$to $z_{l_{1}}.$ The process $\left\{ Y\left(z_{l}(t),\mathbb{C}_{l}\right)\right\} $
will have to pick a number randomly only in certain locations of disc.
See Figure \ref{fig:Non-availability-of-thespace}. The space within
$D\left(z_{l_{0}}(Y),r_{0}(Y)\right)$ is divided into three components,
namely, the set of points due to the removal of the space of the contour
$\gamma_{l}(X,t)$, say, $S_{1}\left[D\left(z_{l_{0}}(Y),r_{0}(Y)\right)\right]$,
the set of points available in the disc to which a contour or piecewise
arcs can be drawn from $z_{l_{0}}(Y)$ to $z_{l_{1}}(Y),$ say, $S_{2}\left[D\left(z_{l_{0}}(Y),r_{0}(Y)\right)\right]$,
and the set of points available, say, $S_{3}\left[D\left(z_{l_{0}}(Y),r_{0}(Y)\right)\right]$
in the disc to choose $z_{l_{1}}(Y)$ but a contour passing from $z_{l_{0}}(Y)$
to $z_{l_{1}}(Y)$ for $z_{l_{1}}(Y)\in D\left(z_{l_{0}}(Y),r_{0}(Y)\right)$
can not be drawn. The removal process has caused disc $D\left(z_{l_{0}}(Y),r_{0}(Y)\right)$
to write below as a union of these three sets
\begin{align}
D\left(z_{l_{0}}(Y),r_{0}(Y)\right)=S_{1}\left[D\left(z_{l_{0}}(Y),r_{0}(Y)\right)\right]\cup S_{2}\left[D\left(z_{l_{0}}(Y),r_{0}(Y)\right)\right]\cup\nonumber \\
S_{3}\left[D\left(z_{l_{0}}(Y),r_{0}(Y)\right)\right]\label{eq:DISjointuniounofD0Y}
\end{align}

Here $S_{1},$$S_{2},$ and $S_{3}$ are disjoint. A point (number)
to which a contour can be drawn from $z_{l_{0}}(Y)$ is located within
the set $S_{2}\left[D\left(z_{l_{0}}(Y),r_{0}(Y)\right)\right].$
Similarly Let be $z_{l_{a}}(Y)$ an arbitrary point available to draw
a contour from a previous iteration and was chosen by $\left\{ Y\left(z_{l}(t),\mathbb{C}_{l}\right)\right\} $
at some $t$. Using $z_{l_{a}}(Y)$ we can draw a disc $D\left(z_{l_{a}}(Y),r_{a}(Y)\right)$
such that

\begin{align*}
D\left(z_{l_{a}}(Y),r_{a}(Y)\right) & =S_{1}\left[D\left(z_{l_{a}}(Y),r_{a}(Y)\right)\right]\cup S_{2}\left[D\left(z_{l_{a}}(Y),r_{a}(Y)\right)\right]\cup\\
 & S_{1}\left[D\left(z_{l_{a}}(Y),r_{a}(Y)\right)\right].
\end{align*}
If $r_{l}(X,t)\cap D\left(z_{l_{a}}(Y),r_{a}(Y)\right)=\phi\text{ (empty set)},$then
\[
S_{1}\left[D\left(z_{l_{a}}(Y),r_{a}(Y)\right)\right]=\phi\text{ and }S_{3}\left[D\left(z_{l_{a}}(Y),r_{a}(Y)\right)\right]=\phi
\]
and
\[
D\left(z_{l_{a}}(Y),r_{a}(Y)\right)=S_{2}\left[D\left(z_{l_{a}}(Y),r_{a}(Y)\right)\right].
\]
Any point in $D\left(z_{l_{a}}(Y),r_{a}(Y)\right)$ can be randomly
chosen by $\left\{ Y\left(z_{l}(t),\mathbb{C}_{l}\right)\right\} $
such that a contour can be drawn within $D\left(z_{l_{a}}(Y),r_{a}(Y)\right).$
Suppose in a given $D\left(z_{l_{a}}(Y),r_{a}(Y)\right)$, the next
iteration point, say, $z_{l_{a'}}(Y)$ for $z_{l_{a}}(Y)$ for $z_{l_{a'}}(Y)\in D\left(z_{l_{a}}(Y),r_{a}(Y)\right)$
lies such that a direct contour from $z_{l_{a}}(Y)$ to $z_{l_{a'}}(Y)$
cannot be drawn due to a deleted space of the contour $\gamma_{l}(X,t).$
Suppose there exists some space outside the deleted space of $\gamma_{l}(X,t)$
within $D\left(z_{l_{a}}(Y),r_{a}(Y)\right)$ so that a piecewise
arcs can be drawn from $z_{l_{a}}(Y)$ to $z_{l_{a'}}(Y)$. In such
situations $Y\left(z_{l}(t),\mathbb{C}_{l}\right)$ will generate
a set of points around the deleted space to draw piecewise arcs to
join $z_{l_{a}}(Y)$ to $z_{l_{a'}}(Y)$ whose distance is 
\begin{equation}
\int_{a_{0}^{l_{d_{a}}}}^{a_{0}^{l_{d_{a'}}}}\left|z[v_{l_{d_{a}}}(Y,\tau)]\right|v'_{l_{d_{a'}}}(\tau)d\tau.\label{eq:distancearounddeletedSpace}
\end{equation}
Here (\ref{eq:distancearounddeletedSpace}) will be the sum of piecewise
arcs such that
\begin{align}
\int_{a_{0}^{l_{d_{a}}}}^{a_{0}^{l_{d_{a'}}}}\left|z[v_{l_{d_{a}}}(Y,\tau)]\right|v'_{l_{d_{a'}}}(\tau)d\tau & =\int_{a_{0}^{l_{d_{a}}}}^{a_{0}^{l_{a(i)}}}\left|z[v_{l_{a(i)}}(Y,\tau)]\right|v'_{l_{a(i)}}(\tau)d\tau+\nonumber \\
 & \sum_{i=1}^{k}\int_{a_{0}^{l_{d_{a(i)}}}}^{a_{0}^{l_{a(i+1)}}}\left|z[v_{l_{a(i+1)}}(Y,\tau)]\right|v'_{l_{a(i+1)}}(\tau)d\tau+\nonumber \\
 & \int_{a_{0}^{l_{a(k+1)}}}^{a_{0}^{l_{d_{a'}}}}\left|z[v_{l_{d_{a'}}}(Y,\tau)]\right|v'_{l_{d_{a'}}}(\tau)d\tau\label{eq:distancepiecearcsD}
\end{align}
The process of creation of $\gamma_{l}(Y,t)$ for $t>t_{c}$ continues
as described above. Contour $\gamma_{l}(Y,t)$ can have a space of
$\gamma_{l}(X,t)$ that does not gets deleted due to the removal process
of $X\left(z_{l}(t),\mathbb{C}_{l}\right).$ All the points of the
set $\gamma_{l}(X,t)\cap\gamma_{l}(Y,t)$ such that 
\[
\gamma_{l}(X,t)\cap\gamma_{l}(Y,t)\neq\phi\text{ (empty set)}
\]
can get deleted due to the removal process of $X\left(z_{l}(t),\mathbb{C}_{l}\right).$
The process $Y\left(z_{l}(t),\mathbb{C}_{l}\right)$ until a removal
process for $Y\left(z_{l}(t),\mathbb{C}_{l}\right)$ is introduced,
will not influence the differential equation describing the loss of
space in bundle $B_{\mathbb{R}}(\mathbb{C})$ in (\ref{eq:Db/dtforinfinotelength}).
Let us introduce removal process for $\gamma_{l}(Y,t)$ at $t=t_{g}$
i.e. starting at the interval $(t_{g-1},t_{g}].$ Suppose a length
of this contour equivalent to 
\[
\int_{a_{0}^{l_{c}}}^{a_{0}^{l_{g}}}\left|z[v_{l_{g}}(Y,\tau)]\right|v'_{l_{g}}(\tau)d\tau
\]
 is always maintained between the new contour formation and removal
location of this contour such that removal rate never becomes zero.
Suppose $z_{lg(1)}$ be the point chosen in $D\left(z_{l_{g}}(Y),r_{g}(Y)\right)\subset\mathbb{C}_{l}$
where 
\[
z_{lg(1)}(Y)\in S_{2}\left[D\left(z_{l_{g}}(Y),r_{g}(Y)\right)\right]
\]
\begin{align}
D\left(z_{l_{g}}(Y),r_{g}(Y)\right) & =S_{1}\left[D\left(z_{l_{g}}(Y),r_{g}(Y)\right)\right]\cup S_{2}\left[D\left(z_{l_{g}}(Y),r_{g}(Y)\right)\right]\cup\nonumber \\
 & S_{3}\left[D\left(z_{l_{g}}(Y),r_{g}(Y)\right)\right].\label{eq:disjointdiscsinY}
\end{align}
and

The length from $z_{lg}(Y)$ to $z_{lg(1)}(Y)$ is 
\begin{equation}
\int_{a_{0}^{l_{g}}}^{a_{0}^{l_{g(1)}}}\left|z[v_{l_{g(1)}}(Y,\tau)]\right|v'_{l_{g(1)}}(\tau)d\tau.\label{eq:lengthneededforpsi}
\end{equation}
The removal rate of $Y\left(z_{l}(t),\mathbb{C}_{l}\right)$ we denote
here by $\phi(Y,t)$. The value of $\phi(Y,t)$ during $(t_{g},t_{g(1)}]$
is expressed using (\ref{eq:lengthneededforpsi}) as
\[
\phi(Y,t)=\psi(Y,t)\int_{a_{0}^{l_{g}}}^{a_{0}^{l_{g(1)}}}\left|z[v_{l_{g(1)}}(Y,\tau)]\right|v'_{l_{g(1)}}(\tau)d\tau,
\]
 for $0<\psi(Y,t)<1,$ and the value of $\phi(Y,t)$ during $(t_{g},t_{\infty}]$
is expressed using (\ref{eq:lengthneededforpsi}) as 
\[
\phi(Y,t)=\psi(Y,t)\int_{a_{0}^{l_{g}}}^{a_{0}^{l_{\infty}}}\left|z[v_{l_{\infty}}(Y,\tau)]\right|v'_{l_{\infty}}(\tau)d\tau
\]
The dynamics in bundle $B_{\mathbb{R}}(\mathbb{C})$ due to removal
of set of points in $\gamma_{l}(X,t)$ $\left(t\in[t_{b},\infty\right)$
and in $\gamma_{l}(Y,t)$ $\left(t\in[t_{g},\infty\right)$described
above can be divided into below four parts:

$(i)$ Removal of data points due to the removal process introduced
on contour $\gamma_{l}(X,t),$

$(ii)$ Removal of data points due to the removal process introduced
on contour $\gamma_{l}(Y,t),$

$(iii)$ Removal of data points in the set $\gamma_{l}(X,t)\cap\gamma_{l}(Y,t)$
for $\gamma_{l}(X,t)\cap\gamma_{l}(Y,t)\neq\phi\text{ (empty set),}$due
to the removal process introduced on contour $\gamma_{l}(X,t),$

$(iv)$ Removal of data points in the set $\gamma_{l}(X,t)\cap\gamma_{l}(Y,t)$
for $\gamma_{l}(X,t)\cap\gamma_{l}(Y,t)\neq\phi\text{ (empty set),}$due
to the removal process introduced on contour $\gamma_{l}(Y,t).$

Let $\phi(X,t)$ and $\phi(Y,t)$ represent removal rates for the
points purely on $\gamma_{l}(x,t)$ and $\gamma_{l}(Y,t)$ and not
on $\gamma_{l}(X,t)\cap\gamma_{l}(Y,t)\neq\phi\text{ (empty set).}$Let
$\phi_{3}(X,t)$ represent removal rates for the points purely on
$\gamma_{l}(X,t)\cap\gamma_{l}(Y,t)\neq\phi\text{ (empty set)}$for
the contour initiated by $X\left(z_{l}(t),\mathbb{C}_{l}\right)$
and $\phi_{4}(Y,t)$ represent removal rates for the points purely
on $\gamma_{l}(X,t)\cap\gamma_{l}(Y,t)\neq\phi\text{ (empty set)}$
for the contour initiated by $Y\left(z_{l}(t),\mathbb{C}_{l}\right).$
The dynamics in bundle $B_{\mathbb{R}}(\mathbb{C})$ due to four parts
above is expressed through the differential equation:
\begin{equation}
\begin{array}{cc}
\frac{dB_{\mathbb{R}}(\mathbb{C})}{dt} & =B_{\mathbb{R}}(\mathbb{C})-\phi(X,t)B_{\mathbb{R}}(\mathbb{C})-\phi(Y,t)B_{\mathbb{R}}(\mathbb{C})-\phi_{3}(X,t)B_{\mathbb{R}}(\mathbb{C})-\\
 & \phi_{4}(Y,t)B_{\mathbb{R}}(\mathbb{C})\\
 & =B_{\mathbb{R}}(\mathbb{C})-\left[\psi(X,t)\int_{a_{0}^{l_{0}}}^{a_{0}^{l_{\infty}}}\left|z[v_{l_{\infty}}(X,\tau)]\right|v'_{l_{\infty}}(X,\tau)d\tau\right]B_{\mathbb{R}}(\mathbb{C})\\
 & \nonumber-\left[\psi(Y,t)\int_{a_{0}^{l_{g}}}^{a_{0}^{l_{\infty}}}\left|z[v_{l_{\infty}}(Y,\tau)]\right|v'_{l_{\infty}}(Y,\tau)d\tau\right]B_{\mathbb{R}}(\mathbb{C})-\psi_{3}(X,t)B_{\mathbb{R}}(\mathbb{C})-\\
 & \psi_{4}(Y,t)B_{\mathbb{R}}(\mathbb{C})\\
\\
\end{array}\label{eq:dynamicsduetoXandY}
\end{equation}

for $0<\psi_{3}(X,t)<1$ and $0<\psi_{4}(Y,t)<1$ in the same time
interval in which $\psi(X,t)$ and $\psi(Y,t)$ are implemented. 
\begin{thm}
\label{thm:NoStabilityXY}The differential equation 
\begin{align}
\frac{dB_{\mathbb{R}}(\mathbb{C})}{dt}=B_{\mathbb{R}}(\mathbb{C})-\left[\phi(X,t)+\phi(Y,t)+\phi_{3}(X,t)+\phi_{4}(Y,t)\right]B_{\mathbb{R}}(\mathbb{C})\nonumber \\
\label{eq:dB/dtforXY}
\end{align}
where 
\begin{align*}
\phi(X,t) & =\psi(X,t)\int_{a_{0}^{l_{0}}}^{a_{0}^{l_{\infty}}}\left|z[v_{l_{\infty}}(X,\tau)]\right|v'_{l_{\infty}}(X,\tau)d\tau\text{ \ensuremath{\left(0<\psi(X,t)<1\right)},}\\
\phi(Y,t) & =\psi(Y,t)\int_{a_{0}^{l_{g}}}^{a_{0}^{l_{\infty}}}\left|z[v_{l_{\infty}}(Y,\tau)]\right|v'_{l_{\infty}}(Y,\tau)d\tau\text{ }\left(0<\psi(Y,t)<1\right),\\
\phi_{3}(X,t) & =0<\psi_{3}(X,t)<1,\\
\phi_{4}(Y,t) & =0<\psi_{4}(Y,t)<1,
\end{align*}
will never attain global stability. 
\end{thm}

\begin{proof}
The removal process continuously removes sets of points of contours
$\gamma_{l}(x,t)$ and $\gamma_{l}(y,t)$ described in (\ref{eq:dB/dtforXY}).
As in the proof of the Theorem \ref{thm:NoStabilityX}, the contour
formations happens continuously, and 
\[
\left|B_{\mathbb{R}}(\mathbb{C})\right|<\epsilon
\]
would never arise for every $\epsilon>0$, because of the construction
of $\phi(X,t)$, $\phi(Y,t),$ $\phi_{3}(X,t)$, $\phi_{4}(Y,t)$
given in the statement there will be always space of points in $\phi_{4}(Y,t).$ 
\end{proof}
\begin{rem}
One can also develop an argument similar to to the argument in provided
in the proof of the Theorem \ref{thm:NoStabilityX} to prove Theorem
\ref{thm:NoStabilityXY}. 
\end{rem}

Suppose $X^{a}\left(z_{l}(t),\mathbb{C}_{l}\right)$ represents an
arbitrary random variable out of infinitely many random variables
introduced to form contours with origin in $\mathbb{C}_{l}.$ Let
the removal rate of $\gamma_{l}(X^{a},t)$ be $\phi(X^{a},t)$ and
$\phi_{3}(X^{a},t)$ be the removal rate of the data points at the
intersections of one or more contours. A general differential equation
describing the dynamics due to the removal process of points in the
bundle is
\begin{align}
\frac{dB_{\mathbb{R}}(\mathbb{C})}{dt} & =B_{\mathbb{R}}(\mathbb{C})-\int_{X^{a}}\phi(X^{a},t)dt[B_{\mathbb{R}}(\mathbb{C})]-\nonumber \\
 & \int_{X^{a}}\phi_{3}(X^{a},t)dt[B_{\mathbb{R}}(\mathbb{C})].\label{eq:dynamicsOveralli9nfinity}
\end{align}
In (\ref{eq:dynamicsOveralli9nfinity}), the term $\int_{X^{a}}\phi(X^{a},t)dt$
represent the overall removal rates for the infinitely many contours.
At the each iteration of the equation (\ref{eq:dynamicsOveralli9nfinity}),
the quantity $\int_{X^{a}}\phi(X^{a},t)dt$ will be updated based
on new removal of a certain contour. Similarly, $\int_{X^{a}}\phi_{3}(X^{a},t)dt$
represent removal rates of the sets of points available at the intersections
of the contours. 

\section{\textbf{Islands and Holes in $B_{\mathbb{R}}(\mathbb{C})$}}

\begin{figure}
\includegraphics[scale=0.5]{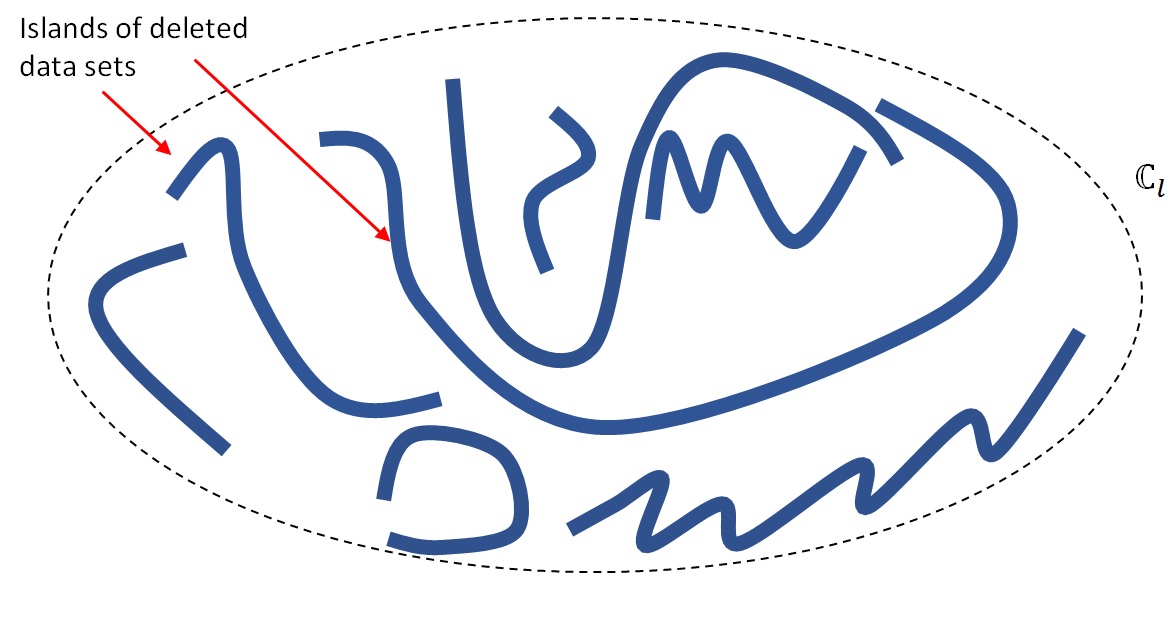}

\caption{\label{fig:Islands-of-deleted}Islands of deleted data due to removal
process on infinitely many contours.}

\end{figure}

The removal process of bundle will create islands of holes due to
overlapping (intersections) of infinitely many contours within $B_{\mathbb{R}}(\mathbb{C}).$
These islands will be never be able to reach again by a newly introduced
two step randomness. See Figure \ref{fig:Islands-of-deleted}.

Suppose we introduce infinitely many random variables of type $X^{a}\left(z_{l}(t),\mathbb{C}_{l}\right)$
that we saw in section 3, but all were introduced at the same time
in $\mathbb{C}_{l}.$ Let the number of these random variables be
such that they are one to one and onto with each member of the complex
plane $\mathbb{C}_{l}.$ Suppose these start creating data for the
formation of infinitely many contours. One contour is assumed not
to block another contour to use its data points. This is explained
further in the following sentences. Suppose $X^{a}\left(z_{l}(t),\mathbb{C}_{l}\right)$
and $X^{b}\left(z_{l}(t),\mathbb{C}_{l}\right)$ be two contours chosen
arbitrarily out of these infinitely many contours that were initiated
at the same time $t_{0}.$ By construction, they have different origins
in $\mathbb{C}_{l}.$ Then the set of complex numbers $\{z_{l}(X^{a},t)\}$
$\left(t\in[t_{0},\infty\right)$ used in the formation of $\gamma_{l}(X^{a},t)$
and the set of numbers used in the formation of $\gamma_{l}(X^{b},t)$
could have a non-empty intersection. That is 
\begin{align}
\left\{ z_{l}:z_{l}\in B_{\mathbb{R}}(\mathbb{C}),\text{\{\ensuremath{z_{l}}(\ensuremath{X^{a}},t)\} \ensuremath{\left(t\in(t_{0},\infty\right)} and }z_{l_{0}}(X^{a})\in\mathbb{C}_{l}\right\} \bigcap\nonumber \\
\left\{ z_{l}:z_{l}\in B_{\mathbb{R}}(\mathbb{C}),\text{\{\ensuremath{z_{l}}(\ensuremath{X^{b}},t)\} \ensuremath{\left(t\in(t_{0},\infty\right)} and }z_{l_{0}}(X^{b})\in\mathbb{C}_{l}\right\} \neq\phi\text{ (empty)}\nonumber \\
\label{eq:4-1}
\end{align}
or 
\begin{align}
\left\{ z_{l}:z_{l}\in B_{\mathbb{R}}(\mathbb{C}),\text{\{\ensuremath{z_{l}}(\ensuremath{X^{a}},t)\} \ensuremath{\left(t\in(t_{0},\infty)\right)} and }z_{l_{0}}(X^{a})\in\mathbb{C}_{l}\right\} \bigcap\nonumber \\
\left\{ z_{l}:z_{l}\in B_{\mathbb{R}}(\mathbb{C}),\text{\{\ensuremath{z_{l}}(\ensuremath{X^{b}},t)\} \ensuremath{\left(t\in(t_{0},\infty)\right)} and }z_{l_{0}}(X^{b})\in\mathbb{C}_{l}\right\} \neq\phi\text{ (empty)}\nonumber \\
\label{eq:4-2}
\end{align}
Any two contours that have different initial values need not be disjoint.
If every $z_{l}(X^{a},t)$ for every $t\in(t_{0},\infty)$ has no
overlap with any element of $\{z_{l}(X^{b},t)\}\left(t\in(t_{0},\infty)\right)$
then that could be purely due to the random environment created in
section 3. Either of these contours or both could be multilevel contours
and have origins in $\mathbb{C}_{l}.$ The formation of multilevel
contours and randomness at $\mathbb{C}_{0}\cap\mathbb{C}_{l}$ described
earlier remains the same. Two contours might have points of intersection
within $B_{\mathbb{R}}(\mathbb{C})$ but such points of intersection
need not behave like common points of $\mathbb{C}_{0}\cap\mathbb{C}_{l}\subset B_{\mathbb{R}}(\mathbb{C}).$
This means the set of points on 
\[
\gamma_{l}(X^{a},t)\cap\gamma_{l}(X^{b},t)
\]
 for which 
\[
\gamma_{l}(X^{a},t)\cap\gamma_{l}(X^{b},t)\neq\phi\text{ (empty)}
\]
 can not be used for changing the plane of the contours. However,
the set of points on 
\begin{equation}
\gamma_{l}(X^{a},t)\cap\gamma_{l}(X^{b},t)\subset\mathbb{C}_{0}\label{eq:4-3}
\end{equation}
 for any two arbitrary random variables $X^{a}\left(z_{l}(t),\mathbb{C}_{l}\right)$
and $X^{b}\left(z_{l}(t),\mathbb{C}_{l}\right)$ could behave similarly
to the points on $\mathbb{C}_{0}\cap\mathbb{C}_{l}.$ The points on
$\mathbb{C}_{0}\cap\mathbb{C}_{p}$ for any arbitrary $\mathbb{C}_{p}\subset B_{\mathbb{R}}(\mathbb{C})$
will have similar properties of forming a multilevel contour as described
in section 3. 

Suppose
\begin{equation}
\gamma_{l}(X^{a},t)\cap\gamma_{l}(X^{b},t)\subset\mathbb{C}_{p}\label{eq:4-4}
\end{equation}
 for some arbitrary plane $\mathbb{C}_{p}\subset B_{\mathbb{R}}(\mathbb{C})$
and $\gamma_{l}(X^{a},t)\cap\gamma_{l}(X^{b},t)$ have origins in
$\mathbb{C}_{l}.$ 

Suppose (\ref{eq:4-4}) satisfied at $t=t_{d},$ then a disc $D\left(z_{l_{d}}(X^{a}),r_{d}(X^{a})\right)$
with center $z_{l_{d}}(X^{a})$ and radius $r_{d}(X^{a})$ is formed
such that 
\begin{equation}
D\left(z_{l_{d}}(X^{a}),r_{d}(X^{a})\right)\subset\mathbb{C}_{p}\label{eq:4-5}
\end{equation}
and next iteration point of $\gamma_{l}(X^{a},t)$ after $z_{l_{d}}(X^{a})$
lies in $\mathbb{C}_{p}$ and not in $\mathbb{C}_{0}\cap\mathbb{C}_{p}$. 

Suppose $z_{l_{d'}}(X^{a})$ be the point generated after $z_{l_{d'}}(X^{a})$
for $z_{l_{d'}}(X^{a})\subset D\left(z_{l_{d}}(X^{a}),r_{d}(X^{a})\right)$,
then 
\[
z_{l_{d'}}(X^{a})\in\mathbb{C}_{p}\text{ and }z_{l_{d'}}(X^{a})\notin\mathbb{C}_{0}\cap\mathbb{C}_{p}.
\]
A contour drawn during $[t_{d},t_{d'}]$ to reach $z_{l_{d'}}(X^{a})$
from $z_{l_{d}}(X^{a})$ with the distance 
\begin{equation}
\int_{a_{0}^{l_{d}}}^{a_{0}^{l_{d'}}}\left|z[u_{l_{d'}}(X^{a},\tau)]\right|u'_{l_{\infty}}(X^{a},\tau)d\tau\text{ ,}\label{eq:4-6}
\end{equation}
lies on $\mathbb{C}_{p}.$ Here the real valued function $u_{l_{d'}}(X^{a},\tau)$
maps $[a_{0}^{l_{d}},a_{0}^{l_{d'}})$ onto the interval $[t_{d},t_{d'}]$.
Because $z_{l_{d}}(X^{a})$ lies on $\gamma_{l}(X^{a},t)$ satisfying
(\ref{eq:4-4}) it could contribute in the next step to form $\gamma_{l}(X^{a},t)$
or $\gamma_{l}(X^{b},t).$ In either situation, the distance in (\ref{eq:4-5})
lies in $\mathbb{C}_{p}.$ Hence a point in $B_{\mathbb{R}}(\mathbb{C})$
if it is in $\mathbb{C}_{0}\cap\mathbb{C}_{p}$ for some arbitrary
$\mathbb{C}_{p}$ have two options to produce a new point on the contour
to continue contour formation. The description of the formation of
contours at the intersection of $\gamma_{l}(X^{a},t)$ and $\gamma_{l}(X^{b},t)$
is also true if there are more than two intersecting contours. we
also note that the lengths of the infinitely many contours up to time
$t_{d}$ which were all introduced at $t_{0}$ could have different
lengths based on the area of the discs formed, and the point chosen
by the corresponding random variable. So the set of lengths
\[
\]
\begin{align}
\left\{ \int_{a_{0}^{l_{d}}}^{a_{0}^{l_{d'}}}\left|z[u_{l_{d'}}(X^{a},\tau)]\right|u'_{l_{d'}}(X^{a},\tau)d\tau,\int_{a_{0}^{l_{d}}}^{a_{0}^{l_{d'}}}\left|z[u_{l_{d'}}(X^{b},\tau)]\right|u'_{l_{d'}}(X^{b},\tau)d\tau,\cdots\right\} \nonumber \\
\label{eq:4-7}
\end{align}
could have different spaces occupied in $B_{\mathbb{R}}(\mathbb{C})$.
The location of each contour after some long time $t_{\infty}$ for
$t_{\infty}>>t_{d'}$ could be anywhere in the bundle and they could
be situated in any plane. The set of lengths 
\[
\]
\begin{align}
\left\{ \int_{a_{0}^{l_{d}}}^{a_{0}^{l_{\infty}}}\left|z[u_{l_{\infty}}(X^{a},\tau)]\right|u'_{l_{\infty}}(X^{a},\tau)d\tau,\int_{a_{0}^{l_{d}}}^{a_{0}^{l_{\infty}}}\left|z[u_{l_{\infty}'}(X^{b},\tau)]\right|u'_{l_{\infty}}(X^{b},\tau)d\tau,\cdots\right\} \nonumber \\
\label{eq:4-8}
\end{align}
and the spaces occupied by 
\[
\gamma_{l}(X^{a},t),\gamma_{l}(X^{b},t),\cdots,
\]
are ever evolving within $B_{\mathbb{R}}(\mathbb{C})$. For a point
$z_{l_{d(2)}}(X^{a})$ $\in\mathbb{C}_{p}$ and $z_{l_{d(2)}}(X^{a})$
$\notin\mathbb{C}_{l}\cap\mathbb{C}_{p}$, the equality
\begin{align}
\int_{a_{0}^{l_{d}}}^{a_{0}^{l_{d(1)}}}\left|z[u_{l_{d(1)}}(X^{a},\tau)]\right|u'_{l_{d(1)}}(X^{a},\tau)d\tau & +\int_{a_{0}^{l_{d(1)}}}^{a_{0}^{l_{d(2)}}}\left|z[u_{l_{d(2)}}(X^{a},\tau)]\right|u'_{l_{d(2)}}(X^{a},\tau)d\tau\nonumber \\
 & =\int_{a_{0}^{l_{d}}}^{a_{0}^{l_{d(2)}}}\left|z[u_{l_{d*(2)}}(X^{a},\tau)]\right|u'_{l_{d*(2)}}(X^{a},\tau)d\tau\label{eq:4-9}
\end{align}
holds for $z_{l_{d}}(X^{a})$ $\in\mathbb{C}_{0}\cap\mathbb{C}_{l}$
and $z_{l_{d(1)}}(X^{a})$ $\in\mathbb{C}_{0}\cap\mathbb{C}_{p}$.
For all such $z_{l_{d}}(X^{b})$ $\in\mathbb{C}_{p}$ and $z_{l_{d(2)}}(X^{b})$
$\notin\mathbb{C}_{l}\cap\mathbb{C}_{p}$, and $z_{l_{d}}(X^{b})\neq z_{l_{d}}(X^{a})$,
$z_{l_{d(1)}}(X^{b})\neq$ $z_{l_{d(1)}}(X^{a})$, and $z_{l_{d(2)}}(X^{b})$
$\neq z_{l_{d(2)}}(X^{a})$ the equality
\begin{align}
\int_{a_{0}^{l_{d}}}^{a_{0}^{l_{d(1)}}}\left|z[u_{l_{d(1)}}(X^{b},\tau)]\right|u'_{l_{d(1)}}(X^{b},\tau)d\tau & +\int_{a_{0}^{l_{d(1)}}}^{a_{0}^{l_{d(2)}}}\left|z[u_{l_{d(2)}}(X^{b},\tau)]\right|u'_{l_{d(2)}}(X^{b},\tau)d\tau\nonumber \\
 & =\int_{a_{0}^{l_{d}}}^{a_{0}^{l_{d(2)}}}\left|z[u_{l_{d*(2)}}(X^{b},\tau)]\right|u'_{l_{d*(2)}}(X^{b},\tau)d\tau\label{eq:4-10}
\end{align}
holds for $z_{l_{d}}(X^{b})$ $\in\mathbb{C}_{0}\cap\mathbb{C}_{l}$
and $z_{l_{d(1)}}(X^{b})$ $\in\mathbb{C}_{0}\cap\mathbb{C}_{p}$.
From (\ref{eq:4-9}) and (\ref{eq:4-10}), we also see that

\begin{equation}
\int_{a_{0}^{l_{d}}}^{a_{0}^{l_{d(2)}}}\left|z[u_{l_{d*(2)}}(X^{a},\tau)]\right|u'_{l_{d*(2)}}(X^{a},\tau)d\tau\label{eq:4-11}
\end{equation}
\[
>\int_{a_{0}^{l_{d(1)}}}^{a_{0}^{l_{d(2)}}}\left|z[u_{l_{d(2)}}(X^{a},\tau)]\right|u'_{l_{d(2)}}(X^{a},\tau)d\tau
\]
and

\begin{equation}
\int_{a_{0}^{l_{d}}}^{a_{0}^{l_{d(2)}}}\left|z[u_{l_{d*(2)}}(X^{b},\tau)]\right|u'_{l_{d*(2)}}(X^{b},\tau)d\tau>\label{eq:4-12}
\end{equation}
\[
\int_{a_{0}^{l_{d(1)}}}^{a_{0}^{l_{d(2)}}}\left|z[u_{l_{d(2)}}(X^{b},\tau)]\right|u'_{l_{d(2)}}(X^{b},\tau)d\tau
\]
because the multilevel contours $\gamma_{l}(X^{a},t)$ and $\gamma_{l}(X^{b},t)$
whose distances are in (\ref{eq:4-9}) and (\ref{eq:4-10}) has to
pass through the plane $\mathbb{C}_{0}\cap\mathbb{C}_{p}.$ For all
sets of three numbers of the type for $z_{l_{d}}(X^{b})$, $z_{l_{d(1)}}(X^{b})$,
and $z_{l_{d(2)}}(X^{b})$ lying in $\mathbb{C}_{0}\cap\mathbb{C}_{l}$,
$\mathbb{C}_{0}\cap\mathbb{C}_{p}$, and $\mathbb{C}_{p},$ for arbitrary
$\mathbb{C}_{l}$ and $\mathbb{C}_{p}$ in $B_{\mathbb{R}}(\mathbb{C})$,
the equality
\begin{equation}
\begin{array}{c}
\sum_{i=1}^{\infty}\int_{a_{0}^{l_{d}}}^{a_{0}^{l_{d(i)}}}\left|z[u_{l_{d(1)}}(X^{b},\tau)]\right|u'_{l_{d(i)}}(X^{b},\tau)d\tau+\\
+\sum_{i=1}^{\infty}\int_{a_{0}^{l_{d(i)}}}^{a_{0}^{l_{d(i+1)}}}\left|z[u_{l_{d(i+1)}}(X^{b},\tau)]\right|u'_{l_{d(i+1)}}(X^{b},\tau)d\tau\\
=\sum_{i=1}^{\infty}\int_{a_{0}^{l_{d}}}^{a_{0}^{l_{d(i)}}}\left|z[u_{l_{d*(i)}}(X^{b},\tau)]\right|u'_{l_{d*(i)}}(X^{b},\tau)d\tau
\end{array}\label{eq:4-13}
\end{equation}
holds, and (\ref{eq:4-13}) leads to 
\begin{equation}
\sum_{i=1}^{\infty}\int_{a_{0}^{l_{d}}}^{a_{0}^{l_{d(i)}}}\left|z[u_{l_{d*(i)}}(X^{b},\tau)]\right|u'_{l_{d*(i)}}(X^{b},\tau)d\tau>\label{eq:4-14}
\end{equation}
\[
\sum_{i=1}^{\infty}\int_{a_{0}^{l_{d(i)}}}^{a_{0}^{l_{d(i+1)}}}\left|z[u_{l_{d(i+1)}}(X^{b},\tau)]\right|u'_{l_{d(i+1)}}(X^{b},\tau)d\tau.
\]
Consider an arbitrary plane $\mathbb{C}_{q}$ lying somewhere above
$\mathbb{C}_{p}$ and $\mathbb{C}_{p}$ lying somewhere above $\mathbb{C}_{l}$
for $\mathbb{C}_{l}$ , $\mathbb{C}_{p}$, and $\mathbb{C}_{q}$ in
$B_{\mathbb{R}}(\mathbb{C}),$Let the five points (numbers) in the
bundle are arranged as follows: 
\[
z_{l_{d}}(X^{a})\in\mathbb{C}_{0}\cap\mathbb{C}_{l},z_{l_{d(1)}}(X^{a})\in\mathbb{C}_{0}\cap\mathbb{C}_{p},
\]
 
\begin{align*}
z_{l_{d(2)}}(X^{a}) & \notin\mathbb{C}_{0}\cap\mathbb{C}_{p}\text{ and }z_{l_{d(2)}}(X^{a})\in\mathbb{C}_{p},z_{l_{d(3)}}(X^{a})\in\mathbb{C}_{0}\cap\mathbb{C}_{p}\text{ }\\
 & \text{(reachable from \ensuremath{z_{l_{d(2)}}}(\ensuremath{X^{a}})),}
\end{align*}
 
\[
z_{l_{d(4)}}(X^{a})\mathbb{C}_{0}\cap\mathbb{C}_{q}\text{, and \ensuremath{z_{l_{d(5)}}}(\ensuremath{X^{a}})\ensuremath{\notin\mathbb{C}_{0}\cap\mathbb{C}_{q}} and \ensuremath{z_{l_{d(5)}}}(\ensuremath{X^{a}})\ensuremath{\in\mathbb{C}_{p}}.}
\]
 Then the equality arising out of these points is
\begin{align}
\int_{a_{0}^{l_{d}}}^{a_{0}^{l_{d(1)}}}\left|z[u_{l_{d(1)}}(X^{a},\tau)]\right|u'_{l_{d(1)}}(X^{a},\tau)d\tau+\nonumber \\
\int_{a_{0}^{l_{d(1)}}}^{a_{0}^{l_{d(2)}}}\left|z[u_{l_{d(2)}}(X^{a},\tau)]\right|u'_{l_{d(2)}}(X^{a},\tau)d\tau+\nonumber \\
\int_{a_{0}^{l_{d(2)}}}^{a_{0}^{l_{d(3)}}}\left|z[u_{l_{d(3)}}(X^{a},\tau)]\right|u'_{l_{d(3)}}(X^{a},\tau)d\tau+\nonumber \\
\int_{a_{0}^{l_{d(3)}}}^{a_{0}^{l_{d(4)}}}\left|z[u_{l_{d(4)}}(X^{a},\tau)]\right|u'_{l_{d(4)}}(X^{a},\tau)d\tau+\nonumber \\
\int_{a_{0}^{l_{d(4)}}}^{a_{0}^{l_{d(5)}}}\left|z[u_{l_{d(5)}}(X^{a},\tau)]\right|u'_{l_{d(5)}}(X^{a},\tau)d\tau\nonumber \\
=\int_{a_{0}^{l_{d}}}^{a_{0}^{l_{d(5)}}}\left|z[u_{l_{d*(5)}}(X^{a},\tau)]\right|u'_{l_{d*(5)}}(X^{a},\tau)d\tau\label{eq:4-15}
\end{align}
\\
 where the real valued function $u_{l_{d*}}$maps appropriate time
intervals after the parametric representation. For all above such
sets of five points in the bundle and for three sets of planes, we
will have 
\begin{align}
\sum_{i=1}^{\infty}\int_{a_{0}^{l_{d}}}^{a_{0}^{l_{d(i)}}}\left|z[u_{l_{d(i)}}(X^{a},\tau)]\right|u'_{l_{d(i)}}(X^{a},\tau)d\tau+\nonumber \\
\sum_{i=1}^{\infty}\int_{a_{0}^{l_{d(1+i)}}}^{a_{0}^{l_{d(2+i)}}}\left|z[u_{l_{d(2)}}(X^{a},\tau)]\right|u'_{l_{d(2)}}(X^{a},\tau)d\tau+\nonumber \\
+\sum_{i=1}^{\infty}\int_{a_{0}^{l_{d(2)}}}^{a_{0}^{l_{d(3)}}}\left|z[u_{l_{d(3)}}(X^{a},\tau)]\right|u'_{l_{d(3)}}(X^{a},\tau)d\tau+\nonumber \\
\sum_{i=1}^{\infty}\int_{a_{0}^{l_{d(3+i)}}}^{a_{0}^{l_{d(4+i)}}}\left|z[u_{l_{d(4)}}(X^{a},\tau)]\right|u'_{l_{d(4)}}(X^{a},\tau)d\tau+\nonumber \\
\sum_{i=1}^{\infty}\int_{a_{0}^{l_{d(4+i)}}}^{a_{0}^{l_{d(5+i)}}}\left|z[u_{l_{d(5)}}(X^{a},\tau)]\right|u'_{l_{d(5)}}(X^{a},\tau)d\tau\nonumber \\
=\sum_{i=1}^{\infty}\int_{a_{0}^{l_{d}}}^{a_{0}^{l_{d(5+i)}}}\left|z[u_{l_{d*(5+i)}}(X^{a},\tau)]\right|u'_{l_{d*(5+i)}}(X^{a},\tau)d\tau\label{eq:4-16}
\end{align}

A hole is formed in $B_{\mathbb{R}}(\mathbb{C})$ due to a loss of
a set of points (numbers) of data in which a piece of a contour was
located (before it got deleted due to a removal process) or a loss
of a set of data points of a group a contours. The line consisting
of points were lost due to a removal process but all other points
around the line or area around a hole could be chosen by a random
variable. We introduce and define two new sets, namely, a hole and
an island. 
\begin{defn}
\textbf{\label{def:Hole:}Hole}: A hole $H$ is a closed set of points
in $B_{\mathbb{R}}(\mathbb{C})$ such that no point of $H$ is available
to be chosen by an arbitrary $X^{a}\left(z_{l}(t),\mathbb{C}_{l}\right)$
for an arbitrary plane $\mathbb{C}_{l}.$ 
\end{defn}

\begin{example}
Let $\left\{ z_{l}(X^{a},t)\right\} $ be the set of numbers of $\gamma_{l}(X^{a},t)$
for $t\in[t_{0},t_{b}${]} that got deleted due to a removal process.
Then $X^{a}\left(z_{l}(t),\mathbb{C}_{l}\right)$ would not be able
to choose a number from $\left\{ z_{l}(X^{a},t)\right\} $ for $t\in[t_{0},t_{b}${]}.
There could be many such holes in $B_{\mathbb{R}}(\mathbb{C}).$ Two
or more contours using the same set of points for a period of time,
then the removal process of one contour could delete the common set
of data so that a hole is formed. We will soon see that the space
created by the set $H$ is dynamic. 
\end{example}

\begin{defn}
\textbf{\label{def:Island:}Island: }Let $S\subset B_{\mathbb{R}}(\mathbb{C})$
and $S\neq\phi\,\text{(empty).}$ The set $S$ is called an island
if any element in $S$ is available to be chosen by a random variable
but no contour can be drawn from an element within $S$ to an element
outside $S,$ say, $S^{c}.$ Here $S,S^{c}\in B_{\mathbb{R}}(\mathbb{C}).$ 
\end{defn}

\begin{example}
Consider a disc $D\left(z_{l_{d}}(X^{a}),r_{d}(X^{a})\right)$ with
a center $z_{l_{d}}(X^{a})$ chosen by $X^{a}\left(z_{l}(t),\mathbb{C}_{l}\right)$
from a previous iteration. Let two contours $\gamma_{l}(X^{b},t)$
and $\gamma_{l}(X^{c},t)$ pass through $D\left(z_{l_{d}}(X^{a}),r_{d}(X^{a})\right)$
and intersects at two locations (say, $z_{l_{i}}(X^{b})$ and $z_{l_{j}}(X^{b})$
as shown in Figure \ref{fig:An-island-is}. Suppose the spaces of
points of $\gamma_{l}(X^{b},t)$ and $\gamma_{l}(X^{c},t)$ passing
through $D\left(z_{l_{d}}(X^{a}),r_{d}(X^{a})\right)$ were lost due
to respective variables' removal processes. Then the space formed
between these two points of intersections including the data on the
contours between $z_{l_{i}}(X^{b})$ and $z_{l_{j}}(X^{b})$ is an
island. 

\begin{figure}
\includegraphics[scale=0.5]{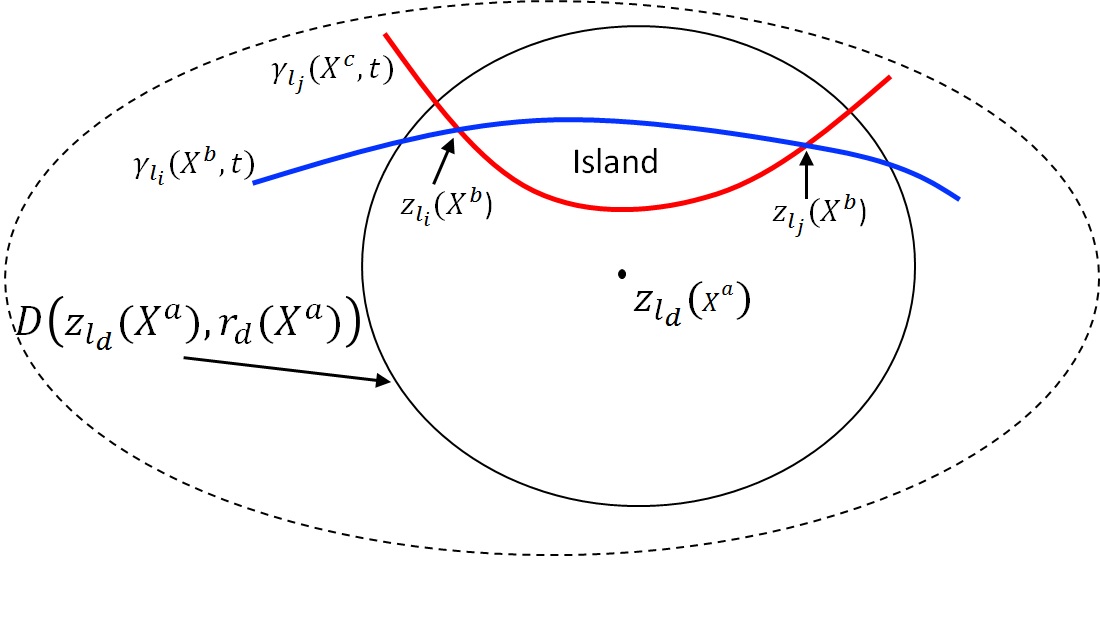}

\caption{\label{fig:An-island-is}An island is created within $D\left(z_{l_{d}}(X^{a}),r_{d}(X^{a})\right).$}
\end{figure}

The two sets $H$ and $S$ are dynamic as spaces created by these
sets could change because of the dynamic nature of contour formation
and removal process described in section 3. A time-dependent versions
of the definitions for \emph{holes} and \emph{islands }can be given
here. A set $S(t)$ for $t\in[t_{0},\infty)$ and satisfying the definition
\ref{def:Island:} can be called an island at $t.$ The set of elements
of $S(t_{c})$ for $t_{c}\in[t_{0},\infty)$ satisfying the Definition
\ref{def:Island:} might lose all its elements in a removal process
and might turn into a \emph{hole} at a time $t_{d}$ for $t_{d}>t_{c}.$
If $S(t)$is an island then no one cannot draw a contour from the
elements of $S(t)$ to an element in $S^{c}(t)$ where 
\begin{equation}
S^{c}(t)=\left\{ z:z\notin S(t)\subset B_{\mathbb{R}}(\mathbb{C})\text{ and }z\in B_{\mathbb{R}}(\mathbb{C})\right\} .\label{eq:Scompliment}
\end{equation}
Similarly, a contour cannot be drawn from an element (point) of $S^{c}(t)$
to an element in $S(t).$ The spaces of $S(t)$ and $S^{c}(t)$ are
separated by $H.$ The area of a hole could change over a time as
more data points removed are added to a specific hole. 
\end{example}

\begin{thm}
\label{thm:ISLAND}Suppose a disc $D\left(z_{l_{d}}(X^{a}),r_{d}(X^{a})\right)$
formed out of $X^{a}\left(z_{l}(t),\mathbb{C}_{l}\right)$ is given.
Let $z_{l_{1}}$ and $z_{l_{2}}$ be two points in the boundary set
of $D.$ A contour $C_{1}$ formed by an arbitrary process $X^{b}\left(z_{l}(t),\mathbb{C}_{l}\right)$
enters the disc through $z_{l_{1}}$ and leaves the disc from $z_{l_{2}}$
and another contour $C_{2}$ formed by an arbitrary process $X^{c}\left(z_{l}(t),\mathbb{C}_{l}\right)$
enters the disc through $z_{l_{1}}$ and leaves the disc from $z_{l_{2}}.$
The paths of $C_{1}$ and $C_{2}$ never meet except at the points
$z_{l_{1}}$ and $z_{l_{2}}$ and the center $z_{l_{d}}(X^{a})$ lies
in between the contours. Suppose the removal process of two contours
$C_{1}$ and $C_{2}$ introduced such that $C_{1}\cup C_{2}$ form
a hole at $t$. Then, the set of points lying between $C_{1}$ and
$C_{2}$ forms an island. 
\end{thm}

\begin{proof}
Given that $z_{l_{1}}$ and $z_{l_{2}}$are located on the boundary
of the disc $D\left(z_{l_{d}}(X^{a}),r_{d}(X^{a})\right)$ and $z_{l_{d}}(X^{a})$
is located in between $C_{1}$ and $C_{2}.$ See Figure \ref{fig:Formation-of-ISLAND}
for a description of given information and locations of $z_{l_{1}}$
and $z_{l_{2}}.$ 

\begin{figure}
\includegraphics[scale=0.5]{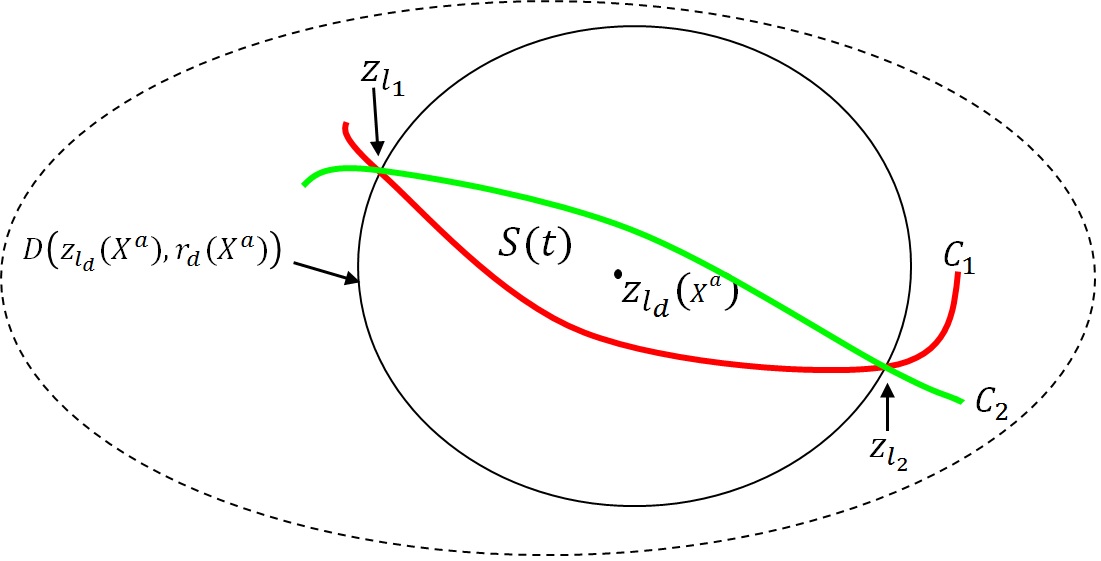}

\caption{\label{fig:Formation-of-ISLAND}Formation of an island due to a removal
process of Theorem \ref{thm:ISLAND}.}

\end{figure}

Let $S(t)$ be a set in $D\left(z_{l_{d}}(X^{a}),r_{d}(X^{a})\right)$
that consists of all points in between $C_{1}$ and $C_{2}$ as shown
in Figure \ref{fig:Formation-of-ISLAND}. The set $D\left(z_{l_{d}}(X^{a}),r_{d}(X^{a})\right)\backslash S(t)$
will consists of points as in (\ref{eq:4-19}), 
\begin{align}
D\left(z_{l_{d}}(X^{a}),r_{d}(X^{a})\right)\backslash S(t)=\left\{ z_{l}:z_{l}\in D\left(z_{l_{d}}(X^{a}),r_{d}(X^{a})\right)\text{ and }z_{l}\notin S(t)\right\} .\nonumber \\
\label{eq:4-19}
\end{align}
The disc $D\left(z_{l_{d}}(X^{a}),r_{d}(X^{a})\right)$ can be partitioned
into disjoint union of three sets as below
\begin{align}
D\left(z_{l_{d}}(X^{a}),r_{d}(X^{a})\right)=\left[D\left(z_{l_{d}}(X^{a}),r_{d}(X^{a})\right)\backslash S(t)\right]\cup S(t)\cup\left[C_{1}\cup C_{2}\right].\nonumber \\
\label{eq:4-20}
\end{align}
By the construction, we cannot draw a contour from a point in $S(t)$
to a point in 
\[
D\left(z_{l_{d}}(X^{a}),r_{d}(X^{a})\right)\backslash S(t).
\]
Hence, $S(t)$ is an island. 
\end{proof}
\begin{thm}
The union of collection of all holes within $B_{\mathbb{R}}(\mathbb{C})$
is compact.
\end{thm}

\begin{proof}
Each hole is a closed set of points from a contour or a collection
of contours. Let $H_{\alpha}$ be an arbitrary hole. An infinite union
of such holes,
\begin{equation}
\bigcup_{\alpha}H_{\alpha}.\label{eq:4-21}
\end{equation}
 is closed. Each $H_{\alpha}$ at time $t$ is bounded by the length
of the contour formed until the time $t+1$, because the removal process
follows contour formation with some lag in the time. So
\begin{equation}
\left|\bigcup_{\alpha}H_{\alpha}\right|<\int_{a_{0}^{l_{d}}}^{a_{0}^{l_{d(t+i)}}}\left|z[u_{l_{d(t+i)}}(X^{a},\tau)]\right|u'_{l_{d(5+i)}}(X^{a},\tau)d\tau,\label{4-22}
\end{equation}
 for an arbitrary process $X^{a}\left(z_{l}(t),\mathbb{C}_{l}\right)$.
So the union of a collection of holes is bounded. Hence, such a collection
is compact. 
\end{proof}
\begin{thm}
Suppose infinitely many random variables of the type $X^{a}\left(z_{l}(t),\mathbb{C}_{l}\right)$
are introduced one by one as in Section 3. If every point within an
island $S(t)$ is part of a contour at time $t$, then due to a removal
process $S(t)$ becomes $H(t_{c})$ for $t_{c}>t.$
\end{thm}

\begin{proof}
Given there is a $S(t)$ at time $t.$ Let there be finitely many
contours passing through the region $S(t)$ such that all the points
of $S(t)$ are in one or more of the contours. Once a removal process
is introduced at $t_{c}>t$ then all the points of $S(t)$ will not
be available for a new random variable. Hence, $S(t)$ will asymptotically
become a hole. 
\end{proof}
Let us introduce a removal process for the infinitely many contours
introduced earlier in this section. These contours were all started
at the same time $t_{0}$. The number of such contours is equivalent
to the set of elements in $\mathbb{C}_{l}.$ Some of the elements
in $\mathbb{C}_{l}$ are also in $\mathbb{C}_{0}\cap\mathbb{C}_{l}.$
Note that,
\begin{equation}
\mathbb{C}_{l}=\left(\mathbb{C}_{0}\cap\mathbb{C}_{l}\right)\cup\mathbb{C}_{l}=\left\{ z_{l}:z_{l}\in\mathbb{C}_{l}\text{ or }z_{l}\in\mathbb{C}_{0}\cap\mathbb{C}_{l}\right\} .\label{eq:4-23}
\end{equation}
 So removing a contour data that was formed from $t_{0}$ until a
time $t_{c}$ at a rate $\phi(X^{a},t)$ for an arbitrary random variable
$X^{a}\left(z_{l}(t),\mathbb{C}_{l}\right)$ would also remove contour
data that is located in the set $\mathbb{C}_{0}\cap\mathbb{C}_{l}.$
We described earlier how the islands of sets of data could be formed,
and formation of islands could happen in distinct time intervals once
removal process is introduced. 

Remember that all the contours have origin in $\mathbb{C}_{l}$ only.
Now at $t_{c}$ due to a removal process, the points of $\mathbb{C}_{l}$
are not available to be chosen by any of the infinitely many random
variables. The paths of these random variables are not disjoint. Some
of these random variables may be in another plane outside $\mathbb{C}_{l}$at
$t=t_{c}.$ Contour formations of these infinitely many random variables
of type $X^{a}\left(z_{l}(t),\mathbb{C}_{l}\right)$ may continue
even after $t_{c}$ because of their presence in outside $\mathbb{C}_{l}$
at $t=t_{c}.$ So the entire plane $\mathbb{C}_{l}$ is not available
in the bundle $B_{\mathbb{R}}(\mathbb{C})$. So the remaining space
formed will be the set of points $B_{\mathbb{R}}(\mathbb{C})\backslash\mathbb{C}_{l},$
where
\begin{equation}
B_{\mathbb{R}}(\mathbb{C})\backslash\mathbb{C}_{l}=\left\{ z_{l}:z_{l}\in B_{\mathbb{R}}(\mathbb{C})\text{ and }z_{l}\in\mathbb{C}_{l}\right\} \label{eq:4-24}
\end{equation}

\subsection{Consequences of $B_{\mathbb{R}}(\mathbb{C})\backslash\mathbb{C}_{l}$
on Multilevel Contours}

At the time of $t=t_{c},$ there could be one of the following two
situations for the status of infinitely many random variables that
were introduced at $t_{0}.$

$(i)$ All the contours formed until $t_{c}$ are active only in $\mathbb{C}_{l}$
at $t=t_{c}$ and for any arbitrary $X^{a}\left(z_{l}(t),\mathbb{C}_{l}\right)$
the set of points $\left\{ z_{l}\left(X^{a},t\right)\right\} $ for
$t\in[t_{0},t_{c}]$ are in $\mathbb{C}_{l},$ i.e.
\begin{equation}
\left\{ z_{l}\left(X^{a},t\right)\right\} \text{}\left(t\in[t_{0},t_{c}]\right)\subset\mathbb{C}_{l}.\label{eq:4-25}
\end{equation}

This means, no contour has crossed the plane until $t=t_{c}.$

$(ii)$ Only a fraction of infinitely many random variables are active
in $\mathbb{C}_{l}$ at $t=t_{c},$ and the rest of all random variables
are active outside $\mathbb{C}_{l}.$ 

If every contour satisfies the above situation $(i)$ at $t=t_{c},$
then due to the removal process all the points in $\mathbb{C}_{l}$cannot
be reached to form a contour. In such a situation, the plane $\mathbb{C}_{l}$
forms a hole, and the contours formation process halts. Suppose only
a fraction of infinitely many random variables, say, $\alpha_{1}$
are inside $\mathbb{C}_{l}$ at $t_{c}$, then there will be two options
for the location of those random variables (i.e. for the fraction
$(1-\alpha_{1}$) which are outside $\mathbb{C}_{l}$ at $t_{c}.$
A fraction of them, say, $\alpha_{2}$ will be above $\mathbb{C}_{l}$
and the remaining fraction $(1-\alpha_{2})$ of the random variables
will be somewhere in a plane below $\mathbb{C}_{l}.$ Let $T$ be
the set of all contours active $t_{c}$ and $\alpha_{1}$ be the fraction
of them active and located in $\mathbb{C}_{l},$ then $(1-\alpha_{1})\left|T\right|$
be the set of contours that are active and outside $\mathbb{C}_{l}.$
Then, $(1-\alpha_{1})\left|T\right|$ are further divided as 
\begin{equation}
(1-\alpha_{1})\left|T\right|=\alpha_{2}\left[(1-\alpha_{1})\left|T\right|\right]+(1-\alpha_{2})\left[(1-\alpha_{1})\left|T\right|\right],\label{eq:4-26}
\end{equation}
 where $\alpha_{2}$ is the fraction of contours that are active at
$t_{c}$ and are located somewhere in a plane above $\mathbb{C}_{l},$
and $1-\alpha_{2}$ is the fraction of contours that are active at
$t_{c}$ and are located somewhere in a plane below $\mathbb{C}_{l}.$
Because the plane $\mathbb{C}_{l}$ became a\emph{ hole}, the set
of contours, say, $T_{\alpha_{1}}$denoted to represent\emph{ $\alpha_{1}\left|T\right|$
}number of contours will stop further formation for $t>t_{c}.$ The
set of contours, say, $T_{1-\alpha_{1}}$ will be active outside $\mathbb{C}_{l}$,
where
\begin{equation}
T_{1-\alpha_{1}}=T_{\alpha_{2}}\cup T_{1-\alpha_{2}}.\label{eq:4-27}
\end{equation}
Here $T_{\alpha_{2}}$ represent the set of contours that are located
somewhere in a plane above $\mathbb{C}_{l},$ $T_{1-\alpha_{2}}$
represent the set of contours that are located somewhere in a plane
below $\mathbb{C}_{l}.$ The carnality of $T_{1-\alpha_{1}}$is constant
for $t>t_{c},$ so as the cardinals of $T_{\alpha_{2}}$ and $T_{1-\alpha_{2}}$
at $t>t_{c}.$We also partition the set of elements in $B_{\mathbb{R}}(\mathbb{C})$
at $t>t_{c}$ as

\begin{equation}
\left.B_{\mathbb{R}}(\mathbb{C})\right|_{t>t_{c}}=B_{\mathbb{R}}(\mathbb{C},\alpha_{2})\cup B_{\mathbb{R}}(\mathbb{C},1-\alpha_{2}),\label{eq:4-28}
\end{equation}
where as
\begin{equation}
B_{\mathbb{R}}(\mathbb{C})=B_{\mathbb{R}}(\mathbb{C},\alpha_{2})\cup\mathbb{C}_{l}\cup B_{\mathbb{R}}(\mathbb{C},1-\alpha_{2}),\label{eq:4-29}
\end{equation}
where $B_{\mathbb{R}}(\mathbb{C},\alpha_{2})\subset B_{\mathbb{R}}(\mathbb{C})$
is the set of planes which are above $\mathbb{C}_{l},$ and $B_{\mathbb{R}}(\mathbb{C},1-\alpha_{2})\subset B_{\mathbb{R}}(\mathbb{C})$
is the set of planes which are above $\mathbb{C}_{l}.$ The set of
contours that are active after $t_{c}$ are located in the set, say
$B_{\mathbb{R}}(\mathbb{C},1-\alpha_{1})$ and from (\ref{eq:4-28}),
we have
\begin{equation}
B_{\mathbb{R}}(\mathbb{C},1-\alpha_{2})=\left.B_{\mathbb{R}}(\mathbb{C})\right|_{t>t_{c}}.\label{eq:4-30}
\end{equation}
As described above, contour formations and removal processes of the
contours in $B_{\mathbb{R}}(\mathbb{C},1-\alpha_{2})$ will continue.
Due to presence of the \emph{hole} $\mathbb{C}_{l},$ the active contours
of $B_{\mathbb{R}}(\mathbb{C},\alpha_{2})$ and $B_{\mathbb{R}}(\mathbb{C},1-\alpha_{2})$
will not have any further intersecting points. The tails of the contours
remaining in $B_{\mathbb{R}}(\mathbb{C},\alpha_{2})$ and $B_{\mathbb{R}}(\mathbb{C},1-\alpha_{2})$
will be eventually lost for some time after $t_{c}.$ Let $\gamma_{l}(X^{a},t)$
be an arbitrary contour that is active in $B_{\mathbb{R}}(\mathbb{C},\alpha_{2})$,
and was created by $X^{a}.$ Suppose the set of points touched by
the contour $\gamma_{l}(X^{a},t)$ prior to $t_{c}$ were located
in $B_{\mathbb{R}}(\mathbb{C},\alpha_{2})$, $\mathbb{C}_{l}$, and
$B_{\mathbb{R}}(\mathbb{C},1-\alpha_{2}),$ this contour is described
by $z_{l}(X^{a},t)$ $\left(t\in[t_{0},\infty)\right)$ and $t=u_{l_{c}}(X^{a},\tau)$
$\left(a_{0}\leq\tau\leq a_{c}\right)$ is the parametric representation
for $\gamma_{l}(X^{a},t)$ with a real values function $u_{l}(X^{a},\tau)$
mapping $[a_{0},a_{c}]$ onto the interval $[t_{0},t_{c}].$ Let $L\left(\gamma_{l}(X^{a},t)\left(t\in[t_{0},t_{c}]\right)\right)$
represent the length of $\gamma_{l}(X^{a},t)\left(t\in[t_{0},t_{c}]\right)$
up to $t_{c},$ then

\begin{equation}
L\left(\gamma_{l}(X^{a},t)\left(t\in[t_{0},t_{c}]\right)\right)=\int_{a_{0}}^{a_{c}}\left|z[u_{l_{c}}(X^{a},\tau)]\right|u'_{l_{c}}(X^{a},\tau)d\tau\label{eq:4-31}
\end{equation}
had covered points from each disjoint set of $B_{\mathbb{R}}(\mathbb{C})$
in (\ref{eq:4-29}). Let us assume that $\gamma_{l}(X^{a},t)\left(t\in[t_{0},t_{c}]\right)$
has visited a multiple number of times through each of the sets of
(\ref{eq:4-29}) before in remained active in $B_{\mathbb{R}}(\mathbb{C},\alpha_{2})$
at $t=t_{c}.$ Then $L\left(\gamma_{l}(X^{a},t)\left(t\in[t_{0},t_{c}]\right)\right)$
in (\ref{eq:4-31}) can be expressed as three components where each
component is made up of several contour integrals. Since $\gamma_{l}(X^{a},t)$
had visited each portion in (\ref{eq:4-31}) several times, the length
$L\left(\gamma_{l}(X^{a},t)\left(t\in[t_{0},t_{c}]\right)\right)$
in (\ref{eq:4-31}) is distributed into corresponding parts.The first
part consists of the sum of all the lengths of piecewise contours
of (\ref{eq:4-31}) lying in $\mathbb{C}_{l}$, say, $L\left(\gamma_{l}(X^{a},t,\alpha_{1})\left(t\in[t_{0},t_{c}]\right)\right)$,
and can be computed using
\begin{align}
L\left(\gamma_{l}(X^{a},t,\alpha_{1})\left(t\in[t_{0},t_{c}]\right)\right) & =\int_{a(0)}^{a(1)}\left|z[u_{l_{a(1)}}(X^{a},\tau)]\right|u'_{l_{a(1)}}(X^{a},\tau)d\tau+\nonumber \\
 & \sum_{i\in A(\alpha_{1})}\int_{a(i)}^{a(i+1)}\left|z[u_{l_{a(i+1)}}(X^{a},\tau)]\right|u'_{l_{a(i+1)}}(X^{a},\tau)d\tau,\label{eq:4-32}
\end{align}
where $u_{l_{a(1)}}(X^{a},\tau)$ and $u_{l_{a(i+1)}}(X^{a},\tau)$
are the real valued functions used in parametric representations with
corresponding onto mappings. The notation $i\in A(\alpha_{1})$ indicates
summing the length over all the piecewise contours in the set $A(\alpha_{1}).$
The set $A(\alpha_{1})$ consists of all the piecewise contours of
$X^{a}\left(z_{l}(t),\mathbb{C}_{l}\right)$ until $t_{c}$ that are
lying in $\mathbb{C}_{l}.$ The first integral on the R.H.S of (\ref{eq:4-32})
is the length of the piecewise contour from its origin to the entry
point either in $B_{\mathbb{R}}(\mathbb{C},\alpha_{2})$ or in $B_{\mathbb{R}}(\mathbb{C},1-\alpha_{2})$.
The sum of integrals on the R.H.S of (\ref{eq:4-32}) is the total
length of the piecewise contours due to $A(\alpha_{1}).$ The second
part in (\ref{eq:4-31}) consists of piecewise contours in $B_{\mathbb{R}}(\mathbb{C},\alpha_{2})$
whose total length, say, $L\left(\gamma_{l}(X^{a},t,\alpha_{2})\left(t\in[t_{0},t_{c}]\right)\right)$,
is computed as
\begin{align}
L\left(\gamma_{l}(X^{a},t,\alpha_{2})\left(t\in[t_{0},t_{c}]\right)\right) & =\sum_{i\in A(\alpha_{2})}\int_{a(i)}^{a(i+1)}\left|z[u_{l_{a(i+1)}}(X^{a},\tau)]\right|u'_{l_{a(i+1)}}(X^{a},\tau)d\tau+\nonumber \\
 & \int_{a_{c-1}}^{a_{c}}\left|z[u_{l_{c}}(X^{a},\tau)]\right|u'_{l_{c}}(X^{a},\tau)d\tau\label{eq:4-33}
\end{align}
The set $A(\alpha_{2})$ in (\ref{eq:4-33}) consists of all the piecewise
contours of $X^{a}\left(z_{l}(t),\mathbb{C}_{l}\right)$ until $t_{c}$
that are lying in $B_{\mathbb{R}}(\mathbb{C},\alpha_{2}).$ The second
integral in the R.H.S. of (\ref{eq:4-33}) consists of length of the
last piece of the contour $\gamma_{l}(X^{a},t)$ until $t=t_{c}$
in $B_{\mathbb{R}}(\mathbb{C},\alpha_{2}).$ Here $u_{l_{c}}(X^{a},\tau)$
is the real valued function used in parametric representations with
corresponding onto mappings. The third part in (\ref{eq:4-31}) consists
of piecewise contours in $B_{\mathbb{R}}(\mathbb{C},1-\alpha_{2})$
whose total length, say, $L\left(\gamma_{l}(X^{a},t,1-\alpha_{2})\left(t\in[t_{0},t_{c}]\right)\right)$,
is computed as
\begin{equation}
\begin{array}{c}
L\left(\gamma_{l}(X^{a},t,1-\alpha_{2})\left(t\in[t_{0},t_{c}]\right)\right)\\
\sum_{i\in A(1-\alpha_{2})}\int_{a(i)}^{a(i+1)}\left|z[u_{l_{a(i+1)}}(X^{a},\tau)]\right|u'_{l_{a(i+1)}}(X^{a},\tau)d\tau.
\end{array}\label{eq:4-34}
\end{equation}
Hence the length in (\ref{eq:4-31}) can be expressed using (\ref{eq:4-32}),
(\ref{eq:4-33}), and (\ref{eq:4-34}) as
\begin{align}
L\left(\gamma_{l}(X^{a},t)\left(t\in[t_{0},t_{c}]\right)\right) & =\int_{a(0)}^{a(1)}\left|z[u_{l_{a(1)}}(X^{a},\tau)]\right|u'_{l_{a(1)}}(X^{a},\tau)d\tau+\nonumber \\
 & \sum_{i\in A(\alpha_{1})}\int_{a(i)}^{a(i+1)}\left|z[u_{l_{a(i+1)}}(X^{a},\tau)]\right|u'_{l_{a(i+1)}}(X^{a},\tau)d\tau+\nonumber \\
 & =\sum_{i\in A(\alpha_{2})}\int_{a(i)}^{a(i+1)}\left|z[u_{l_{a(i+1)}}(X^{a},\tau)]\right|u'_{l_{a(i+1)}}(X^{a},\tau)d\tau+\nonumber \\
 & \int_{a_{c-1}}^{a_{c}}\left|z[u_{l_{c}}(X^{a},\tau)]\right|u'_{l_{c}}(X^{a},\tau)d\tau+\nonumber \\
 & \sum_{i\in A(1-\alpha_{2})}\int_{a(i)}^{a(i+1)}\left|z[u_{l_{a(i+1)}}(X^{a},\tau)]\right|u'_{l_{a(i+1)}}(X^{a},\tau)d\tau.\label{eq:4-35}
\end{align}
Due to the hole $\mathbb{C}_{l}$ created in the bundle $B_{\mathbb{R}}(\mathbb{C})$,
the remaining length of the contour that will be subjected to removal
process is obtained by removing the sum of piecewise contour lengths
in $\mathbb{C}_{l},$ and is given by

\begin{align}
L\left(\gamma_{l}(X^{a},t)\left(t\in[t_{0},t_{c}]\right)\right) & =\sum_{i\in A(\alpha_{2})}\int_{a(i)}^{a(i+1)}\left|z[u_{l_{a(i+1)}}(X^{a},\tau)]\right|u'_{l_{a(i+1)}}(X^{a},\tau)d\tau+\nonumber \\
 & \int_{a_{c-1}}^{a_{c}}\left|z[u_{l_{c}}(X^{a},\tau)]\right|u'_{l_{c}}(X^{a},\tau)d\tau+\nonumber \\
 & \sum_{i\in A(1-\alpha_{2})}\int_{a(i)}^{a(i+1)}\left|z[u_{l_{a(i+1)}}(X^{a},\tau)]\right|u'_{l_{a(i+1)}}(X^{a},\tau)d\tau.\label{eq:4-36}
\end{align}
Since $\gamma_{l}(X^{a},t)$ is active in $B_{\mathbb{R}}(\mathbb{C},\alpha_{2})$,
the formation of the contour will continue forever and the sum of
the pieces of the lengths of $\gamma_{l}(X^{a},t)$ that is there
in $B_{\mathbb{R}}(\mathbb{C},1-\alpha_{2})$ will be deleted from
$B_{\mathbb{R}}(\mathbb{C},1-\alpha_{2}).$ This deletion could be
according to a removal function $\phi(X^{a},t,1-\alpha_{2})$ similar
to the procedure explained in (\ref{eq:phiinintegral}). The differential
equation to model the space of points lost due to a removal process
$\phi(X^{a},t,1-\alpha_{2})$ is 
\[
\]
\begin{align}
\left.\frac{dB_{\mathbb{R}}(\mathbb{C},1-\alpha_{2})}{dt}\right|_{\phi(X^{a},t,1-\alpha_{2})}=\left.B_{\mathbb{R}}(\mathbb{C},1-\alpha_{2})\right|_{t=t_{c}}-\phi(X^{a},t)\left.B_{\mathbb{R}}(\mathbb{C},1-\alpha_{2})\right|_{t=t_{c}}\nonumber \\
\label{eq:4-37}
\end{align}
for 
\begin{align}
\phi(X^{a},t,1-\alpha_{2})=\psi(X^{a},t,1-\alpha_{2})\times\nonumber \\
\sum_{i\in A(1-\alpha_{2})}\int_{a(i)}^{a(i+1)}\left|z[u_{l_{a(i+1)}}(X^{a},\tau)]\right|u'_{l_{a(i+1)}}(X^{a},\tau)d\tau\text{ \ensuremath{\left(t\in(t_{c},t_{c'}]\right)}}.\nonumber \\
\label{eq:4-38}
\end{align}
for $0<\psi(X^{a},t,1-\alpha_{2})<1.$ Suppose $\gamma_{l}(X^{a},t)$
is active in $B_{\mathbb{R}}(\mathbb{C},1-\alpha_{2})$ instead of
in $B_{\mathbb{R}}(\mathbb{C},\alpha_{2}).$ As in above, the set
of points touched by the contour $\gamma_{l}(X^{a},t)$ prior to $t_{c}$
would have located in $B_{\mathbb{R}}(\mathbb{C},\alpha_{2})$, $\mathbb{C}_{l}$,
and $B_{\mathbb{R}}(\mathbb{C},1-\alpha_{2}).$ We described this
contour by $z_{l}(X^{a},t)$ $\left(t\in[t_{0},\infty)\right)$ and
$t=w_{l_{c}}(X^{a},\tau)$ $\left(a'_{0}\leq\tau\leq a'_{c}\right)$
is the parametric representation for $\gamma_{l}(X^{a},t)$ with a
real values function $w_{l_{c}}(X^{a},\tau)$ mapping $[a'_{0},a'_{c}]$
onto the interval $[t_{0},t_{c}].$ Let $L\left(\gamma_{l}(X^{a},t,1-\alpha_{2})\left(t\in[t_{0},t_{c}]\right)\right)$
represent the length of $\gamma_{l}(X^{a},t)\left(t\in[t_{0},t_{c}]\right)$
up to $t_{c},$ then
\begin{align}
L\left(\gamma_{l}(X^{a},t,1-\alpha_{2})\left(t\in[t_{0},t_{c}]\right)\right)=\int_{a'_{0}}^{a'_{c}}\left|z[w_{l_{c}}(X^{a},\tau)]\right|w'_{l_{c}}(X^{a},\tau)d\tau\nonumber \\
\label{eq:4-39}
\end{align}

had covered points from each disjoint set of $B_{\mathbb{R}}(\mathbb{C})$
in (\ref{eq:4-29}). As in above, let us assume that $\gamma_{l}(X^{a},t)\left(t\in[t_{0},t_{c}]\right)$
has visited a multiple number of times through each of the sets of
(\ref{eq:4-29}) before in remained active in $B_{\mathbb{R}}(\mathbb{C},1-\alpha_{2})$
at $t=t_{c}.$ Then $L\left(\gamma_{l}(X^{a},t,1-\alpha_{2})\left(t\in[t_{0},t_{c}]\right)\right)$
in (\ref{eq:4-39}) can be expressed as three components where each
component is made up of several contour integrals. Since $\gamma_{l}(X^{a},t)$
had visited each portion in (\ref{eq:4-39}) several times, the length
$L\left(\gamma_{l}(X^{a},t,1-\alpha_{2})\left(t\in[t_{0},t_{c}]\right)\right)$
in (\ref{eq:4-39}) is distributed into corresponding parts.The first
part consists of the sum of all the lengths of piecewise contours
of (\ref{eq:4-39}) lying in $\mathbb{C}_{l}$, say, $L\left(\gamma_{l}(X^{a},t,\alpha_{1})\left(t\in[t_{0},t_{c}]\right)\right)$,
and can be computed using
\begin{align}
L\left(\gamma_{l}(X^{a},t,\alpha_{1})\left(t\in[t_{0},t_{c}]\right)\right) & =\int_{a'(0)}^{a'(1)}\left|z[w_{l_{a'(1)}}(X^{a},\tau)]\right|w'_{l_{a'(1)}}(X^{a},\tau)d\tau+\nonumber \\
 & \sum_{i\in A'(\alpha_{1})}\int_{a'(i)}^{a'(i+1)}\left|z[w_{l_{a(i+1)}}(X^{a},\tau)]\right|w'_{l_{a(i+1)}}(X^{a},\tau)d\tau,\label{eq:4-40}
\end{align}
where $w_{l_{a'(1)}}(X^{a},\tau)$ and $w_{l_{a'(i+1)}}(X^{a},\tau)$
are the real valued functions used in parametric representations with
corresponding onto mappings. The notation $i\in A'(\alpha_{1})$ indicates
summing the length over all the piecewise contours in the set $A'(\alpha_{1}).$
The set $A'(\alpha_{1})$ consists of all the piecewise contours of
$X^{a}\left(z_{l}(t),\mathbb{C}_{l}\right)$ until $t_{c}$ that are
lying in $\mathbb{C}_{l}.$ The first integral on the R.H.S of (\ref{eq:4-40})
is the length of the piecewise contour from its origin to the entry
point either in $B_{\mathbb{R}}(\mathbb{C},\alpha_{2})$ or in $B_{\mathbb{R}}(\mathbb{C},1-\alpha_{2})$.
The sum of integrals on the R.H.S of (\ref{eq:4-40}) is the total
length of the piecewise contours due to $A'(\alpha_{1}).$ The second
part in (\ref{eq:4-39}) consists of piecewise contours in $B_{\mathbb{R}}(\mathbb{C},1-\alpha_{2})$
whose total length, say, $L\left(\gamma_{l}(X^{a},t,1-\alpha_{2})\left(t\in[t_{0},t_{c}]\right)\right)$,
is computed as
\begin{equation}
\begin{array}{c}
L\left(\gamma_{l}(X^{a},t,1-\alpha_{2})\left(t\in[t_{0},t_{c}]\right)\right)=\\
\sum_{i\in A'(1-\alpha_{2})}\int_{a'(i)}^{a'(i+1)}\left|z[w_{l_{a'(i+1)}}(X^{a},\tau)]\right|w'_{l_{a'(i+1)}}(X^{a},\tau)d\tau+\\
\int_{a'_{c-1}}^{a'_{c}}\left|z[w_{l_{c}}(X^{a},\tau)]\right|w'_{l_{c}}(X^{a},\tau)d\tau
\end{array}\label{eq:4-41}
\end{equation}
Note that the contour is active in $B_{\mathbb{R}}(\mathbb{C},1-\alpha_{2})$.
The set $A'(\alpha_{2})$ in (\ref{eq:4-41}) consists of all the
piecewise contours of $X^{a}\left(z_{l}(t),\mathbb{C}_{l}\right)$
until $t_{c}$ that are lying in $B_{\mathbb{R}}(\mathbb{C},1-\alpha_{2}).$
The second integral in the R.H.S. of (\ref{eq:4-41}) consists of
length of the last piece of the contour $\gamma_{l}(X^{a},t)$ until
$t=t_{c}$ in $B_{\mathbb{R}}(\mathbb{C},1-\alpha_{2}).$ Here $w_{l_{c}}(X^{a},\tau)$
is the real valued function used in parametric representations with
corresponding onto mappings. The third part in (\ref{eq:4-39}) consists
of piecewise contours in $B_{\mathbb{R}}(\mathbb{C},\alpha_{2})$
whose total length, say, $L\left(\gamma_{l}(X^{a},t,\alpha_{2})\left(t\in[t_{0},t_{c}]\right)\right)$,
is computed as
\begin{align}
L\left(\gamma_{l}(X^{a},t,\alpha_{2})\left(t\in[t_{0},t_{c}]\right)\right) & =\nonumber \\
 & \sum_{i\in A'(1-\alpha_{2})}\int_{a(i)}^{a(i+1)}\left|z[w_{l_{a'(i+1)}}(X^{a},\tau)]\right|w'_{l_{a'(i+1)}}(X^{a},\tau)d\tau.\nonumber \\
\label{eq:4-42}
\end{align}
The length in (\ref{eq:4-39}) can be expressed using (\ref{eq:4-40}),
(\ref{eq:4-41}), and (\ref{eq:4-42}) as
\begin{align}
L\left(\gamma_{l}(X^{a},t,1-\alpha_{2})\left(t\in[t_{0},t_{c}]\right)\right) & =\nonumber \\
\int_{a'(0)}^{a'(1)}\left|z[w_{l_{a'(1)}}(X^{a},\tau)]\right|w'_{l_{a'(1)}}(X^{a},\tau)d\tau+\nonumber \\
\sum_{i\in A'(\alpha_{1})}\int_{a'(i)}^{a'(i+1)}\left|z[w_{l_{a(i+1)}}(X^{a},\tau)]\right|w'_{l_{a(i+1)}}(X^{a},\tau)d\tau+\nonumber \\
\sum_{i\in A'(1-\alpha_{2})}\int_{a'(i)}^{a'(i+1)}\left|z[w_{l_{a'(i+1)}}(X^{a},\tau)]\right|w'_{l_{a'(i+1)}}(X^{a},\tau)d\tau+\nonumber \\
\int_{a'_{c-1}}^{a'_{c}}\left|z[w_{l_{c}}(X^{a},\tau)]\right|w'_{l_{c}}(X^{a},\tau)d\tau\nonumber \\
\sum_{i\in A'(1-\alpha_{2})}\int_{a(i)}^{a(i+1)}\left|z[w_{l_{a'(i+1)}}(X^{a},\tau)]\right|w'_{l_{a'(i+1)}}(X^{a},\tau)d\tau\label{eq:4-43}
\end{align}
Due to the hole $\mathbb{C}_{l}$ created in the bundle $B_{\mathbb{R}}(\mathbb{C})$,
the remaining length of the contour that will be subjected to removal
process is obtained by removing the sum of piecewise contour lengths
in $\mathbb{C}_{l},$ and is given by
\begin{equation}
\begin{array}{c}
L\left(\gamma_{l}(X^{a},t,1-\alpha_{2})\left(t\in[t_{0},t_{c}]\right)\right)=\\
\sum_{i\in A'(1-\alpha_{2})}\int_{a'(i)}^{a'(i+1)}\left|z[w_{l_{a'(i+1)}}(X^{a},\tau)]\right|w'_{l_{a'(i+1)}}(X^{a},\tau)d\tau+\\
\int_{a'_{c-1}}^{a'_{c}}\left|z[w_{l_{c}}(X^{a},\tau)]\right|w'_{l_{c}}(X^{a},\tau)d\tau\\
\sum_{i\in A'(1-\alpha_{2})}\int_{a(i)}^{a(i+1)}\left|z[w_{l_{a'(i+1)}}(X^{a},\tau)]\right|w'_{l_{a'(i+1)}}(X^{a},\tau)d\tau
\end{array}\label{eq:4-44}
\end{equation}

Since $\gamma_{l}(X^{a},t)$ is active in $B_{\mathbb{R}}(\mathbb{C},1-\alpha_{2})$,
the formation of the contour will continue forever. the tail part,
that is the sum of the pieces of the lengths of $\gamma_{l}(X^{a},t)$
that is there in $B_{\mathbb{R}}(\mathbb{C},\alpha_{2})$ will be
deleted from $B_{\mathbb{R}}(\mathbb{C},\alpha_{2}).$ This deletion
could be according to a removal function $\phi(X^{a},t,\alpha_{2})$
similar to the procedure explained in (\ref{eq:4-37}). The differential
equation to model the space of points lost due to a removal process
$\phi(X^{a},t,\alpha_{2})$ is 
\[
\]
\begin{align}
\left.\frac{dB_{\mathbb{R}}(\mathbb{C},\alpha_{2})}{dt}\right|_{\phi(X^{a},t,\alpha_{2})} & =\left.B_{\mathbb{R}}(\mathbb{C},\alpha_{2})\right|_{t=t_{c}}-\phi(X^{a},t,\alpha_{2})\left.B_{\mathbb{R}}(\mathbb{C},\alpha_{2})\right|_{t=t_{c}}\nonumber \\
\label{eq:4-45}
\end{align}
where
\begin{align}
\phi(X^{a},t,\alpha_{2})=\psi(X^{a},t,\alpha_{2})\times\nonumber \\
\sum_{i\in A'(1-\alpha_{2})}\int_{a'(i)}^{a'(i+1)}\left|z[w_{l_{a'(i+1)}}(X^{a},\tau)]\right|w'_{l_{a'(i+1)}}(X^{a},\tau)d\tau\text{ \ensuremath{\left(t\in(t_{c},t_{c'}]\right)}}\nonumber \\
\label{eq:4-46}
\end{align}
for $0<\psi(X^{a},t,\alpha_{2})<1.$

\begin{figure}
\textbf{\includegraphics[scale=0.5]{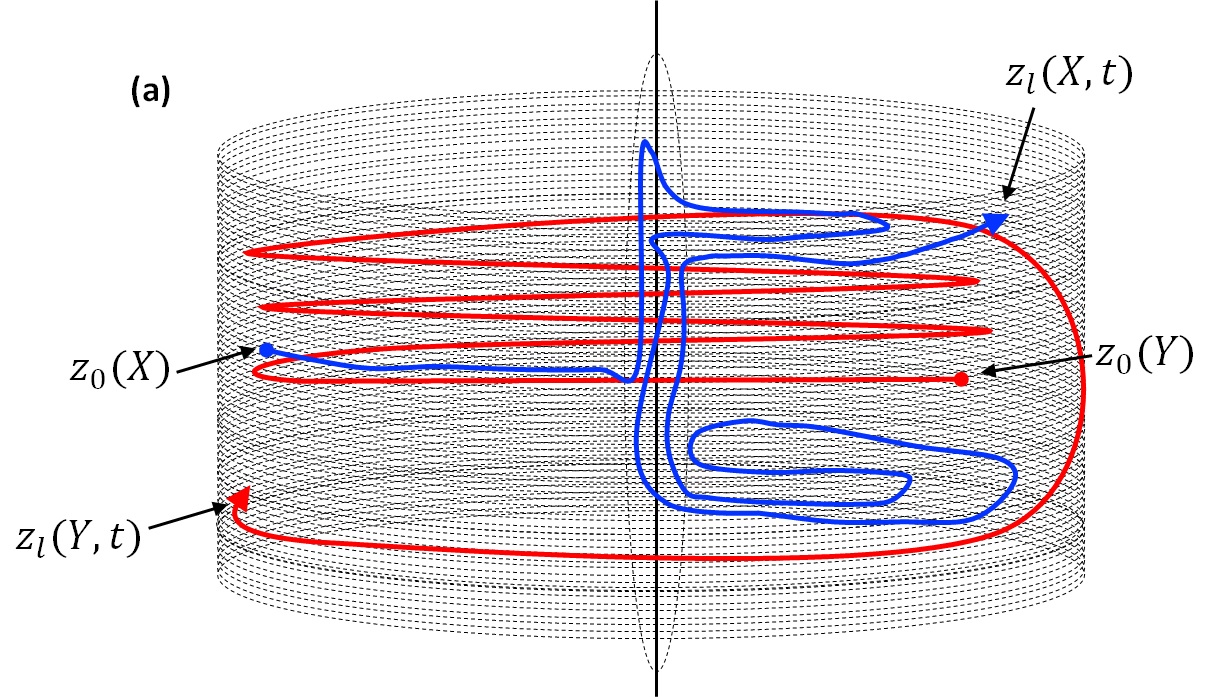}}

$ $

\textbf{\includegraphics[scale=0.5]{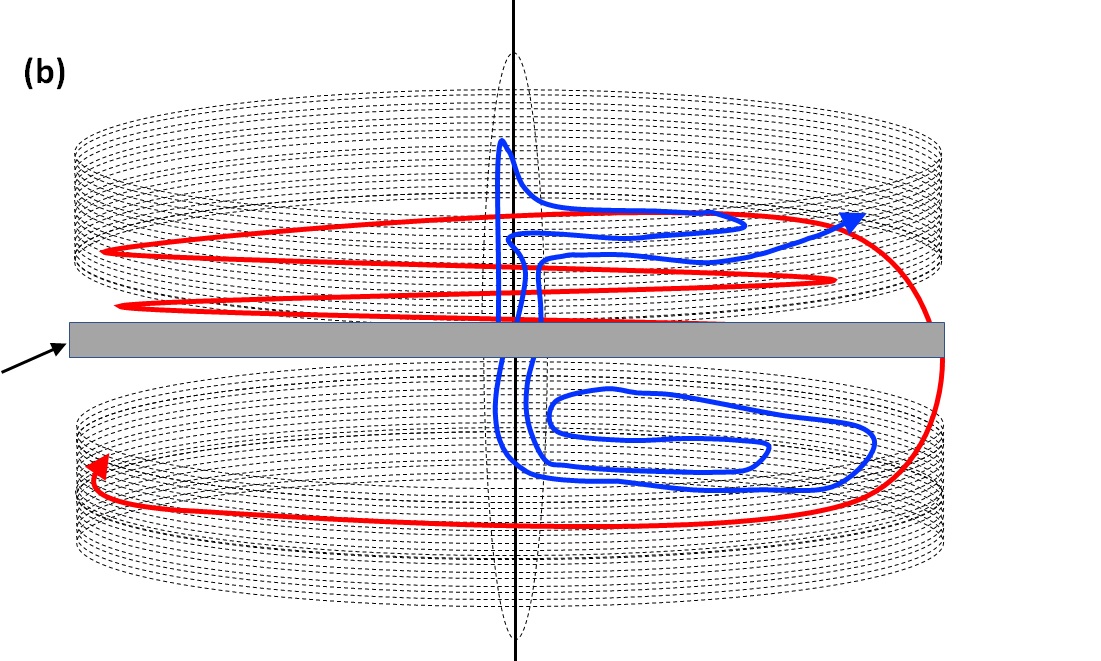}}

\caption{\label{fig:Formation-of-holeinBundle}\textbf{(a)} Formation of contours
in a bundle, and \textbf{(b)} formation of a hole in the bundle due
to a removal process of an arbitrary plane in which contours originated. }

\end{figure}

\subsection{PDE's for the Dynamics of Lost Space}

Suppose we consider simultaneously another arbitrary random variable
$X^{b}\left(z_{l}(t),\mathbb{C}_{l}\right)$ along with $X^{a}\left(z_{l}(t),\mathbb{C}_{l}\right)$
in $\mathbb{C}_{l}$. As in $X^{b}\left(z_{l}(t),\mathbb{C}_{l}\right)$,
we will develop two models, one for understanding the dynamics of
the space lost by first assuming that $\gamma_{l}(X^{b},t)$ at $t=t_{c}$
is active in $B_{\mathbb{R}}(\mathbb{C},\alpha_{2})$ and then second
time by assuming $\gamma_{l}(X^{b},t)$ at $t=t_{c}$ is active in
$B_{\mathbb{R}}(\mathbb{C},1-\alpha_{2})$. When $\gamma_{l}(X^{b},t)$
at $t=t_{c}$ is active in $B_{\mathbb{R}}(\mathbb{C},\alpha_{2})$,
the set of points touched by the contour $\gamma_{l}(X^{b},t)$ prior
to $t_{c}$ would have located in $B_{\mathbb{R}}(\mathbb{C},\alpha_{2})$,
$\mathbb{C}_{l}$, and $B_{\mathbb{R}}(\mathbb{C},1-\alpha_{2}).$
The contour $\gamma_{l}(X^{b},t)$ is described by $z_{l}(X^{b},t)$
$\left(t\in[t_{0},\infty)\right)$ and $t=O_{l_{c}}(X^{b},\tau,\alpha_{2})$
$\left(b_{0}\leq\tau\leq b{}_{c}\right)$ is the parametric representation
for $\gamma_{l}(X^{b},t,\alpha_{2})$ with a real values function
$O_{l_{c}}(X^{b},\tau,\alpha_{2})$ mapping $[b_{0},b_{c}]$ onto
the interval $[t_{0},t_{c}].$ Let $L\left(\gamma_{l}(X^{b},t,\alpha_{2})\left(t\in[t_{0},t_{c}]\right)\right)$
represent the length of $\gamma_{l}(X^{b},t,\alpha_{2})\left(t\in[t_{0},t_{c}]\right)$
up to $t_{c},$ then

\begin{align}
L\left(\gamma_{l}(X^{b},t,\alpha_{2})\left(t\in[t_{0},t_{c}]\right)\right)=\int_{b_{0}}^{b_{c}}\left|z[O_{l_{c}}(X^{b},\tau,\alpha_{2})]\right|O'_{l_{c}}(X^{b},\tau,\alpha_{2})d\tau.\nonumber \\
\label{eq:4-47}
\end{align}
By following the procedure explained in section 4.1, we partition
\[
L\left(\gamma_{l}(X^{b},t,\alpha_{2})\left(t\in[t_{0},t_{c}]\right)\right)
\]
into a sum of three parts of piecewise contour lengths. The final
differential equation to model the space of points lost due to only
the removal process $\phi(X^{b},t,1-\alpha_{2})$ is 
\begin{align}
\left.\frac{dB_{\mathbb{R}}(\mathbb{C},1-\alpha_{2})}{dt}\right|_{\phi(X^{b},t,1-\alpha_{2})} & =\left.B_{\mathbb{R}}(\mathbb{C},1-\alpha_{2})\right|_{t=t_{c}}-\nonumber \\
 & =\phi(X^{b},t,1-\alpha_{2})\left.B_{\mathbb{R}}(\mathbb{C},1-\alpha_{2})\right|_{t=t_{c}}\label{eq:4-48}
\end{align}
where
\begin{align}
\phi(X^{b},t,1-\alpha_{2})=\psi(X^{b},t,1-\alpha_{2})\times\nonumber \\
\sum_{i\in B(1-\alpha_{2})}\int_{b(i)}^{b(i+1)}\left|z[O_{l_{b(i+1)}}(X^{b},\tau)]\right|O'_{l_{b(i+1)}}(X^{b},\tau)d\tau\text{ \ensuremath{\left(t\in(t_{c},t_{c'}]\right)}}.\nonumber \\
\label{eq:4-50}
\end{align}
 for $0<\psi(X^{b},t,1-\alpha_{2})<1.$ Meaning of the set $B(1-\alpha_{2})$
and the procedure to obtain the integral in (\ref{eq:4-50}) are similar
to corresponding model in (\ref{eq:4-38}). Alternatively, when $\gamma_{l}(X^{b},t)$
at $t=t_{c}$ is active in $B_{\mathbb{R}}(\mathbb{C},1-\alpha_{2})$,
the set of points touched by the contour $\gamma_{l}(X^{b},t)$ prior
to $t_{c}$ would have located in $B_{\mathbb{R}}(\mathbb{C},\alpha_{2})$,
$\mathbb{C}_{l}$, and $B_{\mathbb{R}}(\mathbb{C},1-\alpha_{2}).$
The the contour $\gamma_{l}(X^{b},t)$ can be described by $z_{l}(X^{b},t)$
$\left(t\in[t_{0},\infty)\right)$ and $t=Q_{l_{c}}(X^{b},\tau,1-\alpha_{2})$
$\left(b'_{0}\leq\tau\leq b'{}_{c}\right)$ is the parametric representation
for $\gamma_{l}(X^{b},t,1-\alpha_{2})$ with a real values function
$Q_{l_{c}}(X^{b},\tau,1-\alpha_{2})$ mapping $[b'_{0},b'_{c}]$ onto
the interval $[t_{0},t_{c}].$ Let $L\left(\gamma_{l}(X^{b},t,1-\alpha_{2})\left(t\in[t_{0},t_{c}]\right)\right)$
represent the length of $\gamma_{l}(X^{b},t,1-\alpha_{2})\left(t\in[t_{0},t_{c}]\right)$
up to $t_{c},$ then
\begin{equation}
\begin{array}{c}
L\left(\gamma_{l}(X^{b},t,1-\alpha_{2})\left(t\in[t_{0},t_{c}]\right)\right)=\\
\int_{b'_{0}}^{b'_{c}}\left|z[Q_{l_{c}}(X^{b},\tau,1-\alpha_{2})]\right|Q'_{l_{c}}(X^{b},\tau,1-\alpha_{2})d\tau.
\end{array}\label{eq:4-51}
\end{equation}
By following the procedure explained in section 4.1, we partition
\[
L\left(\gamma_{l}(X^{b},t,1-\alpha_{2})\left(t\in[t_{0},t_{c}]\right)\right)
\]
into a sum of three parts of piecewise contour lengths. The corresponding
differential equation to model the space of points lost due to only
the removal process $\phi(X^{b},t,\alpha_{2})$ is 
\begin{align}
\left.\frac{dB_{\mathbb{R}}(\mathbb{C},\alpha_{2})}{dt}\right|_{\phi(X^{b},t,\alpha_{2})} & =\left.B_{\mathbb{R}}(\mathbb{C},1-\alpha_{2})\right|_{t=t_{c}}-\nonumber \\
 & \phi(X^{b},t,\alpha_{2})\left.B_{\mathbb{R}}(\mathbb{C},1-\alpha_{2})\right|_{t=t_{c}}\label{eq:4-52}
\end{align}
where
\begin{align}
\phi(X^{b},t,\alpha_{2})=\psi(X^{b},t,\alpha_{2})\times\nonumber \\
\sum_{i\in B'(1-\alpha_{2})}\int_{b'(i)}^{b'(i+1)}\left|z[Q_{l_{b(i+1)}}(X^{b},\tau)]\right|Q'_{l_{b(i+1)}}(X^{b},\tau)d\tau\text{ \ensuremath{\left(t\in(t_{c},t_{c'}]\right)}}.\nonumber \\
\label{eq:4-53}
\end{align}
for $0<\psi(X^{b},t,1-\alpha_{2})<1.$ The set $B'(1-\alpha_{2})$
consists of piecewise contours in $B_{\mathbb{R}}(\mathbb{C},1-\alpha_{2})$
as described in section 4.1 and the procedure to obtain the integral
in (\ref{eq:4-53}) are similar to corresponding model in (\ref{eq:4-46}).
The differential equations (\ref{eq:4-48}) and (\ref{eq:4-52}) are
good when the removal process of $X^{a}(z_{l}(t),t)$ was not introduced
simultaneously. The differential equations (\ref{eq:4-37}) and (\ref{eq:4-45})
does not provide dynamics of lost spaces when the removal process
of $X^{b}(z_{l}(t),t)$ was not introduced simultaneously. Under the
simultaneous existence of the two contours $\gamma_{l}(X^{a},t)$
and $\gamma_{l}(X^{b},t)$, and the corresponding removal processes,
we will have four situations. Suppose $\gamma_{l}(X^{a},t)$ and $\gamma_{l}(X^{b},t)$
are active in $B_{\mathbb{R}}(\mathbb{C},\alpha_{2})$ at $t=t_{c}$,
then the partial differential equation describing the dynamics of
removal of the spaces created by $\gamma_{l}(X^{a},t)$ and $\gamma_{l}(X^{b},t)$
in $B_{\mathbb{R}}(\mathbb{C},1-\alpha_{2})$ until $t_{c}$ are
\[
\]
\begin{align}
\left.\frac{\partial^{2}B_{\mathbb{R}}(\mathbb{C},1-\alpha_{2})}{\partial X^{a}\partial t}\right|_{\phi(X^{a},t,1-\alpha_{2})} & =\left.B_{\mathbb{R}}(\mathbb{C},1-\alpha_{2})\right|_{t=t_{c}}-\nonumber \\
 & \frac{\partial\phi(X^{a},t,1-\alpha_{2})\left[L(X^{a}+X^{b},1-\alpha_{2})\right]}{\partial X^{a}},\label{eq:4-54}
\end{align}
\begin{align}
\left.\frac{\partial^{2}B_{\mathbb{R}}(\mathbb{C},1-\alpha_{2})}{\partial X^{b}\partial t}\right|_{\phi(X^{b},t,1-\alpha_{2})} & =\left.B_{\mathbb{R}}(\mathbb{C},1-\alpha_{2})\right|_{t=t_{c}}-\nonumber \\
 & \frac{\partial\phi(X^{b},t,1-\alpha_{2})\left[L(X^{a}+X^{b},1-\alpha_{2})\right]}{\partial X^{b}},\label{eq:4-55}
\end{align}
 for the removal rates $\phi(X^{a},t,1-\alpha_{2})$, $\phi(X^{b},t,1-\alpha_{2})$
and length of contour formed due to $\gamma_{l}(X^{a},t)$ and $\gamma_{l}(X^{b},t)$
in $B_{\mathbb{R}}(\mathbb{C},1-\alpha_{2})$ is $L(X^{a}+X^{b},1-\alpha_{2}).$
Suppose $\gamma_{l}(X^{a},t)$ is active in $B_{\mathbb{R}}(\mathbb{C},\alpha_{2})$
and $\gamma_{l}(X^{b},t)$ is active in $B_{\mathbb{R}}(\mathbb{C},1-\alpha_{2})$
at $t=t_{c}$, then the partial differential equations describing
the dynamics of removal of the space created by $\gamma_{l}(X^{a},t)$
in $B_{\mathbb{R}}(\mathbb{C},1-\alpha_{2})$ and $\gamma_{l}(X^{b},t)$
in $B_{\mathbb{R}}(\mathbb{C},\alpha_{2})$ until $t_{c}$ are 
\begin{align}
\left.\frac{\partial^{2}B_{\mathbb{R}}(\mathbb{C},1-\alpha_{2})}{\partial X^{a}\partial t}\right|_{\phi(X^{a},t,1-\alpha_{2})} & =\left.B_{\mathbb{R}}(\mathbb{C},1-\alpha_{2})\right|_{t=t_{c}}-\nonumber \\
 & \frac{\partial\phi(X^{a},t,1-\alpha_{2})\left[L(X^{a},1-\alpha_{2})+L(X^{b},\alpha_{2})\right]}{\partial X^{a}},\label{eq:4-56}
\end{align}
\begin{align}
\left.\frac{\partial^{2}B_{\mathbb{R}}(\mathbb{C},\alpha_{2})}{\partial X^{b}\partial t}\right|_{\phi(X^{b},t,1-\alpha_{2})} & =\left.B_{\mathbb{R}}(\mathbb{C},\alpha_{2})\right|_{t=t_{c}}-\nonumber \\
 & \frac{\partial\phi(X^{b},t,\alpha_{2})\left[L(X^{a},1-\alpha_{2})+L(X^{b},\alpha_{2})\right]}{\partial X^{a}},\label{eq:4-57}
\end{align}
where $\phi(X^{a},t,1-\alpha_{2})$ and $\phi(X^{b},t,\alpha_{2})$
are the removal rates of $\gamma_{l}(X^{a},t)$ and $\gamma_{l}(X^{b},t)$
in $B_{\mathbb{R}}(\mathbb{C},1-\alpha_{2})$, and $B_{\mathbb{R}}(\mathbb{C},\alpha_{2})$,
respectively. The sum of the piecewise lengths in $B_{\mathbb{R}}(\mathbb{C},1-\alpha_{2})$,
and $B_{\mathbb{R}}(\mathbb{C},\alpha_{2})$ are represented by $L(X^{a},1-\alpha_{2})$
and $L(X^{b},\alpha_{2}),$ respectively. Suppose $\gamma_{l}(X^{a},t)$
and $\gamma_{l}(X^{b},t)$ are active in $B_{\mathbb{R}}(\mathbb{C},1-\alpha_{2})$
at $t=t_{c}$, then the partial differential equation describing the
dynamics of removal of the spaces created by $\gamma_{l}(X^{a},t)$
and $\gamma_{l}(X^{b},t)$ in $B_{\mathbb{R}}(\mathbb{C},\alpha_{2})$
until $t_{c}$ are 
\begin{align}
\left.\frac{\partial^{2}B_{\mathbb{R}}(\mathbb{C},\alpha_{2})}{\partial X^{a}\partial t}\right|_{\phi(X^{a},t,\alpha_{2})} & =\left.B_{\mathbb{R}}(\mathbb{C},\alpha_{2})\right|_{t=t_{c}}-\nonumber \\
 & \frac{\partial\phi(X^{a},t,\alpha_{2})\left[L(X^{a}+X^{b},\alpha_{2})\right]}{\partial X^{a}},\label{eq:4-58}
\end{align}
\begin{align}
\left.\frac{\partial^{2}B_{\mathbb{R}}(\mathbb{C},\alpha_{2})}{\partial X^{b}\partial t}\right|_{\phi(X^{b},t,\alpha_{2})} & =\left.B_{\mathbb{R}}(\mathbb{C},\alpha_{2})\right|_{t=t_{c}}-\nonumber \\
 & \frac{\partial\phi(X^{b},t,\alpha_{2})\left[L(X^{a}+X^{b},\alpha_{2})\right]}{\partial X^{b}},\label{eq:4-59}
\end{align}
where $\phi(X^{a},t,\alpha_{2})$ and $\phi(X^{b},t,\alpha_{2})$
are the removal rates and the sums of the piecewise lengths of contours
formed due to $\gamma_{l}(X^{a},t)$ and $\gamma_{l}(X^{b},t)$ in
$B_{\mathbb{R}}(\mathbb{C},\alpha_{2})$ is $L(X^{a}+X^{b},\alpha_{2}).$
Suppose $\gamma_{l}(X^{a},t)$ is active in $B_{\mathbb{R}}(\mathbb{C},1-\alpha_{2})$
and $\gamma_{l}(X^{b},t)$ is active in $B_{\mathbb{R}}(\mathbb{C},\alpha_{2})$
at $t=t_{c}$, then the partial differential equations describing
the dynamics of removal of the space created by $\gamma_{l}(X^{a},t)$
in $B_{\mathbb{R}}(\mathbb{C},\alpha_{2})$ and $\gamma_{l}(X^{b},t)$
in $B_{\mathbb{R}}(\mathbb{C},1-\alpha_{2})$ until $t_{c}$ are 
\begin{align}
\left.\frac{\partial^{2}B_{\mathbb{R}}(\mathbb{C},\alpha_{2})}{\partial X^{a}\partial t}\right|_{\phi(X^{a},t,\alpha_{2})} & =\left.B_{\mathbb{R}}(\mathbb{C},\alpha_{2})\right|_{t=t_{c}}-\nonumber \\
 & \frac{\partial\phi(X^{a},t,\alpha_{2})\left[L(X^{a},\alpha_{2})+L(X^{b},1-\alpha_{2})\right]}{\partial X^{a}},\label{eq:4-60}
\end{align}
\begin{align}
\left.\frac{\partial^{2}B_{\mathbb{R}}(\mathbb{C},1-\alpha_{2})}{\partial X^{b}\partial t}\right|_{\phi(X^{b},t,1-\alpha_{2})} & =\left.B_{\mathbb{R}}(\mathbb{C},1-\alpha_{2})\right|_{t=t_{c}}-\nonumber \\
 & \frac{\partial\phi(X^{b},t,\alpha_{2})\left[L(X^{a},\alpha_{2})+L(X^{b},1-\alpha_{2})\right]}{\partial X^{a}},\label{eq:4-61}
\end{align}
where $\phi(X^{a},t,\alpha_{2})$ and $\phi(X^{b},t,1-\alpha_{2})$
are the removal rates of $\gamma_{l}(X^{a},t)$ and $\gamma_{l}(X^{b},t)$
in $B_{\mathbb{R}}(\mathbb{C},\alpha_{2})$, and $B_{\mathbb{R}}(\mathbb{C},1-\alpha_{2})$,
respectively. The sum of the piecewise lengths in $B_{\mathbb{R}}(\mathbb{C},\alpha_{2})$,
and $B_{\mathbb{R}}(\mathbb{C},1-\alpha_{2})$ are represented by
$L(X^{a},\alpha_{2})$ and $L(X^{b},1-\alpha_{2}),$ respectively.
The two partitions in $B_{\mathbb{R}}(\mathbb{C},\alpha_{2})$, and
$B_{\mathbb{R}}(\mathbb{C},1-\alpha_{2})$ behave independently in
terms of further removal and growth of contours. 

The PDEs represented in (\ref{eq:4-54}) through (\ref{eq:4-61})
are simultaneous removal of two contours at a time. These can be treated
as an example of the removal process of the dynamics of multiple removal
processes. The plan is to remove all the infinitely many contours
at $t_{c}$. 

\section{\textbf{Concluding Remarks}}

Multilevel contours passing through bundle $B_{\mathbb{R}}(\mathbb{C})$
of complex planes could demonstrate interesting properties. The random
environment created brings the dynamic nature of the bundle through
the removal process introduced. The continuous-time Markov properties,
differential equations, and topological analysis on the bundle give
scope to further investigate them using functional approximations.
The transportation of information through contours for different complex
planes could be extended further for practical situations arising
out of transportation problems. There are several applications of
complex analysis that are out of scope to discuss in this article.
A wide range of literature is available for interested readers, see
for example \cite{Rao-Krantz-handbook,editedbookMexico,Ponnuswamy,Campos,pathak,Cohen,Krantz-applicationsbook}.
One can also introduce several forms of parametric contour formations
by assuming functional growth rates of contours. That could be an
independent approach towards modeling the behavior of the contours
concerning a given functional form of contour formation. Similarly,
the removal rates can be assumed to follow certain closed-form approximations
(special forms of Harmonic functions, and Poisson integrals). Information
carried between various complex planes that were discussed in Section
2 would be obstructed if the set of points that a given connected
contour gets deleted (lost) due to a removal process. 

\section*{\textbf{Acknowledgments}}

I wish to thank and appreciate our children (daughter: Sheetal Rao,
son: GopalKrishna Rao, son: Raghav Rao) whose several weekends play
time with me was sacrificed by them while I was occupied with this
project during the summer of 2021. 

\pagebreak


\begin{thebibliography}{10}
\bibitem[1]{Ahlfors-book1978}Ahlfors, Lars V. Complex analysis. An
introduction to the theory of analytic functions of one complex variable.
Third edition. International Series in Pure and Applied Mathematics.
McGraw-Hill Book Co., New York, 1978. xi+331

\bibitem[2]{Churchil-Brwon-book}Churchill, Ruel V.; Brown, James
Ward Complex variables and applications. Fourth edition. McGraw-Hill
Book Co., New York, 1984. x+339 pp.

\bibitem[3]{Krantz-geometric}Krantz, Steven G. Complex analysis:
the geometric viewpoint. Second edition. Carus Mathematical Monographs,
23. Mathematical Association of America, Washington, DC, 2004. xviii+219
pp.

\bibitem[4]{Rudin-Real-Complex}Rudin, Walter Real and complex analysis.
Third edition. McGraw-Hill Book Co., New York, 1987. xiv+416 pp. ISBN:
0-07-054234-1

\bibitem[5]{Rao-Krantz-handbook}Rao A.S.R.S. and Krantz S.G. (2021).
Rao distances and conformal mappings,\emph{ Information Geometry,}
\emph{Volume 45: Handbook of Statistics}, Elsevier/North-Holland,
Amsterdam

\bibitem[6]{good}Good, I. J. The frequency count of a Markov chain
and the transition to continuous time. Ann. Math. Statist. 32 (1961),
41--48.

\bibitem[7]{bRBhat}Bhat, B. R.; Deshpande, Sunita K. Likelihood ratio
test for testing order of continuous time finite Markov chains. Comm.
Statist. A---Theory Methods 15 (1986), no. 6, 1751--1771. 

\bibitem[8]{ChenMa}Chen, Mu Fa On three classical problems for Markov
chains with continuous time parameters. J. Appl. Probab. 28 (1991),
no. 2, 305--320.

\bibitem[9]{Gani}Gani, J.; Stals, L. A continuous time Markov chain
model for a plantation-nursery system. Environmetrics 16 (2005), no.
8, 849--861.

\bibitem[10]{GoswamiRao}Goswami, A.; Rao, B. V. A course in applied
stochastic processes. Texts and Readings in Mathematics, 40. Hindustan
Book Agency, New Delhi, 2006.

\bibitem[11]{editedbookMexico} Sagun Chanillo, Paulo D. Cordaro,
Nicholas Hanges, Jorge Hounie and Abdelhamid Meziani (Edited). Geometric
analysis of PDE and several complex variables. Dedicated to Fran�ois
Treves. Including papers from the workshop held in S�o Paulo, August
2003.. Contemporary Mathematics, 368. American Mathematical Society,
Providence, RI, 2005

\bibitem[12]{Ponnuswamy}Ponnusamy, S.; Silverman, Herb Complex variables
with applications. Birkh�user Boston, Inc., Boston, MA, 2006.

\bibitem[13]{Campos}Campos, L. M. B. C. Complex analysis with applications
to flows and fields. Mathematics and Physics for Science and Technology.
CRC Press, Boca Raton, FL, 2011.

\bibitem[14]{pathak}Pathak, H.K. Complex analysis and applications.
Springer, Singapore, 2019

\bibitem[15]{Cohen}Cohen, Harold Complex analysis with applications
in science and engineering. Second edition. Springer, New York, 2007.

\bibitem[16]{Krantz-applicationsbook}Krantz, Steven G. Complex variables.
A physical approach with applications and MATLAB�. Textbooks in Mathematics.
Chapman \& Hall/CRC, Boca Raton, FL, 2008.
\end{thebibliography}
\end{document}